\documentclass{amsart}

\usepackage[leqno]{amsmath}
\usepackage{latexsym,amsthm,amsfonts,amssymb,hyperref}
\usepackage[all]{xy} \SelectTips{eu}{}

\newcommand{\numberseries}{\bfseries}   

\newlength{\thmtopspace}                
\newlength{\thmbotspace}                
\newlength{\thmheadspace}               
\newlength{\thmindent}                  

\renewcommand{\subparagraph}{\vspace*{\thmbotspace}}

\newtheoremstyle{bfupright head,slanted body}
                {\thmtopspace}{\thmbotspace}
                {\slshape}{\thmindent}{\bfseries}{.}{\thmheadspace}
                {{\numberseries \thmnumber{#2\;}}\thmnote{#3}}

\newtheoremstyle{bfupright head,upright body}
                {\thmtopspace}{\thmbotspace}
                {\upshape}{\thmindent}{\bfseries}{.}{\thmheadspace}
                {{\numberseries \thmnumber{#2\;}}\thmnote{#3}}

\newtheoremstyle{bfit head,upright body}
                {\thmtopspace}{\thmbotspace}
                {\upshape}{\thmindent}{\upshape}{.}{\thmheadspace}
                {{\numberseries\thmnumber{#2\;}}
                {\bfseries\itshape\thmnote{\negthickspace#3}}}

\newtheoremstyle{it head,upright body}
                {\thmtopspace}{\thmbotspace}
                {\upshape}{\thmindent}{\upshape}{.}{\thmheadspace}
                {{\numberseries\thmnumber{#2\;}}
                {\itshape\thmnote{\negthickspace#3}}}

\newtheoremstyle{fixed bf head,slanted body}
                {\thmtopspace}{\thmbotspace}{\slshape}
                {\thmindent}{\bfseries}{.}{\thmheadspace}
                {{\numberseries \thmnumber{#2\;}}\thmname{#1}\thmnote{ (#3)}}

\newtheoremstyle{fixed bf head,upright body}
                {\thmtopspace}{\thmbotspace}{\upshape}
                {\thmindent}{\bfseries}{.}{\thmheadspace}
                {{\numberseries \thmnumber{#2\;}}\thmname{#1}\thmnote{ (#3)}}

\newtheoremstyle{fixed bfit head,upright body}
                {\thmtopspace}{\thmbotspace}{\upshape}
                {\thmindent}{\bfseries\itshape}{.}{\thmheadspace}
                {{\numberseries \thmnumber{#2\;}}\thmname{#1}\thmnote{ (#3)}}

\newtheoremstyle{sc head,small body}
                {\thmtopspace}{\thmbotspace}
                {\small\upshape}{\thmindent}{\scshape}{.}{\thmheadspace}
                {\thmname{#1}}

\newtheoremstyle{numbered paragraph}
                {\thmtopspace}{\thmbotspace}{\upshape}
                {\thmindent}{\upshape}{}{\thmheadspace}
                {{\numberseries \thmnumber{#2.}}}

\newtheoremstyle{unnumbered paragraph}
                {\thmtopspace}{\thmbotspace}{\upshape}
                {\parindent}{\upshape}{}{0pt}

\theoremstyle{bfupright head,slanted body}
\newtheorem{res}{}[section]             \newtheorem*{res*}{}

\theoremstyle{bfit head,upright body}
                 \newtheorem*{com*}{}

\theoremstyle{bfupright head,upright body}
\newtheorem{bfhpg}[res]{}               \newtheorem*{bfhpg*}{}

\theoremstyle{it head,upright body}
               \newtheorem*{ithpg*}{}

\theoremstyle{sc head,small body}

\theoremstyle{fixed bf head,slanted body}
\newtheorem{thm}[res]{Theorem}          \newtheorem*{thm*}{Theorem}
\newtheorem{prp}[res]{Proposition}      \newtheorem*{prp*}{Proposition}
        \newtheorem*{cor*}{Corollary}
\newtheorem{lem}[res]{Lemma}            \newtheorem*{lem*}{Lemma}

\theoremstyle{fixed bf head,upright body}
       \newtheorem*{dfn*}{Definition}
     \newtheorem*{con*}{Construction}
      \newtheorem*{obs*}{Observation}
\newtheorem{rmk}[res]{Remark}           \newtheorem*{rmk*}{Remark}
\newtheorem{exa}[res]{Example}          \newtheorem*{exa*}{Example}
         \newtheorem*{exe*}{Exercise}
\newtheorem{stp}[res]{Setup}            \newtheorem*{stp*}{Setup}
             \newtheorem*{fct*}{Fact}
       \newtheorem*{cnj*}{Conjecture}

\theoremstyle{numbered paragraph}

\theoremstyle{unnumbered paragraph}
\newtheorem{ipg*}{}

\newlength{\thmlistleft}        
\newlength{\thmlistright}       
\newlength{\thmlistpartopsep}   
\newlength{\thmlisttopsep}      
\newlength{\thmlistparsep}      
\newlength{\thmlistitemsep}     

\newcounter{eqc} 
  {\end{list}}%

\newcounter{prt}
  {\end{list}}%

\newcounter{rqm}
  {\end{list}}%

  {\end{list}}%

  {\end{list}}%

\newcounter{exercise}
  {\end{list}}%

\newenvironment{prf*}[1][Proof]{%
  \begin{proof}[\bf #1]
    \setcounter{equation}{0}
    \renewcommand{\theequation}{\arabic{equation}}}
  {\end{proof}
}

\newcommand{\pgref}[1]{{\rm \ref{#1}}}

\newcommand{\stpref}[2][Setup~]{#1\ref{stp:#2}}

\newcommand{\thmref}[2][Theorem~]{#1\pgref{thm:#2}}

\newcommand{\prpref}[2][Proposition~]{#1\pgref{prp:#2}}
\newcommand{\lemref}[2][Lemma~]{#1\pgref{lem:#2}}

\newcommand{\rmkref}[2][Remark~]{#1\pgref{rmk:#2}}

\newcommand{\secref}[2][Section~]{#1\ref{sec:#2}}

\renewcommand{\eqref}[1]{(\pgref{eq:#1})}

\newcommand{\thmcite}[2][?]{\cite[Theorem~#1]{#2}}

\newcommand{\corcite}[2][?]{\cite[Corollary~#1]{#2}}
\newcommand{\prpcite}[2][?]{\cite[Proposition~#1]{#2}}
\newcommand{\lemcite}[2][?]{\cite[Lemma~#1]{#2}}

\newcommand{\chpcite}[2][?]{\cite[Chapter~#1]{#2}}
\newcommand{\seccite}[2][?]{\cite[Section~#1]{#2}}

\newcommand{\Sym}[2][z]{\ensuremath{#2}%
\index{symbols}{#1@\ensuremath{\protect#2}}}

\def\urltilda{\kern -.15em\lower .7ex\hbox{\~{}}\kern .04em} 
\renewcommand{\emptyset}{\varnothing}
\newcommand{\kd}[1]{\delta_{#1}}
\newcommand{\LHS}{\text{\small LHS}}
\newcommand{\RHS}{\text{\small RHS}}

\newcommand{\set}[2][\mspace{1mu}]{\{#1 #2 #1\}}
\newcommand{\setof}[3][\mspace{1mu}]{\{#1#2 \mid #3#1\}}

\newcommand{\ZZ}{\mathbb{Z}}
\newcommand{\PP}{\mathbb{P}}
\newcommand{\WW}{\mathbb{W}}

\newcommand{\qtext}[1]{\quad\text{#1}\quad}
\newcommand{\qqtext}[1]{\qquad\text{#1}\qquad}
\newcommand{\qand}{\qtext{and}}
\newcommand{\qqand}{\qqtext{and}}
\newcommand{\deq}{\:=\:}
\renewcommand{\le}{\leqslant}
\renewcommand{\ge}{\geqslant}
\newcommand{\dis}{\:\is\:}
\newcommand{\dle}{\:\le\:}

\newcommand{\mfA}{\mathfrak{A}}
\newcommand{\mfM}{\mathfrak{M}}
\newcommand{\mfa}{\mathfrak{a}}
\newcommand{\mfm}{\mathfrak{m}}

\newcommand{\is}{\cong}\newcommand{\lra}{\longrightarrow}
\newcommand{\xra}[2][]{\xrightarrow[#1]{\;#2\;}}

\newcommand{\tp}[3][R]{\nobreak{#2\otimes_{#1}#3}}

\newcommand{\Tor}[4][R]{\operatorname{Tor}^{#1}_{#2}(#3,#4)}

\newcommand{\grd}[2]{\operatorname{grade}_{#1}(#2)}

 \newcommand{\pfT}[1][\calT]{\operatorname{pf}_#1}

\newcommand{\fa}{\mathfrak{a}}

\newcommand{\gra}{\alpha}
\newcommand{\grb}{\beta}
\newcommand{\grg}{\gamma}
\newcommand{\grt}{\tau}
\newcommand{\grr}{\rho}
\newcommand{\gro}{\omega}
\newcommand{\grs}{\sigma}
\newcommand{\calA}{\mathcal{A}}
\newcommand{\calR}{\mathcal{R}}
\newcommand{\calB}{\mathcal{B}}
\newcommand{\calF}{\mathcal{F}}
\newcommand{\calI}{\mathcal{I}}

\newcommand{\calL}{\mathcal{L}}

\newcommand{\calS}{\mathcal{S}}
\newcommand{\calT}{\mathcal{T}}
\newcommand{\calU}{\mathcal{U}}

\newcommand{\SO}[1]{\operatorname{SO}(#1)}
\newcommand{\W}[1]{\operatorname{W}(#1)}
\newcommand{\Spin}[1]{\operatorname{Spin}(#1)}
\renewcommand{\Sym}[2][k]{\operatorname{Sym}_{#1}(#2)}

\newcommand{\OGr}[1]{\operatorname{OGr}(#1)}
\newcommand{\ESUB}[1]{\mathcal{PE}_{#1}}

\numberwithin{equation}{res}

\newcommand{\pff}[1]{\operatorname{\mathcal{P}}[#1]}
\newcommand{\sgn}[2]{\operatorname{sgn}%
  \big(\begin{smallmatrix}#1\\#2 \end{smallmatrix}\big)}
\newcommand{\Sign}[1]{(-1)^{#1}}

\def\widebardisplay#1{%
  \setbox0=\hbox{$\displaystyle #1$}
  \dimen0=\wd0%
  \advance\dimen0 by -2pt
  \vbox{%
    \nointerlineskip%
    \moveright 1pt 
    \vbox{\hrule width \dimen0}%
    \nointerlineskip%
    \kern 1pt
    \box0%
    }%
  }

\def\widebartext#1{%
  \setbox0=\hbox{$#1$}
  \dimen0=\wd0%
  \advance\dimen0 by -2pt
  \vbox{%
    \nointerlineskip%
    \moveright 1pt 
    \vbox{\hrule width \dimen0}%
    \nointerlineskip%
    \kern 1pt
    \box0%
    }%
  }

\def\widebarscript#1{%
  \setbox0=\hbox{$\scriptstyle #1$}
  \dimen0=\wd0%
  \advance\dimen0 by -3pt
  \vbox{%
    \nointerlineskip%
    \moveright 1.5pt 
    \vbox{\hrule width \dimen0}%
    \nointerlineskip%
    \kern .8pt
    \box0%
    }%
  }

\def\widebarscriptscript#1{%
  \setbox0=\hbox{$\scriptscriptstyle #1$}
  \dimen0=\wd0%
  \advance\dimen0 by -3pt
  \vbox{%
    \nointerlineskip%
    \moveright 1.5pt 
    \vbox{\hrule width \dimen0}%
    \nointerlineskip%
    \kern .6pt
    \box0%
    }%
  }

\def\widebar#1{\mathchoice%
  {\widebardisplay{#1}}%
  {\widebartext{#1}}%
  {\widebarscript{#1}}%
  {\widebarscriptscript{#1}}%
  }

\makeatletter

\newcommand{\Pf}[2][]{\operatorname{Pf}_{#1}(#2)}
\newcommand{\Pfbar}[2][U]{\operatorname{Pf}_{\overline{#2}}(#1)}
\newcommand{\pf}[2][T]{\operatorname{pf}_{#1}{(#2)}}
\newcommand{\pfbar}[2][\calT]{\operatorname{pf}_{#1}(\overline{#2})}

\begin{document}

\title[Almost complete intersection ideals of grade
  3]{Three takes on\\ almost complete intersection ideals of grade 3}

\author[L.\,W. Christensen]{Lars Winther Christensen}

\address{Texas Tech University, Lubbock, TX 79409, U.S.A.}

\email{lars.w.christensen@ttu.edu}

\urladdr{http://www.math.ttu.edu/\urltilda lchriste}

\author[O. Veliche]{Oana Veliche}

\address{Northeastern University, Boston, MA~02115, U.S.A.}

\email{o.veliche@northeastern.edu}

\urladdr{https://web.northeastern.edu/oveliche}

\author[J. Weyman]{Jerzy Weyman}

\address{Jagiellonian University, 30-348 Krak\'ow, Poland and\newline
  \phantom{x}\hspace{4pt} University of Connecticut, Storrs, CT~06269,
  U.S.A.}  \email{jerzy.weyman@uj.edu.pl}

\urladdr{http://www.math.uconn.edu/\urltilda weyman}

\thanks{L.W.C.\ was partly supported by Simons Foundation
  collaboration grant 428308 and J.W.\ was partly supported by NSF DMS
  grant 1802067.}

\date{\today}



\begin{abstract}
  We are interested in the structure of almost complete intersection
  ideals of grade 3.  We give three constructions of these ideals and
  their free resolutions: one from the commutative algebra point of
  view, an equivariant construction giving a nice canonical form, and
  finally an interpretation in terms of open sets in certain Schubert
  varieties.
\end{abstract}


\maketitle

\allowdisplaybreaks
\thispagestyle{empty}

\section*{Contents}
\contentsline {section}{\tocsection {}{1}{Almost complete intersections following Avramov and Brown}}{2}{section.1}%
\contentsline {section}{\tocsection {}{2}{Generic almost complete intersections}}{4}{section.2}%
\contentsline {section}{\tocsection {}{3}{The equivariant form of the format $(1,4,n,n-3)$}}{13}{section.3}%
\contentsline {section}{\tocsection {}{4}{Schubert varieties in orthogonal Grassmannians}}{17}{section.4}%
\contentsline {section}{\tocsection {Appendix}{A}{Pfaffian identities following Knuth}}{27}{appendix.A}%
\contentsline {section}{\tocsection {Appendix}{B}{Minors via Pfaffians following Brill}}{32}{appendix.B}%
\contentsline {section}{\tocsection {Appendix}{C}{Generic almost complete intersections: the proofs}}{35}{appendix.C}%
\contentsline {section}{\tocsection {Appendix}{}{Acknowledgments}}{46}{section*.45}%
\contentsline {section}{\tocsection {Appendix}{}{References}}{46}{section*.46}%


\section*{Introduction}
\label{sec:introduction}

\noindent
Let $R$ be a commutative noetherian local ring. A celebrated result of
Buchsbaum and Eisenbud \cite{DABDEs77} states that every Gorenstein
ideal in $R$ of grade $3$ is generated by the $2m\times 2m$ Pfaffians
of a $(2m+1)\times (2m+1)$ skew symmetric matrix. Later, Avramov
\cite{LLA81a} and Brown \cite{AEB87} independently proved a similar
result for almost complete intersections. Their proofs are based on
the fact that an almost complete intersection ideal is linked to a
Gorenstein ideal; this means that their descriptions of the
resolutions depend on certain choices, so they are not coordinate
free.

In this paper we take three approaches to almost complete intersection
ideals of grade 3. They involve different languages, so they can be
appreciated by different audiences.  However we show how these three
approaches intertwine and influence each other.

The first approach uses only commutative and linear algebra. The main
theorems about grade $3$ almost complete intersection ideals in the
local ring $R$ are stated in \secref{local}. They are proved in
\secref{generic} and Appendices \secref[]{knuth}--\secref[]{odd} by
specialization from the generic case. In the generic case we use the
Buchsbaum--Eisenbud Acyclicity Criterion and a computation with
Pfaffians inspired by the Buchsbaum--Eisenbud Structure Theorems, see
\rmkref{how}, to construct the minimal free resolutions. We emphasize
that our description of the resolutions in the generic case does not
depend on linkage; this avoids an implicit change of basis present in
\cite{AEB87}, see \rmkref{Brown} Under this first and purely algebraic
approach all statements are given full proofs; the next two approaches
offer interpretations of the same statements.

The second approach, taken in \secref{equiv}, is to provide canonical
equivariant forms of almost complete intersections.  The ideals one
obtains depend on a skew symmetric matrix and three vectors.  This
view of almost complete intersections was reached by analyzing the
generic ring ${\hat R}_{gen}$ constructed by Weyman~\cite{JWm18}. The
idea was to look for an open set in ${\hat R}_{gen}$ of points where
the corresponding resolution is a resolution of a perfect ideal. This
set can be explicitly described as the points where localization of
certain complex over ${\hat R}_{gen}$ is split exact. Calculating this
``splitting form'' of an ideal of grade $3$ with four generators led
to our form of almost complete intersection. One could use the
geometric technique of calculating syzygies to prove the acyclicity of
these complexes but they are identical to those from the commutative
algebra approach so we do not follow through on that.  The advantage
of this method is that one can give a geometric interpretation of the
zero set of almost complete intersection ideals. Moreover the fact,
first noticed in \cite{LLA81a}, that the skew symmetric matrix
associated to an almost complete intersection ideal can be chosen with
a $3 \times 3$ block of zeros on the diagonal is particularly natural
under this approach.

Finally, in \secref{schu}, we give a geometric interpretation of both
Gorenstein ideals and almost complete intersections of grade $3$. It
turns out that they are intersections of the so-called big open cell
with two Schubert varieties of codimension 3 in the connected
component of the orthogonal Grassmannian $\OGr{n, 2n}$ of isotropic
subspaces of dimension $n$ in a $2n$-dimensional orthogonal space. It
is interesting that in this construction the two Schubert varieties
appear together with a regular sequence by which they are linked. This
pattern generalizes from the $D_n$ root system to $E_6$, $E_7$ and
$E_8$; see Sam and Weyman \cite{SVSJWm21}.  We show that the defining
ideals are exactly the same as in commutative algebra approach, but we
indicate how one could see the graded format of the finite free
resolutions just from representation theory viewpoint.  Also, the fact
about three submaximal Pfaffians forming a regular sequence get a
clear geometric interpretation, as one can see geometrically that
their zero set has codimension 3.


\section{Almost complete intersections following Avramov and Brown}
\label{sec:local}

\noindent
For a grade $3$ perfect ideal $\mfa$ in a commutative noetherian local
ring $(R,\mfm,k)$, the minimal free resolution of the quotient ring
$R/\mfa$ has the form
\begin{equation*}
  F \deq 0 \lra F_3 \lra F_2 \lra F_1 \lra F_0 \:,
\end{equation*}
and we refer to the rank of $F_3$ as the \emph{type} of $R/\mfa$; if
$R$ is Cohen--Macaulay, then this is indeed the Cohen--Macaulay
type. Throughout the paper we treat quotients of odd and even type
separately.

By a result of Buchsbaum and Eisenbud \cite{DABDEs77} the minimal free
resolution $F$ has a structure of a skew commutative differential
graded algebra. This structure is not unique, but the induced skew
commutative algebra structure on $\Tor[R]{*}{R/\mfa}{k}$ is unique. It
provides for a classification of quotients $R/\mfa$ as worked out
Weyman~\cite{JWm89} and by Avramov, Kustin, and Miller \cite{AKM-88}.

To state the main theorems about grade $3$ almost complete
intersection ideals in local rings we introduce some matrix-related
notation.

\begin{bfhpg}[Notation]
  \label{notation}
  Let $M$ be an $m \times n$ matrix with entries in a commutative
  ring. For subsets
  \begin{equation*}
    I \deq \set{i_1,\ldots,i_k} \:\subseteq\: \set{1,\ldots,m} \qqand
    J \deq \set{j_1,\ldots, j_l} \:\subseteq\: \set{1,\ldots,n}
  \end{equation*}
  with $i_1<\cdots < i_k$ and $j_1<\cdots < j_l$ we write
  $M[i_1\ldots i_k;j_1\ldots j_l]$ for the submatrix of $M$ obtained
  by taking the rows indexed by $I$ and the columns indexed by $J$.
  At times, it is more convenient to specify a submatrix in terms of
  removal of rows and columns: The symbol
  $M[\overline{i_1\ldots i_k};\overline{j_1\ldots j_l}]$ specifies the
  submatrix of $M$ obtained by removing the rows indexed by $I$ and
  the columns indexed by $J$. These notations can also be combined:
  For example, $M[\overline{i_1\ldots i_k};j_1\ldots j_l]$ is the
  submatrix obtained by taking the rows indexed by the complement of
  $I$ and the columns indexed by $J$.

  For an $n \times n$ skew symmetric matrix $T$, the Pfaffian of $T$
  is written $\Pf{T}$. For a subset
  $\set{i_1,\ldots, i_k} \subseteq \set{1,\ldots, n}$ the Pfaffian of
  the submatrix $T[i_1\ldots i_k;i_1\ldots i_k]$ is written
  $\Pf[i_1\ldots i_k]{T}$ while the Pfaffian of
  $T[\overline{i_1 \ldots i_k};\overline{i_1 \ldots i_k}]$ is written
  $\Pfbar[T]{i_1\ldots i_k}$.
\end{bfhpg}

\begin{thm}
  \label{thm:odd-local}
  Let $n \ge 5$ be an odd number.  Let $(R,\mfm,k)$ be a local ring
  and $\mfa \subset R$ a grade $3$ almost complete intersection ideal
  such that $R/\fa$ is of type $n-3$. There exists an $n\times n$ skew
  symmetric block matrix
  \begin{equation*}
    \label{eq:Teven}
    U \deq 
    \left(
      \begin{array}[c]{c|c}
        O & B\\[1pt]
        \hline
        \phantom{\raisebox{3pt}{$|$}}\mspace{-10mu}-B^T\mspace{-10mu}
        \phantom{\raisebox{3pt}{$|$}} & A
      \end{array}
    \right)
    \deq
    \left(
      \begin{array}{ccc|ccc}
        0 & 0 & 0 & t_{14} & t_{15} & \hspace{-1ex} \dots\\[.3ex]
        0 & 0 & 0 & t_{24} & t_{25} & \hspace{-1ex} \dots\\[.3ex]
        0 & 0 & 0 & t_{34} & t_{35} & \hspace{-1ex} \dots\\[.3ex] \hline        
        -t_{14} & -t_{24} & -t_{34} & 0 & t_{45} & \hspace{-1ex} \dots \\[.3ex]
        -t_{15} & -t_{25} & -t_{35} & -t_{45} & 0 & \hspace{-1ex} \dots \\[-.6ex]        
        \vdots & \vdots & \vdots & \vdots & \vdots & \hspace{-1ex} \smash{\ddots}        
      \end{array}
    \right)
  \end{equation*}
  with entries in $\mfm$ such that the minimal free resolution of
  $R/\mfa$,
  \begin{equation*}
    F \deq 0 \lra R^{n-3} \xra{\partial_3} R^n \xra{\partial_2} R^4 \xra{\partial_1} R \:,
  \end{equation*}
  has differentials
  \begin{gather*}
    \partial_3 \deq \left(
      \begin{array}[c]{c} \:B \\
        \hline A \:
      \end{array} \right) \:,
    \\[.5ex]
    \partial_2 \deq \setlength\arraycolsep{4pt}
    \begin{pmatrix}
      \Pf{A} & 0 & 0 & -\Pfbar{234} & \Pfbar{235}
      & \hspace{-2pt}\cdots\hspace{-2pt} & \Pfbar{23n} \\
      0 & \Pf{A} & 0 & -\Pfbar{134} & \Pfbar{135}
      & \hspace{-2pt}\cdots\hspace{-2pt} & \Pfbar{13n} \\
      0 & 0 & \Pf{A} & -\Pfbar{124} & \Pfbar{125}
      & \hspace{-2pt}\cdots\hspace{-2pt} & \Pfbar{12n} \\
      \Pfbar{1} & -\Pfbar{2} & \Pfbar{3} & -\Pfbar{4} & \Pfbar{5} &
      \hspace{-2pt}\cdots\hspace{-2pt} & \Pfbar{n}
    \end{pmatrix} \:,
  \end{gather*}
  and
  \begin{equation*}
    \setlength\arraycolsep{4pt}
    \partial_1 \deq
    \begin{pmatrix}
      -\Pfbar{1} & \Pfbar{2} & -\Pfbar{3} & \Pf{A}
    \end{pmatrix} \:.
  \end{equation*}
  In particular, $\fa$ is generated by $\Pfbar{1}$, $\Pfbar{2}$,
  $\Pfbar{3}$, and $\Pf{A}$. Moreover, the multiplicative structure on
  $\Tor{*}{R/\fa}{k}$ is of class $\mathbf{H}(3,2)$ if $R/\mfa$ is of
  type $2$ and otherwise of class $\mathbf{H}(3,0)$.
\end{thm}

The proof of this theorem is given in \ref{proof-odd} and the next
theorem is proved in \ref{proof-even}.

\begin{thm}
  \label{thm:even-local}
  Let $n \ge 6$ be an even number.  Let $(R,\mfm,k)$ be a local ring
  and $\mfa \subset R$ a grade $3$ almost complete intersection ideal
  such that $R/\fa$ is of type $n-3$. There exists an $n\times n$ skew
  symmetric block matrix $U$ as in \thmref{odd-local} such that the
  minimal free resolution of $R/\mfa$,
  \begin{equation*}
    F \deq 0 \lra R^{n-3} \xra{\partial_3} R^n \xra{\partial_2}
    R^4 \xra{\partial_1} R \:,
  \end{equation*}
  has differentials
  \begin{gather*}
    \partial_3 \deq \left(
      \begin{array}[c]{c}
        \:B \\
        \hline A \:
      \end{array}
    \right) \:,
    \\[.5ex]
    \partial_2 \deq \setlength\arraycolsep{2.5pt}
    \begin{pmatrix}
      0 & 0 & 0 & \hspace{-8pt}-\Pfbar{1234} & \Pfbar{1235}
      & \cdots\hspace{-4pt} & -\Pfbar{123n} \\
      \hspace{-2pt}\Pfbar{13} & \hspace{-1pt}-\Pfbar{23}\hspace{-1pt}
      & 0 & \hspace{-2pt}\Pfbar{34} & -\Pfbar{35}
      & \cdots\hspace{-4pt} & \Pfbar{3n} \\
      \hspace{-2pt}-\Pfbar{12} & 0 & \Pfbar{23} &
      \hspace{-2pt}-\Pfbar{24} & \Pfbar{25}
      & \cdots\hspace{-4pt} & -\Pfbar{2n} \\
      0 & \hspace{2pt}\Pfbar{12} & -\Pfbar{13} &
      \hspace{-2pt}\Pfbar{14} & -\Pfbar{15}& \cdots\hspace{-4pt} &
      \Pfbar{1n}
    \end{pmatrix} \:,
  \end{gather*}
  and
  \begin{equation*}
    \partial_1 \deq
    \setlength\arraycolsep{4pt}
    \begin{pmatrix}
      \Pf{U} & \Pfbar{12} & \Pfbar{13} & \Pfbar{23}
    \end{pmatrix} \:.
  \end{equation*}
  In particular, $\fa$ is generated by $\Pf{U}$, $\Pfbar{12}$,
  $\Pfbar{13}$, and $\Pfbar{23}$. Moreover, the multiplicative
  structure on $\Tor{*}{R/\fa}{k}$ is of class $\mathbf{T}$.
\end{thm}

The proofs \thmref[Theorems~]{odd-local} and \thmref[]{even-local}
have been deferred to the next section because we obtain them by
specialization of statements about generic almost complete
intersections.


\section{Generic almost complete intersections}
\label{sec:generic}

\noindent
In this section and the appendices we deal extensively with relations
between Pfaffians of submatrices $T[i_1\ldots i_k;i_1\ldots i_k]$ of a
fixed skew symmetric matrix $T$. It is, therefore, convenient to have
the following variation on the notation from \ref{notation}:
\begin{equation}
  \label{eq:notation}
  \pf[T]{i_1\ldots i_k} \deq \Pf[i_1 \ldots i_k]{T} \qqand
  \pfbar[T]{i_1\ldots i_k} \deq \Pfbar[T]{i_1 \ldots i_k} \:.
\end{equation}
It emphasizes the subset, which changes, over the matrix, which is
fixed; for homogeneity we set $\pfT[T] = \Pf{T}$.

\begin{stp}
  \label{stp:A}
  Let $n$ be a natural number and
  $\calR = \ZZ[\grt_{ij}\mid 1 \le i < j \le n]$ the polynomial
  algebra in indeterminates $\grt_{ij}$ over $\ZZ$.  Let $\calT$ be
  the $n \times n$ skew symmetric matrix with entries
  $\calT[i;j] = \grt_{ij} = -\calT[j;i]$ for $1 \le i < j \le n$ and
  zeros on the diagonal. It looks like this:
  \begin{equation*}
    \calT  \deq \begin{pmatrix}
      0 & \grt_{12} & \grt_{13} & \grt_{14} &
      \hspace{-2pt}\dots\hspace{-2pt} \\[.3ex]
      -\grt_{12} & 0 & \grt_{23} & \grt_{24} & \hspace{-1ex}\dots \\[.3ex]
      -\grt_{13} & -\grt_{23} & 0 & \grt_{34}  & \hspace{-1ex}\dots \\[.3ex]
      -\grt_{14} & -\grt_{24} & -\grt_{34} & 0  & \hspace{-1ex}\dots \\[-.6ex]
      \vdots & \vdots & \vdots & \vdots & \hspace{-1ex} \smash{\ddots}
    \end{pmatrix}.
  \end{equation*}
\end{stp}

\begin{lem}
  \label{lem:d3}
  Adopt the setup from \stpref[]{A} and denote by $\mfM$ the ideal
  generated by the indeterminates $\grt_{ij}$. Let
  $\partial \colon \calR^{n-3} \lra \calR^{n}$ be the linear map given
  by the matrix $\calT[1\ldots n;4\ldots n]$. One has
  $\partial(\calR^{n-3}) \cap \mfM^2\calR^n=
  \partial(\mfM\calR^{n-3})$.
\end{lem}

\begin{prf*}
  Let $g_4,\ldots, g_n$ and $f_1,\ldots f_{n}$ be the standard bases
  for the free modules $\calR^{n-3}$ and $\calR^{n}$. Let
  $x = \sum_{i=4}^{n} a_ig_i$ be an element of $\calR^{n-3}$; one has
  \begin{equation*}
    \partial(x) \deq \sum_{i=4}^{n} a_i\Big(\sum_{j=1}^{i-1}\grt_{ji}f_j -
    \sum_{j=i+1}^{n}\grt_{ij}f_j\Big)
    \deq \sum_{j=1}^n \Big(\sum_{i=4}^{j-1}-a_i\grt_{ij}
    + \sum_{i=j+1}^{n}a_i\grt_{ji}\Big)f_j \:.
  \end{equation*}
  Thus, if $\partial(x)$ is contained in $\mfM^2\calR^n$, then all the
  elements $a_i\grt_{ij}$ and $a_i\grt_{ji}$ belong to $\mfM^2$, which
  implies that the every $a_i$ is in $\mfM$. Thus,
  $\partial(\calR^{n-3}) \cap \mfM^2\calR^n$ is contained in
  $\partial(\mfM\calR^{n-3})$, and the opposite containment is
  trivial.
\end{prf*}

\subsection*{Quotients of even type}

\begin{thm}
  \label{thm:odd}
  Let $n \ge 5$ be an odd number; consider the ring $\calR$ and the
  $n\times n$ matrix $\calT$ from \stpref[]{A}. The homomorphisms
  given by the matrices
  \begin{gather*}
    \partial_3 \deq \calT[1\ldots n;4\ldots n] \:,
    \\[0.5ex]
    \partial_2 \deq \setlength\arraycolsep{4pt}
    \begin{pmatrix}
      \pfbar{123} & 0 & 0 & -\pfbar{234} & \pfbar{235}
      & \hspace{-2pt}\cdots\hspace{-2pt} & \pfbar{23n} \\
      0 & \pfbar{123} & 0 & -\pfbar{134} & \pfbar{135}
      & \hspace{-2pt}\cdots\hspace{-2pt} & \pfbar{13n} \\
      0 & 0 & \pfbar{123} & -\pfbar{124} & \pfbar{125}
      & \hspace{-2pt}\cdots\hspace{-2pt} & \pfbar{12n} \\
      \pfbar{1} & -\pfbar{2} & \pfbar{3} & -\pfbar{4} & \pfbar{5} &
      \hspace{-2pt}\cdots\hspace{-2pt}& \pfbar{n}
    \end{pmatrix} \:,
  \end{gather*}
  and
  \begin{equation*}
    \setlength\arraycolsep{4pt}
    \partial_1 \deq
    \begin{pmatrix}
      -\pfbar{1} & \pfbar{2} & -\pfbar{3} & \pfbar{123}
    \end{pmatrix} \:,
  \end{equation*}
  define an exact sequence
  \begin{equation*}
    \calF \deq 0 \lra \calR^{n-3} \xra{\partial_3} \calR^{n} \xra{\partial_2}
    \calR^4 \xra{\partial_1} \calR \:.
  \end{equation*}
  That is, denoting by $\mfA_n$ the ideal generated by the entries of
  $\partial_1$, the complex $\calF$ is a free resolution of
  $\calR/\mfA_n$. Moreover, the ideal $\mfA_n$ is perfect of grade
  $3$.
\end{thm}

The proof of this theorem relies on a series of technical results that
we defer to \secref[Appendix~]{odd}. The proof shows how they come
together.

\begin{prf*}
  It follows from \lemref{setup-odd} that $\calF$ is a complex. The
  expected ranks of the homomorphisms $\partial_3$, $\partial_2$, and
  $\partial_1$ are $n-3$, $3$, and $1$. To show that the complex is
  exact at $\calR^{n-3}$, $\calR^{n}$, and $\calR^{4}$ it suffices by
  the Buchsbaum--Eisenbud Acyclicity Criterion \cite{DABDEs73} to
  verify the inequalities
  \begin{equation*}
    \grd{\calR}{I_{n-3}(\partial_3)} \ge 3\:,\quad
    \grd{\calR}{I_{3}(\partial_2)} \ge 2\:,\qand
    \grd{\calR}{I_{1}(\partial_1)} \ge 1\:,
  \end{equation*}
  where as usual $I_r(\partial)$ denotes the ideal generated by the
  $r \times r$ minors of $\partial$. By \lemref{regseq-odd} the ideal
  $I_1(\partial_1) = \mfA_n$ has grade at least $3$. By
  \lemref{d3-minors-odd} the radical $\sqrt{\mfA_n}$ is contained in
  $\sqrt{I_{n-3}(\partial_3)}$, so
  $\grd{\calR}{I_{n-3}(\partial_3)} \ge \grd{\calR}{\mfA_n} \ge 3$
  holds. By \prpref{minors-of-d2-odd}, the generators of
  $I_{3}(\partial_2)$ are products of generators of the ideals
  $I_{n-3}(\partial_3)$ and $I_1(\partial_1) = \mfA_n$.  It follows
  that the radical $\sqrt{I_3(\partial_2)}$ contains $\sqrt{\mfA_n}$,
  so one also has $\grd{\calR}{I_{3}(\partial_2)} \ge 3$. Thus,
  $\calF$ is a free resolution of $\calR/\mfA_n$; in particular, the
  projective dimension of $\calR/\mfA_n$ is at most $3$. As the grade
  of $\mfA_n$ is at least $3$, it follows that $\mfA_n$ is perfect of
  grade $3$.
\end{prf*}

The following commentary also applies to the proof of \thmref{even}.
\begin{rmk}
  \label{rmk:how}
  The proof of \thmref{odd} is based on establishing containments
  among radicals to ensure that the rank conditions in the
  Buchsbaum--Eisenbud Acyclicity Criterion are met; per
  \thmcite[2.1]{DABDEs74} the conclusion that $\calF$ is exact implies
  that the radicals $\sqrt{I_{n-3}(\partial_3)}$,
  $\sqrt{I_3(\partial_2)}$, and $\sqrt{I_1(\partial_1)}$ agree.

  The inspiration for the pivotal \prpref{minors-of-d2-odd} came from
  the same paper, namely from the Buchsbaum--Eisenbud Structure
  Theorems which, in the guise of \corcite[5.1]{DABDEs74}, say that
  for $\calF$ to be a resolution the equality
  $I_{n-3}(\partial_3)I_1(\partial_1) = I_3(\partial_2)$ must
  hold. The vehicle for the proof of \prpref{minors-of-d2-odd} is a
  relation between the sub-Pfaffians and general minors of a skew
  symmetric matrix; it was first discovered by Brill \cite{JBr04} and
  reproved by us in \cite{CVW-20a} using Knuth's \cite{DEK96}
  combinatorial approach to Pfaffians in the same vein as in
  Appendices \secref[]{knuth}--\secref[]{odd}.
\end{rmk}

As noticed in \cite{LLA81a} one can replace the upper left $3\times 3$
block in $\calT$ with a block of zeros without changing the ideal
$\mfA_n$.

\begin{lem}
  \label{lem:odd-bis}
  Let $n \ge 5$ be an odd number and $T = (t_{ij})$ an $n \times n$
  skew symmetric matrix with entries in a commutative ring $R$. Let
  $U$ be the matrix obtained from $T$ by replacing the upper left
  $3 \times 3$ block by a block of zeros; i.e.\
  \begin{equation*}
    U  \deq \begin{pmatrix}
      0 & 0 & 0 & t_{14} & t_{15} &
      \hspace{-2pt}\dots\hspace{-2pt} \\[.3ex]
      0 & 0 & 0 & t_{24} & t_{25} & \hspace{-1ex}\dots \\[.3ex]
      0 & 0 & 0 & t_{34} & t_{35}  & \hspace{-1ex}\dots \\[.3ex]
      -t_{14} & -t_{24} & -t_{34} & 0  & t_{45} & \hspace{-1ex}\dots \\
      -t_{15} & -t_{25} & -t_{35} & -t_{45} & 0 &\hspace{-1ex}\dots \\[-.6ex]
      \vdots & \vdots & \vdots & \vdots & \vdots & \hspace{-1ex} \smash{\ddots}
    \end{pmatrix}.
  \end{equation*}
  There is an equality of ideals in $R$,
  \begin{equation*}
    \bigl(\Pfbar[T]{1}\,,\: \Pfbar[T]{2}\,,\: \Pfbar[T]{3}\,,\: \Pfbar[T]{123}\bigr) \deq
    \bigl(\Pfbar{1}\,,\: \Pfbar{2}\,,\: \Pfbar{3}\,,\: \Pfbar{123}\bigr)\:.
  \end{equation*}
\end{lem}

\begin{prf*}
  Notice that $\Pfbar[T]{12i} = \Pfbar{12i}$ holds for
  $i \in \set{3,\ldots,n}$. \lemref{exp} applied with
  $u_1\ldots u_k = 2\ldots n$ and $\ell=1$ now yields
  \begin{equation*}
    \begin{aligned}
      \Pfbar[T]{1} & \deq \sum_{i=3}^{n}t_{2i}(-1)^{i-1}\Pfbar[T]{12i} \\
      & \deq t_{23}\Pfbar{123} +
      \sum_{i=4}^{n}t_{2i}(-1)^{i-1}\Pfbar{12i} \deq t_{23}\Pfbar{123}
      + \Pfbar{1}\:.
    \end{aligned}
  \end{equation*}
  Similarly, one gets
  \begin{equation*}
    \begin{aligned}
      \Pfbar[T]{2} & \deq \sum_{i=3}^{n}t_{1i}(-1)^{i-1}\Pfbar[T]{12i}  \\
      & \deq t_{13}\Pfbar{123} +
      \sum_{i=4}^{n}t_{1i}(-1)^{i-1}\Pfbar{12i} \deq t_{13}\Pfbar{123}
      + \Pfbar{2}
    \end{aligned}
  \end{equation*}
  and
  \begin{equation*}
    \begin{aligned}
      \Pfbar[T]{3} & \deq t_{12}\Pfbar[T]{123}
      + \sum_{i=4}^{n}t_{1i}(-1)^{i-1}\Pfbar[T]{13i} \\
      & \deq t_{12}\Pfbar{123} +
      \sum_{i=4}^{n}t_{1i}(-1)^{i-1}\Pfbar{13i} \deq t_{12}\Pfbar{123}
      + \Pfbar{3} \:.
    \end{aligned}
  \end{equation*}
  The asserted equality of ideals is immediate from these three
  expressions.
\end{prf*}

\begin{prp}
  \label{prp:odd-bis}
  Let $n \ge 5$ be an odd number. Consider the ring $\calR$ and the
  $n\times n$ matrix $\calT$ from \stpref[]{A} as well as the ideal
  $\mfA_n$ from \thmref{odd}. Let $\calU$ be the matrix obtained from
  $\calT$ by replacing the upper left $3 \times 3$ block by a block of
  zeros; i.e.\
  \begin{equation*}
    \calU  \deq \begin{pmatrix}
      0 & 0 & 0 & \grt_{14} & \grt_{15} &
      \hspace{-2pt}\dots\hspace{-2pt} \\[.3ex]
      0 & 0 & 0 & \grt_{24} & \grt_{25} & \hspace{-1ex}\dots \\[.3ex]
      0 & 0 & 0 & \grt_{34} & \grt_{35}  & \hspace{-1ex}\dots \\[.3ex]
      -\grt_{14} & -\grt_{24} & -\grt_{34} & 0  & \grt_{45} & \hspace{-1ex}\dots \\
      -\grt_{15} & -\grt_{25} & -\grt_{35} & -\grt_{45} & 0 &\hspace{-1ex}\dots \\[-.6ex]
      \vdots & \vdots & \vdots & \vdots & \vdots & \hspace{-1ex} \smash{\ddots}
    \end{pmatrix}.
  \end{equation*}
  The homomorphisms given by the matrices
  \begin{gather*}
    \partial_3 \deq \calU[1\ldots n;4\ldots n] \:,
    \\[0.5ex]
    \partial_2 \deq \setlength\arraycolsep{4pt}
    \begin{pmatrix}
      \pfbar[\calU]{123} & 0 & 0 & -\pfbar[\calU]{234} &
      \pfbar[\calU]{235}
      & \hspace{-2pt}\cdots\hspace{-2pt} & \pfbar[\calU]{23n} \\
      0 & \pfbar[\calU]{123} & 0 & -\pfbar[\calU]{134} &
      \pfbar[\calU]{135}
      & \hspace{-2pt}\cdots\hspace{-2pt} & \pfbar[\calU]{13n} \\
      0 & 0 & \pfbar[\calU]{123} & -\pfbar[\calU]{124} &
      \pfbar[\calU]{125}
      & \hspace{-2pt}\cdots\hspace{-2pt} & \pfbar[\calU]{12n} \\
      \pfbar[\calU]{1} & -\pfbar[\calU]{2} & \pfbar[\calU]{3} &
      -\pfbar[\calU]{4} & \pfbar[\calU]{5} &
      \hspace{-2pt}\cdots\hspace{-2pt}& \pfbar[\calU]{n}
    \end{pmatrix} \:,
  \end{gather*}
  and
  \begin{equation*}
    \setlength\arraycolsep{4pt}
    \partial_1 \deq
    \begin{pmatrix}
      -\pfbar[\calU]{1} & \pfbar[\calU]{2} & -\pfbar[\calU]{3} &
      \pfbar[\calU]{123}
    \end{pmatrix} \:,
  \end{equation*}
  define a free resolution
  $\calL = 0 \lra \calR^{n-3} \xra{\partial_3} \calR^{n}
  \xra{\partial_2} \calR^4 \xra{\partial_1} \calR$ of
  $\calR/\mfA_n$. In particular, the ideal $\mfA_n$ is generated by
  $\pfbar[\calU]{1}$, $\pfbar[\calU]{2}$, $\pfbar[\calU]{3}$, and
  $\pfbar[\calU]{123}$.
\end{prp}

\begin{prf*}
  By \lemref{odd-bis} one has
  \begin{equation*}
    \bigl(\pfbar[\calU]{1},\, \pfbar[\calU]{2},\,
    \pfbar[\calU]{3},\, \pfbar[\calU]{123}\bigr) \deq \bigl(\pfbar{1},\, \pfbar{2},\,
    \pfbar{3},\, \pfbar{123}\bigr) \deq \mfA_n \:.
  \end{equation*}
  Next we show how to obtain the free resolution $\calL$ from the
  resolution $\calF$ from \thmref{odd}. To distinguish the
  differentials on the resolutions we introduce superscripts $\calL$
  and $\calF$. Consider the matrix
  \begin{equation*}
    S \deq
    \begin{pmatrix}
      1 & 0 & 0 & 0 \\
      0 & 1 & 0 & 0 \\
      0 & 0 & 1 & 0 \\
      \grt_{23} & -\grt_{13} & \grt_{12} & 1
    \end{pmatrix}
    \qtext{with inverse}
    S^{-1} \deq
    \begin{pmatrix}
      1 & 0 & 0 & 0 \\
      0 & 1 & 0 & 0 \\
      0 & 0 & 1 & 0 \\
      -\grt_{23} & \grt_{13} & -\grt_{12} & 1
    \end{pmatrix} \:.
  \end{equation*}
  As in the proof of \lemref{odd-bis} one has
  \begin{align*}
    \tag{1}
    \pfbar{1} & \deq \grt_{23}\pfbar[\calU]{123} + \pfbar[\calU]{1} \:, \\
    \tag{2}
    \pfbar{2} & \deq \grt_{13}\pfbar[\calU]{123} + \pfbar[\calU]{2} \:, \\
    \tag{3}
    \pfbar{3} & \deq \grt_{12}\pfbar[\calU]{123} + \pfbar[\calU]{3} \:,
                \quad\text{and} \\
    \tag{4} \pfbar{123} & \deq \pfbar[\calU]{123} \:.
  \end{align*}
  These identities show that one has
  $\partial_1^\calL = \partial_1^\calF S$. Thus, the matrices
  \begin{equation*}
    \partial_3^\calF \:,\quad S^{-1}\partial_2^\calF\:, \qand
    \partial_1^\calF S
  \end{equation*}
  determine a free resolution of $\calR/\mfA_n$. As the submatrices
  $\partial_3^\calL = \calU[1\ldots n;4\ldots n]$ and
  $\partial_3^\calF = \calT[1\ldots n;4\ldots n]$ agree, it suffices
  to show that $\partial_2^\calL = S^{-1}\partial_2^\calF$ holds.

  For indices $1 \le i \le n$ one has
  \begin{equation*}
    \pfbar{12i} \deq \pfbar[\calU]{12i}\:,\quad
    \pfbar{13i} \deq \pfbar[\calU]{13i}\:,\qand
    \pfbar{23i} \deq \pfbar[\calU]{23i}\:.    
  \end{equation*}
  It follows that the first three rows in the matrices
  $\partial_2^\calL$, $\partial_2^\calF$, and $S^{-1}\partial_2^\calF$
  agree.  We now focus of the fourth rows of $\partial_2^\calL$ and
  $S^{-1}\partial_2^\calF$. The first three entries in the fourth rows
  agree by the identities $(1)$, $(2)$, and $(3)$. Now fix
  $j \in \set{4,\dots,n}$. Another application of \lemref{exp} yields
  \begin{equation*}
    \begin{aligned}
      \pfbar{j} & \deq \grt_{12}\pfbar{12j} - \grt_{13}\pfbar{13j} \\
      & \hspace{3pc} {} + \sum_{i=4}^{j-1}\Sign{i}\grt_{1i}\pfbar{1ij}
      + \sum_{i=j+1}^{n}\Sign{i-1}\grt_{1i}\pfbar{1ij} \:.
    \end{aligned}
  \end{equation*}
  One has $\pfbar{1ij} = \grt_{23}\pfbar{123ij} + \pfbar[\calU]{1ij}$,
  again by \lemref{exp}, and therefore
  \begin{align*}
    \pfbar{j} & - \grt_{12}\pfbar{12j} + \grt_{13}\pfbar{13j}
    \\
              & \deq {} \grt_{23}\Big(
                \sum_{i=4}^{j-1}\Sign{i}\grt_{1i}\pfbar{123ij} +
                \sum_{i=j+1}^{n}\Sign{i-1}\grt_{1i}\pfbar{123ij} \Big) \\
              &  \hspace{3pc} {} +
                \sum_{i=1}^{j-1}\Sign{i}\grt_{1i}\pfbar[\calU]{1ij}
                + \sum_{i=j+1}^{n}\Sign{i-1}\grt_{1i}\pfbar[\calU]{1ij}
    \\
              & \deq  \grt_{23}\pfbar{23j} + \pfbar[\calU]{j} \:,
  \end{align*}
  where the last equality follows from two applications of
  \lemref{exp}. This identity shows that the fourth row entries of
  $\partial_2^\calL$ and $S^{-1}\partial_2^\calF$ agree in column $j$.
\end{prf*}

One could also establish \prpref{odd-bis} as follows: After invoking
\lemref{odd-bis}, repeat the proof of \thmref{odd} noticing at every
step that the conclusions remain valid after evaluation at
$\grt_{12} = \grt_{13} = \grt_{23} =0$.

\begin{bfhpg}[Proof of Theorem \ref{thm:odd-local}]
  \label{proof-odd}
  An almost complete intersection ideal of grade $3$ is by
  \prpcite[5.2]{DABDEs77} linked to a Gorenstein ideal of grade $3$,
  and Brown \prpcite[4.3]{AEB87} uses this to show that there exits a
  skew symmetric matrix $T$ with entries in $\mfm$ such that $\mfa$ is
  generated by the Pfaffians $\Pfbar[T]{1}$, $\Pfbar[T]{2}$,
  $\Pfbar[T]{3}$, and $\Pfbar[T]{123}$. \lemref{odd-bis} shows than
  one can replace the upper left $3\times 3$ block in $T$ with zeroes
  and arrive at the asserted block matrix $U$.

  Adopt \stpref{A}. Let $\calR \to R$ be given by
  $\tau_{ij} \mapsto t_{ij}$; it makes $R$ an $\calR$-algebra and maps
  Pfaffians of submatrices of $\calT$ to the corresponding Pfaffians
  of submatrices of $U$, i.e.\ $\pfbar{123}$ maps to $\Pfbar{123}$
  etc. Let $\calF$ be the free resolution of $\calR/\mfA_n$ from
  \thmref{odd}. As one has $R/\mfa = \calR/\mfA_n \otimes_\calR R$ and
  $\mfa$ has grade $3$ it follows from Bruns and Vetter
  \thmcite[3.5]{bruvet} that
  \begin{equation*}
    \tag{0}
    F \deq \tp[\calR]{\calF}{R}
  \end{equation*}
  is a free resolution of $R/\mfa$ over $R$, and it is minimal as the
  differentials are given by matrices with entries in $\mfm$.

  We now establish parts of a multiplicative structure on $F$: just
  enough to determine the multiplicative structure on the $k$-algebra
  $\Tor{*}{R/\fa}{k}$. Let $e_1,\ldots,e_4$, $f_1,\ldots,f_n$, and
  $g_1,\ldots,g_{n-3}$ be the standard bases for the free modules
  $F_1$, $F_2$, and $F_3$. From the three obvious Koszul relations one
  gets
  \begin{align*}
    \partial_2(e_4e_1) & \deq  \Pfbar{123}e_1 + \Pfbar{1}e_4
                         \deq \partial_2(f_1) \:, \\
    \partial_2(e_4e_2) & \deq  \Pfbar{123}e_2 - \Pfbar{2}e_4
                         \deq \partial_2(f_2) \:,\text{ and}\\
    \partial_2(e_4e_3) & \deq  \Pfbar{123}e_3 + \Pfbar{3}e_4
                         \deq \partial_2(f_3) \:.
  \end{align*}
  Thus one can set
  \begin{equation*}
    \tag{1}
    e_4e_1 \deq f_1 \:,\quad e_4e_2 \deq f_2 \:,\qand
    e_4e_3 \deq f_3 \:.
  \end{equation*}
  These three products in $F$ induce non-trivial products in
  $\Tor{*}{R/\fa}{k}$. Applying \lemref{exp} the same way as in the
  proof of \lemref{odd-bis} one gets:
  \begin{align*}
    \partial_2(e_1e_2) & \deq  -\Pfbar{1}e_2 - \Pfbar{2}e_1 
                         \deq \partial_2\Big(\sum_{i=4}^nt_{3i}f_i\Big) \:,
    \\
    \partial_2(e_2e_3) & \deq  \Pfbar{2}e_3 + \Pfbar{3}e_2 
                         \deq \partial_2\Big(\sum_{i=4}^nt_{1i}f_i\Big) \:,\text{ and}
    \\
    \partial_2(e_3e_1) & \deq  -\Pfbar{3}e_1 + \Pfbar{1}e_3
                         \deq \partial_2\Big(\sum_{i=4}^n t_{2i}f_i\Big) \:.
  \end{align*}
  Thus one can set
  \begin{equation*}
    \tag{2}
    e_1e_2 \deq \sum_{i=4}^nt_{3i}f_i \:,\quad e_2e_3 \deq \sum_{i=4}^n t_{1i}f_i \:,\qand
    e_3e_1 \deq \sum_{i=4}^n t_{2i}f_i \:.
  \end{equation*}
  The products $(2)$ induce trivial products in
  $\Tor{*}{R/\fa}{k}$. We have now accounted for all products of
  elements from $F_1$, so $R/\mfa$ is of class $\mathbf{H}(3,q)$,
  where $q$ denotes the dimension of the subspace
  $\Tor{1}{R/\fa}{k}\cdot \Tor{2}{R/\fa}{k}$ of $\Tor{3}{R/\fa}{k}$;
  see \thmcite[2.1]{AKM-88}.
  
  For $j \in \set{1,2,3}$ one has
  \begin{equation*}
    \partial_3(e_4f_j) \deq \Pfbar{123}f_j - e_4\bigl(\Pfbar{123}e_j
    \pm \Pfbar{j}e_4\bigr) \deq 0
  \end{equation*}
  by $(1)$. Since $\partial_3$ is injective, one has
  \begin{equation*}
    e_4f_1 \deq     e_4f_2 \deq     e_4f_3 \deq  0 \:.
  \end{equation*}
  For $j\in\set{4,\ldots,n}$ one gets
  \begin{equation*}
    \tag{3}
    \begin{aligned}
      \partial_3(e_4f_j) & \deq \Pfbar{123}f_j \\
      & \hspace{3pc} {} + \Sign{j}e_4 \bigl(\Pfbar{23j}e_1 +
      \Pfbar{13j}e_2 + \Pfbar{12j}e_3 \Pfbar{j}e_4\bigr)
      \\
      & \deq \Sign{j}\bigl(\Pfbar{23j}f_1 + \Pfbar{13j}f_2 +
      \Pfbar{12j}f_3\bigr) - \Pfbar{123}f_j \:,
    \end{aligned}
  \end{equation*}
  again by $(1)$. Thus, for $n \ge 7$ one has
  $\partial_3(e_4f_j) \in \mfm^2F_2$. In this case it follows from
  \lemref{d3} and $(0)$ that there is an element $x_j\in\mfm F_3$ with
  $\partial_3(x_j) = \partial_3(e_4f_j)$, so by injectivity of
  $\partial_3$ one has $e_4f_j = x_j$; in particular this product
  induces a trivial product in $\Tor{*}{R/\fa}{k}$.  For $n=5$ one has
  $j\in\set{4,5}$ and $(3)$ specializes to
  \begin{align*}
    \partial_3(e_4f_4)  &\deq t_{15}f_1 + t_{25}f_2 + t_{35}f_3
                          - t_{45}f_4 \deq \partial_3(g_2) \quad\text{and}\\
    \partial_3(e_4f_5)  &\deq -(t_{14}f_1 + t_{24}f_2 + t_{34}f_3
                          + t_{45}f_5) \deq -\partial_3(g_1) \:.
  \end{align*}
  As $\partial_3$ is injective, this shows that in this case one has
  \begin{equation*}
    e_4f_4 \deq g_2 \qqand e_4f_5 \deq -g_1 \:.
  \end{equation*}
  To prove the assertion about the multiplicative structure on
  $\Tor{*}{R/\fa}{k}$, it suffices to show that there are no further
  non-zero products in $\Tor{1}{R/\fa}{k}\cdot \Tor{2}{R/\fa}{k}$. To
  this end it suffices by \lemref{d3} and $(0)$ to shows that
  $\partial_3(e_if_j)$ belongs to $\mfm^2F_2$ for indices
  $1\le i \le 3$ and $1 \le j \le n$.  For $1 \le i,j \le 3$ one has
  \begin{equation*}
    \partial_3(e_if_j) \deq \Sign{i}\Pfbar{i}f_j - e_i\bigl(\Pfbar{123}e_j +
    \Sign{j-1}\Pfbar{j}e_4\bigr) \:.
  \end{equation*}
  This is indeed in $\mfm^2F_2$ as $\Pfbar{i}$ and $\Pfbar{j}$ belong
  to $\mfm^2$, and $e_ie_j \in \mfm F_2$ by $(2)$.  For indices
  $1 \le i \le 3$ and $4 \le j \le n$ one has
  \begin{align*}
    \partial_3(e_if_j) &\deq \Sign{i}\Pfbar{i}f_j \\
                       & \hspace{1.5pc} - e_i\Sign{j-1}\bigl(\Pfbar{23j}e_1 + \Pfbar{13j}e_2 + \Pfbar{12j}e_3 +
                         \Pfbar{j}e_4\bigr) \:.
  \end{align*}
  As above $\Pfbar{i}$ and $\Pfbar{j}$ belong to $\mfm^2$, and the
  products $e_ie_1$, $e_ie_2$, and $e_ie_3$ belong to $\mfm F_2$ by
  $(2)$. \qed
\end{bfhpg}

\subsection*{Quotients of odd type}

\begin{rmk}
  \label{rmk:Brown}
  In Brown's work \cite{AEB87}, the statements to the effect that all
  almost complete intersection ideals come from skew symmetric
  matrices---Propositions~4.2 and 4.3 in \emph{loc.\ cit.\ }as cited
  in our proofs of \thmref[Theorems~]{odd-local} and
  \thmref[]{even-local}---are separated from the descriptions of the
  free resolutions: Propositions~3.2 and 3.3 in \emph{loc.\ cit.} The
  proofs of all four statements in \cite{AEB87} rely on the fact from
  \cite{DABDEs77} that almost complete intersection ideals are linked
  to Gorenstein ideals, but compare the proofs of \cite[Propositions
  3.2 and 4.2]{AEB87} for almost complete intersections of odd type:
  The linking sequence used in the description of the free resolution
  is different from the one used to associate a skew symmetric matrix;
  a change of basis argument is thus required to reconcile the
  two. Our explicit construction of the free resolution in the generic
  case, \thmref[Theorems~]{odd} and \thmref[]{even}, allows us to
  avoid such issues in the proofs of \thmref[Theorems~]{odd-local} and
  \thmref[]{even-local}.
\end{rmk}

\begin{thm}
  \label{thm:even}
  Let $n \ge 6$ be an even number; consider the ring $\calR$ and the
  $n\times n$ matrix $\calT$ from \stpref[]{A}. The homomorphisms
  given by the matrices
  \begin{gather*}
    \partial_3 \deq \calT[1\ldots n;4\ldots n] \:,
    \\[0.5ex]
    \partial_2 \deq \setlength\arraycolsep{2.5pt}
    \begin{pmatrix}
      0 & 0 & 0 & \hspace{-8pt}-\pfbar{1234} & \pfbar{1235}
      & \cdots\hspace{-4pt} & -\pfbar{123n} \\
      \hspace{-2pt}\pfbar{13} & \hspace{-1pt}-\pfbar{23}\hspace{-1pt}
      & 0 & \hspace{-2pt}\pfbar{34} & -\pfbar{35}
      & \cdots\hspace{-4pt} & \pfbar{3n} \\
      \hspace{-2pt}-\pfbar{12} & 0 & \pfbar{23} &
      \hspace{-2pt}-\pfbar{24} & \pfbar{25}
      & \cdots\hspace{-4pt} & -\pfbar{2n} \\
      0 & \hspace{-1pt}\pfbar{12}\hspace{-1pt} & -\pfbar{13} &
      \hspace{-2pt}\pfbar{14} & -\pfbar{15}& \cdots\hspace{-4pt} &
      \pfbar{1n}
    \end{pmatrix} \:,
  \end{gather*}
  and
  \begin{equation*}
    \setlength\arraycolsep{4pt}
    \partial_1 \deq
    \begin{pmatrix}
      \pfT & \pfbar{12} & \pfbar{13} & \pfbar{23}
    \end{pmatrix}
  \end{equation*}
  define an exact sequence
  \begin{equation*}
    \calF \deq 0 \lra \calR^{n-3} \xra{\partial_3} \calR^{n} \xra{\partial_2}
    \calR^4 \xra{\partial_1} \calR \:.
  \end{equation*}
  That is, denoting by $\mfA_n$ the ideal generated by the entries of
  $\partial_1$, the complex $\calF$ is a free resolution of
  $\calR/\mfA_n$. Moreover, the ideal $\mfA_n$ is perfect of grade
  $3$.
\end{thm}

\begin{prf*}
  The proof of \thmref{odd} applies, one only needs to replace the
  references to \lemref[]{setup-odd}--\prpref[]{minors-of-d2-odd} with
  references to \lemref[]{setup-even}--\prpref[]{minors-of-d2-even}.
\end{prf*}

\begin{lem}
  \label{lem:even-bis}
  Let $n \ge 6$ be an even number and $T = (t_{ij})$ an $n \times n$
  skew symmetric matrix with entries in a commutative ring $R$. Let
  $U$ be the matrix obtained from $T$ by replacing the upper left
  $3 \times 3$ block by a block of zeros; see \lemref{odd-bis}.  There
  is an equality of ideals in $R$
  \begin{equation*}
    \bigl(\Pf{T}\,,\: \Pfbar[T]{12}\,,\: \Pfbar[T]{13}\,,\: \Pfbar[T]{23} \bigr)
    \deq \bigl(\Pf{U}\,,\: \Pfbar{12}\,,\: \Pfbar{13}\,,\: \Pfbar{23} \bigr)\:.
  \end{equation*}
\end{lem}

\begin{prf*}
  First notice that one has
  \begin{equation*}
    \tag{1}
    \Pfbar[T]{12} \deq \Pfbar{12}\:,\quad
    \Pfbar[T]{13} \deq \Pfbar{13}\:,\qand
    \Pfbar[T]{23} \deq \Pfbar{23}\:.
  \end{equation*}
  \lemref{exp} applied with $u_1\ldots u_k = 1\ldots n$ and $\ell=1$
  yields
  \begin{equation*}
    \tag{2}
    \begin{aligned}
      \Pf{T} & \deq \sum_{i=2}^{n}t_{1i}(-1)^{i}\Pfbar[T]{1i} \\
      & \deq t_{12}\Pfbar[T]{12} - t_{13}\Pfbar[T]{13} +
      \sum_{i=4}^{n}t_{1i}(-1)^{i}\Pfbar[T]{1i} \:.
    \end{aligned}
  \end{equation*}
  For $i \ge 4$ the same lemma applied with
  $u_1\ldots u_k = 2\ldots n \setminus i$ and $\ell = 2$ yields
  \begin{equation*}
    \tag{3}
    \Pfbar[T]{1i} \deq t_{23}\Pfbar[T]{123i} + \Pfbar{1i} \:.
  \end{equation*}
  From $(2)$, $(3)$, and further applications of \lemref[]{exp} one
  now gets
  \begin{equation*}
    \tag{4}
    \begin{aligned}
      \Pf{T} & - t_{12}\Pfbar[T]{12} + t_{13}\Pfbar[T]{13} \\
      & \deq t_{23}\sum_{i=4}^{n}t_{1i}(-1)^{i}\Pfbar[T]{123i}
      + \sum_{i=4}^{n}t_{1i}(-1)^{i}\Pfbar{1i} \\
      & \deq t_{23}\Pfbar[T]{23} + \Pf{U}\:.
    \end{aligned}
  \end{equation*}
  The asserted equality of ideals is immediate from $(1)$ and $(4)$.
\end{prf*}

\begin{prp}
  \label{prp:even-bis}
  Let $n \ge 6$ be an even number. Consider the ring $\calR$ and the
  \mbox{$n\times n$} matrix $\calT$ from \stpref[]{A} as well as the
  ideal $\mfA_n$ from \thmref{even}.  Let $\calU$ be the matrix
  obtained from $\calT$ by replacing the upper left $3 \times 3$ block
  by a block of zeros; see \prpref{odd-bis}.  The homomorphisms given
  by the matrices
  \begin{gather*}
    \partial_3 \deq \calU[1,\ldots,n;4,\ldots,n] \:,
    \\[0.5ex]
    \partial_2 \deq \setlength\arraycolsep{2.5pt}
    \begin{pmatrix}
      0 & 0 & 0 & \hspace{-8pt}-\pfbar[\calU]{1234} &
      \pfbar[\calU]{1235}
      & \cdots\hspace{-4pt} & -\pfbar[\calU]{123n} \\
      \hspace{-2pt}\pfbar[\calU]{13} &
      \hspace{-1pt}-\pfbar[\calU]{23}\hspace{-1pt} & 0 &
      \hspace{-2pt}\pfbar[\calU]{34} & -\pfbar[\calU]{35}
      & \cdots\hspace{-4pt} & \pfbar[\calU]{3n} \\
      \hspace{-2pt}-\pfbar[\calU]{12} & 0 & \pfbar[\calU]{23} &
      \hspace{-2pt}-\pfbar[\calU]{24} & \pfbar[\calU]{25}
      & \cdots\hspace{-4pt} & -\pfbar[\calU]{2n} \\
      0 & \hspace{-1pt}\pfbar[\calU]{12}\hspace{-1pt} &
      -\pfbar[\calU]{13} & \hspace{-2pt}\pfbar[\calU]{14} &
      -\pfbar[\calU]{15}& \cdots\hspace{-4pt} & \pfbar[\calU]{1n}
    \end{pmatrix}
  \end{gather*}
  and
  \begin{equation*}
    \partial_1 \deq
    \setlength\arraycolsep{4pt}
    \begin{pmatrix}
      \pfT[\calU] & \pfbar[\calU]{12} & \pfbar[\calU]{13} &
      \pfbar[\calU]{23}
    \end{pmatrix}
  \end{equation*}
  define a free resolution
  $\calL = 0 \lra \calR^{n-3} \xra{\partial_3} \calR^{n}
  \xra{\partial_2} \calR^4 \xra{\partial_1} \calR$ of
  $\calR/\mfA_n$. In particular, the ideal $\mfA_n$ is generated by
  $\pfT[\calU]$, $\pfbar[\calU]{12}$, $\pfbar[\calU]{13}$, and
  $\pfbar[\calU]{23}$.
\end{prp}

\begin{prf*}
  By \lemref{even-bis} one has
  \begin{equation*}
    \bigl(\pfT[\calU],\, \pfbar[\calU]{12},\,
    \pfbar[\calU]{13},\, \pfbar[\calU]{23}\bigr) \deq
    \bigl(\pfT,\, \pfbar{12},\, \pfbar{13},\, \pfbar{23}\bigr) \deq \mfA_n \:.
  \end{equation*}
  As in the proof of \prpref{odd-bis}, we proceed to show how the free
  resolution $\calL$ is obtained from the resolution $\calF$ from
  \thmref{even}. To distinguish the differentials on the resolutions
  we introduce superscripts $\calL$ and $\calF$. Consider the matrix
  \begin{equation*}
    S \deq
    \begin{pmatrix}
      1 & 0 & 0 & 0 \\
      -\grt_{12} & 1 & 0 & 0 \\
      \grt_{13} & 0 & 1 & 0 \\
      -\grt_{23} & 0 & 0 & 1
    \end{pmatrix}
    \qtext{with inverse}
    S^{-1} \deq
    \begin{pmatrix}
      1 & 0 & 0 & 0 \\
      \grt_{12} & 1 & 0 & 0 \\
      -\grt_{13} & 0 & 1 & 0 \\
      \grt_{23} & 0 & 0 & 1
    \end{pmatrix} \:.
  \end{equation*}
  Notice that one has
  \begin{equation*}
    \tag{1}
    \pfbar{12} \deq \pfbar[\calU]{12}\:,\quad
    \pfbar{13} \deq \pfbar[\calU]{13}\:,\qand
    \pfbar{23} \deq \pfbar[\calU]{23}\:.
  \end{equation*}
  As in the proof of \lemref{even-bis} one gets
  \begin{equation*}
    \tag{2}
    \pfT  \deq \grt_{23}\pfbar{23} + \pfT[\calU]\:.
  \end{equation*}
  These identities yield $\partial_1^\calL = \partial_1^\calF
  S$. Thus, the matrices $\partial_3^\calF$, $S^{-1}\partial_2^\calF$,
  and $\partial_1^\calF S$ determine a free resolution of
  $\calR/\mfA_n$. As the matrices $\partial_3^\calL$ and
  $\partial_3^\calF$ agree, it suffices to show that
  $\partial_2^\calL = S^{-1}\partial_2^\calF$ holds.

  By $(1)$ the first three columns of the matrices $\partial_2^\calL$,
  $\partial_2^\calF$, and $S^{-1}\partial_2^\calF$ agree. For indices
  $1 \le i \le n$ one has $\pfbar{123i} = \pfbar[\calU]{123i}$, so
  also the top rows in the matrices $\partial_2^\calL$,
  $\partial_2^\calF$, and $S^{-1}\partial_2^\calF$ agree.  Now fix
  $i\in\set{4,\ldots,n}$. The $(4,i)$ entry in the matrix
  $S^{-1}\partial_2^\calF$ is
  \begin{equation*}
    \grt_{23}\Sign{i-1}\pfbar{123i} + \Sign{i}\pfbar{1i} \:.
  \end{equation*}
  To see that this is indeed $\Sign{i-1}\pfbar[\calU]{123i}$, the
  $(4,i)$ entry in $\partial_2^\calL$, apply \lemref{exp} with
  $u_1\ldots u_k = 2\ldots n \setminus i$ and $\ell = 2$ to get
  \begin{equation*}
    \pfbar{1i} \deq \grt_{23}\pfbar{123i} + \pfbar[\calU]{1i} \:.
  \end{equation*}
  Similar applications of \lemref{exp} yield the identities
  \begin{equation*}
    \pfbar{2i} \deq \grt_{13}\pfbar{123i} + \pfbar[\calU]{2i}
    \qand \pfbar{3i} \deq \grt_{12}\pfbar{123i} + \pfbar[\calU]{3i}
    \:,    
  \end{equation*}
  which show that also the $(3,i)$ and $(2,i)$ entries in the two
  matrices agree.
\end{prf*}

\begin{bfhpg}[Proof of Theorem \ref{thm:even-local}]
  \label{proof-even}
  An almost complete intersection ideal of grade $3$ is by
  \prpcite[5.2]{DABDEs77} linked to a Gorenstein ideal of grade $3$,
  and Brown \prpcite[4.2]{AEB87} uses this to show that there exits a
  skew symmetric matrix $T$ with entries in $\mfm$ such that $\mfa$ is
  generated by the Pfaffians $\Pf{T}$, $\Pfbar[T]{12}$,
  $\Pfbar[T]{13}$, and $\Pfbar[T]{23}$. \lemref{even-bis} shows than
  one can replace the upper left $3\times 3$ block in $T$ with zeroes
  and arrive at the asserted block matrix $U$.  Adopt \stpref{A} and
  let $\calF$ be the free resolution of $\calR/\mfA_n$ from
  \thmref{even}. As in the proof of \thmref{odd-local} one sees~that
  \begin{equation*}
    \tag{0}
    F \deq \tp[\calR]{\calF}{R}
  \end{equation*}
  is a minimal free resolution of $R/\mfa$ over $R$.

  As in the proof of \thmref{odd-local} we proceed to determine enough
  of a multiplicative structure on $F$ to recognize the multiplicative
  structure on $\Tor{*}{R/\fa}{k}$. Let $e_1,\ldots,e_4$,
  $f_1,\ldots,f_n$, and $g_1,\ldots,g_{n-3}$ be the standard bases for
  the free modules $F_1$, $F_2$, and $F_3$. From the three obvious
  Koszul relations one gets
  \begin{align*}
    \partial_2(e_2e_3) & \deq  \Pfbar{12}e_3 - \Pfbar{13}e_2
                         \deq \partial_2(-f_1) \:, \\
    \partial_2(e_3e_4) & \deq  \Pfbar{13}e_4 - \Pfbar{23}e_3
                         \deq \partial_2(-f_3) \:,\text{ and}\\
    \partial_2(e_4e_2) & \deq  \Pfbar{23}e_2 - \Pfbar{13}e_4
                         \deq \partial_2(-f_2) \:.
  \end{align*}
  Thus one can set
  \begin{equation*}
    \tag{1}
    e_2e_3 \deq -f_1 \:,\quad e_3e_4 \deq -f_3 \:,\qand
    e_4e_2 \deq -f_2 \:.
  \end{equation*}
  These three products in $F$ induce non-trivial products in
  $\Tor{*}{R/\fa}{k}$. Repeated applications of \lemref{exp} yield:
  \begin{align*}
    \partial_2(e_1e_2) & \deq  \Pf{T}e_2 - \Pfbar{12}e_1  \\
                       & \deq -\Big(\sum_{i=4}^nt_{3i}\Sign{i}\Pfbar{123i}\Big)e_1
                         + \Big(\sum_{i=4}^nt_{3i}\Sign{i}\Pfbar{3i}\Big)e_2 \\
                       & \deq \partial_2\Big(\sum_{i=4}^n t_{3i}f_i\Big)  \:,
    \\
    \partial_2(e_1e_3) & \deq  \Pf{T}e_3 - \Pfbar{13}e_3 
                         \deq \partial_2\Big(\sum_{i=4}^nt_{2i}f_i\Big) \:,\text{ and}
    \\
    \partial_2(e_1e_4) & \deq  \Pf{T}e_4 - \Pfbar{23}e_1
                         \deq \partial_2\Big(\sum_{i=4}^n t_{1i}f_i\Big) \:.
  \end{align*}
  Thus one can set
  \begin{equation*}
    \tag{2}
    e_1e_2 \deq \sum_{i=4}^nt_{3i}f_i \:,\quad e_1e_3 \deq \sum_{i=4}^n t_{2i}f_i
    \:,\qand
    e_1e_4 \deq \sum_{i=4}^n t_{1i}f_i \:.
  \end{equation*}
  The products $(2)$ induce trivial products in
  $\Tor{*}{R/\fa}{k}$. We have now accounted for all products of
  elements from $F_1$. To prove that $R/\mfa$ is of class $\mathbf{T}$
  it suffices to show that all products of the form $e_if_j$ induce
  the zero product in $\Tor{*}{R/\fa}{k}$; see
  \thmcite[2.1]{AKM-88}. To this end, it suffices by \lemref{d3} and
  $(0)$ to shows that $\partial_3(e_if_j)$ belongs to $\mfm^2F_2$ for
  all indices $1\le i \le 4$ and $1 \le j \le n$.  One has
  \begin{equation*}
    \partial_3(e_if_j) \deq \partial_1(e_i)f_j - \partial_2(f_j)e_i \:.
  \end{equation*}
  For all $i$ one has $\partial_1(e_i) \in \mfm^2$, and for
  $1\le j \le 3$ also $\partial_2(f_j)$ belongs to $\mfm^2$. For
  $4 \le j \le n$ one has
  \begin{equation*}
    \partial_2(f_j)e_i \deq \Sign{j-1}\bigl(\Pfbar{123j}e_1 -
    \Pfbar{3j}e_2 + \Pfbar{2j}e_3 - \Pfbar{1j}e_4\bigr)e_i \:.
  \end{equation*}
  This too is in $\mfm^2F_2$ as the coefficients $\Pfbar{1j}$,
  $\Pfbar{2j}$, and $\Pfbar{3j}$ belong to $\mfm^2$ and the product
  $e_1e_i$ is in $\mfm F_2$. \qed
\end{bfhpg}


\section{The equivariant form of the format $(1,4,n,n-3)$}
\label{sec:equiv}

\noindent
In this section we give an equivariant interpretation of generic four
generated perfect ideals of codimension three.  These ideals were
already considered from a purely algebraic point of view in
\secref{generic}, and they will be treated as linear sections of
Schubert varieties in \secref{schu}.

\subsection*{Quotients of even type}
Let $n = 2m + 3$ where $m$ is a natural number.  Consider a
$2m\times 2m$ generic skew symmetric matrix $A=(c_{ij})$ and a
$3\times 2m$ generic matrix $B= (u_{ki})$.  Thus we work over a ring
\begin{equation*}
  \calA \deq \Sym[\ZZ]{\textstyle\bigwedge^2 F\oplus F\otimes G} \dis
  \ZZ[c_{ij}, u_{ki}]
\end{equation*}
where $F= \ZZ^{2m}$ and $G=\ZZ^3$ are free $\ZZ$-modules.  The ring
$\calA$ has an obvious bigrading with $|c_{ij}|= (1,0)$ and
$|u_{ki}|= (0,1)$.

\begin{prp}
  \label{prp:equiv-odd}
  Let $\set{g_1,\ldots, g_{2m}}$ be a basis for $F$ and set
  \begin{equation*}
    C \deq \sum_{1\le i<j\le 2m} c_{ij}g_i\wedge g_j \qand u_{k}
    \deq  \sum_{i=1}^{2m} u_{ki}g_{i} \ \text{ for } \ 1 \le k  \le 3 \:.
  \end{equation*}
  We denote by $C^j$ the $j^\mathrm{th}$ exterior power in
  $\bigwedge^{2j}F$.  The ideal
  \begin{equation*}
    \calI_{n} \deq ( C^m, C^{m-1}\wedge u_{1}\wedge u_{2},
    C^{m-1}\wedge u_{1}\wedge u_{3}, C^{m-1}\wedge u_{2}\wedge u_{3})
  \end{equation*}
  is a grade $3$ almost complete intersection ideal of type
  $2m = n-3$.
\end{prp}

\begin{prf*}
  The exterior powers $C^j$ have a natural description in terms of
  Pfaffians of the matrix $A$,
  \begin{equation*}
    C^j \deq \sum_{1\le i_1<\ldots <i_{2j}\le 2m}
    \Pf[i_1\ldots i_{2j}]{A} \cdot g_{i_1}\wedge\ldots\wedge g_{i_{2j}}\;.
  \end{equation*}
  Plugging these in, we see that we get the generators of the ideals
  described in this proposition from the matrix $\calU$
  in \prpref{odd-bis} via the substitutions $c_{ij}=\tau_{i+3,j+3}$ and
  $u_{i,j}=\tau_{i,j+3}$.
\end{prf*}

Let us work out the minimal free resolution of the ideal defined
above.

\begin{prp}
  \label{prp:equiv-odd-d}
  Let $n = 2m+3$ and $\calI_n$ be the ideal from
  \prpref{equiv-odd}. The minimal graded free resolution of the cyclic
  $\calA$-module $\calA/\calI_{n}$ is
  \begin{gather*}
    \calF_{\bullet} : 0 \lra F\otimes (\bigwedge^{2m}F)^{\otimes
      2}\otimes\bigwedge^3
    G\otimes \calA(-2m+1,-3) \xra{\partial_3}  \\
    (\bigwedge^{2m}F)^{\otimes 2}\otimes\bigwedge^2 G\otimes
    \calA(-2m+1,-2)\oplus\bigwedge^{2m}F\otimes\bigwedge^{2m-1}F\otimes\bigwedge^3G\otimes
    \calA(-2m+2,-3) \\
    \xra{\partial_2} \bigwedge^{2m} F\otimes \calA(-m,0)\oplus
    \bigwedge^{2m} F\otimes\bigwedge^2 G\otimes \calA(-m+1,-2)
    \xra{\partial_1} \calA.
  \end{gather*}
  The differentials $\partial_3$, $\partial_2$, $\partial_1$ are
  described in the proof below.  For every field $k$ the resolution
  $\calF_{\bullet}\otimes_{\ZZ}k$ is minimal over
  \begin{equation*}
    \calA_k \deq \Sym{\textstyle \bigwedge^2 \widebar{F}\oplus \widebar{F}\otimes \widebar{G})} \dis k[c_{ij}, u_{ki}] \:,
  \end{equation*}
  with $\widebar{F} = F \otimes_{\ZZ} k$ and
  $\widebar{G} = G \otimes_{\ZZ} k$. The ideal
  $\calI_{n} \otimes_{\ZZ}k$ is thus perfect of grade three.
\end{prp}

\begin{prf*}
  Let us first describe the differentials in the complex
  $\calF_\bullet$ in this setting.  The last differential $\partial_3$
  is just a $(2m+3)\times 2m$ matrix with the $2m\times 2m$ block
  given by the matrix $A$ and the $3\times 2m$ block given by the
  matrix $B$.  The differential $\partial_2$ can be expressed in block
  form as
  \begin{equation*}
    \partial_2\deq
    \left(
      \begin{array}{c|c}
        A_{11}&A_{12} \\
        \hline A_{21}&A_{22}
      \end{array}
    \right)
  \end{equation*}
  where $A_{21}$ is given by multiplication by the representation
  $\bigwedge^{2m}F$ occurring in the degree $(m,0)$ component of
  $\calA$.  The matrix $A_{22}$ is given by multiplication by the
  representation $\bigwedge^{2m-1}F\otimes G$ occurring in the degree
  $(m-1,1)$ of $\calA$.  The matrix $A_{11}$ is given by
  multiplication by the representation
  $\bigwedge^{2m}F\otimes \bigwedge^2 G$ occurring in the degree
  $(m-1,2)$ component of $\calA$ and $A_{12}$ is given by
  multiplication by the representation
  $\bigwedge^{2m-1}F\otimes \bigwedge^3 G$ occurring in the degree
  $(m-2,3)$ component of $\calA$.  The relations coming from the
  second summand are three Koszul relations between the last generator
  and the three others.  The differential $\partial_1$ is given by the
  generators of $\calI_{n}$.

  The matrices of the differentials are:
  \begin{gather*}
    \partial_3 \deq
    \begin{pmatrix}
      u_{11}&u_{12}&u_{13}&\ldots&u_{1\,2m}\\
      u_{21}&u_{22}&u_{23}&\ldots&u_{2\,2m}\\
      u_{31}&u_{32}&u_{33}&\ldots&u_{3\,2m}\\
      0&c_{12}&c_{13}&\ldots&c_{1\,2m}\\
      -c_{12}&0&c_{23}&\ldots&c_{2\,2m}\\
      \vdots&\vdots&\vdots&\ddots&\vdots\\
      -c_{1\,2m}&-c_{2\,2m}&-c_{3\,2m}&\ldots&0
    \end{pmatrix} \:,
    \\[.5ex]
    \partial_2 \deq
    \begin{pmatrix}
      -x_1&-x_2&-x_3&w_1&w_2&\ldots&w_{2m}\\
      x_4&0&0&v_{\lbrace2,3\rbrace\,1}&v_{\lbrace2,3\rbrace\,2}&\ldots
      &v_{\lbrace2,3\rbrace\,2m}\\
      0&x_4&0&v_{\lbrace1,3\rbrace\,1}&-v_{\lbrace1,3\rbrace\,2}&\ldots
      &-v_{\lbrace1,3\rbrace\,2m}\\
      0&0&x_4&v_{\lbrace1,2\rbrace\,1}&v_{\lbrace1,2\rbrace\,2}&\ldots
      &v_{\lbrace1,2\rbrace\,2m}
    \end{pmatrix} \:,
  \end{gather*}
  and
  \begin{equation*}
    \partial_1 \deq
    \begin{pmatrix}
      x_1&x_2&x_3&x_4
    \end{pmatrix}
  \end{equation*}
  with entries as defined below
  \begin{gather*}
    x_1 \deq C^m, \ x_2 \deq C^{m-1}\!\wedge u_2\wedge u_3\,, \ x_3
    \deq C^{m-1}\!\wedge u_1\wedge u_3\,, \ x_4 \deq C^{m-1}\!\wedge
    u_1\wedge
    u_2 \:, \\
    v_{\lbrace \alpha,\beta\rbrace\,i} \deq \sum_j u_{\gamma
      j}\Pfbar[C]{ij}\,, \qqand w_i \deq \sum_{j,k,l}\pm
    \Delta^{j,k,l}\Pfbar[C]{ijkl} \:.
  \end{gather*}
  Here $\Delta^{j,k,l}$ is a $3\times 3$ minor of the $3\times (2m)$
  matrix $B$ on columns $j,k,l$. Finally $\gamma$ is the complement of
  $\lbrace \alpha, \beta\rbrace$ in the set $\lbrace 1,2,3\rbrace$.

  The exterior powers $C^j$ have a natural description in terms of
  Pfaffians of the matrix $A$,
  \begin{equation*}
    C^j \deq \sum_{1\le i_1<\ldots <i_{2j}\le 2m}
    \Pf[i_1\ldots i_{2j}]{A} \cdot g_{i_1}\wedge\ldots\wedge g_{i_{2j}}\;.
  \end{equation*}
  Plugging these in, we see that we get the generators of the ideals
  described in \prpref{equiv-odd} from the  matrix $\calU$ from 
  \prpref{odd-bis} via the substitution $c_{ij}=\tau_{i+3,j+3}$ and
  $u_{i,j}=\tau_{i,j+3}$.  Using this substitution we see our complex is
  just the complex described in \prpref{odd-bis}.
\end{prf*}

Notice that the representation theory dictates what the differentials
should be, as each component of $\partial_3$, $\partial_2$,
$\partial_1$ is determined by the equivariance property with respect
to $\operatorname{GL}(F)\times \operatorname{GL}(G)$ up to a non-zero
scalar.

\subsection*{Quotients of odd type}
There is a nice analogue in the odd case.  Let $n = 2m + 4$ where $m$
is a natural number.  Consider a $(2m+1)\times (2m+1)$ generic skew
symmetric matrix $A=(c_{ij})$ and a $3\times (2m+1)$ generic matrix
$B= (u_{ki})$.  Thus we work over a ring
$\calA=\operatorname{Sym}_{\ZZ}(\bigwedge^2 F\oplus F\otimes G) \is
\ZZ[c_{ij}, u_{ki}]$ where $F= \ZZ^{2m+1}$ and $G=\ZZ^3$ are free
$\ZZ$-modules.

\begin{prp}
  \label{prp:equiv-even}
  Let $\set{g_1,\ldots, g_{2m+1}}$ be a basis for $F$ and set
  \begin{equation*}
    C \deq \sum_{1\le i<j\le 2m+1} c_{ij}g_i\wedge g_j \qand
    u_{k} \deq \sum_{i=1}^{2m+1} u_{ki}g_{i}  \ \text{ for } \ 1 \le k  \le 3 \:.
  \end{equation*}
  Again we denote by $C^j$ the $j$-th exterior power of $C$ in
  $\bigwedge^{2j}F$.  The ideal
  \begin{equation*}
    \calI_{n} \deq ( C^{m-1}\wedge u_{1}\wedge u_{2}\wedge u_{3}, C^m\wedge u_{1}, C^m\wedge u_{2}, C^m\wedge u_{3})
  \end{equation*}
  is a grade $3$ almost complete intersection of type $2m+1 = n-3$.
\end{prp}

\begin{prf*}
  The exterior powers $C^j$ have a natural description in terms of
  Pfaffians of the matrix $A$.
  \begin{equation*}
    C^j \deq \sum_{1\le i_1<\ldots <i_{2j}\le 2m}
    \Pf[i_1\ldots i_{2j}]{A} \cdot g_{i_1}\wedge\ldots\wedge g_{i_{2j}}\;.
  \end{equation*}
  Plugging these in, we see that we get the generators of the ideals
  described in this proposition from the matrix $\calU$ in
  \prpref{even-bis} via the substitution $c_{ij}=\tau_{i+3,j+3}$ and
  $u_{i,j}=\tau_{i,j+3}$.
\end{prf*}
Let us work out the minimal free resolution of the ideal defined
above.

\begin{prp}
  \label{prop:odd}
  Let $n = 2m+4$ and $\calI_n$ be the ideal from
  \prpref{equiv-even}. The minimal graded free resolution of the
  cyclic $\calA$-module $\calA/\calI_n$ is
  \begin{gather*}
    \calF_\bullet: 0\lra F\otimes (\bigwedge^{2m+1}F)^{\otimes 2}
    \otimes\bigwedge^3 G\otimes \calA(-2m,-3) \xra{\partial_3} \\
    (\bigwedge^{2m+1} F)^{\otimes 2}\otimes\bigwedge^2 G\otimes \calA(-2m,-2)\oplus\bigwedge^{2m+1}F\otimes\bigwedge^{2m}F\otimes\bigwedge^3G\otimes \calA(-2m+1,-3) \\
    \xra{\partial_2} \bigwedge^{2m+1} F\otimes\bigwedge^3 G\otimes
    \calA(-m+1,-3)\oplus\bigwedge^{2m+1} F\otimes G\otimes
    \calA(-m,-1)) \xra{\partial_1} \calA \:.
  \end{gather*}
  The differentials $\partial_3$, $\partial_2$, $\partial_1$ are
  described in the proof below.  For every field $k$ the resolution
  $\calF_{\bullet}\otimes_{\ZZ}k$ is minimal over
  \begin{equation*}
    \calA_k \deq \Sym{\textstyle \bigwedge^2 \widebar{F}\oplus
      \widebar{F}\otimes \widebar{G})} \dis k[c_{ij}, u_{ki}] \:,
  \end{equation*}
  with $\widebar{F} = F \otimes_{\ZZ} k$ and
  $\widebar{G} = G \otimes_{\ZZ} k$.  The ideal
  $\calI_{n} \otimes_{\ZZ}k$ is thus perfect of grade three.
\end{prp}

\begin{prf*}
  Let us describe the differentials of the resolution in this setting.
  The last differential $\partial_3$ is just a $(2m+3)\times 2m$
  matrix with the $2m\times 2m$ block given by the matrix $A$ and the
  $3\times 2m$ block given by three vectors.  The differential
  $\partial_2$ can be expressed in block form as
  \begin{equation*}
    \partial_2\deq
    \left(
      \begin{array}{c|c}
        A_{11}&A_{12} \\
        \hline A_{21}&A_{22}
      \end{array}
    \right)
  \end{equation*}
  where $A_{21}$ is zero.  The matrix $A_{22}$ is given by
  multiplication by the representation $\bigwedge^{2m}F$ occurring in
  the degree $(m,0)$ component of $\calA$.  The matrix $A_{11}$ is
  given by multiplication by the representation
  $\bigwedge^{2m+1}F\otimes G$ occurring in the degree $(m,1)$
  component of $\calA$ and $A_{12}$ is given by multiplication by the
  representation $\bigwedge^{2m}F\otimes\bigwedge^2G$ occurring in the
  degree $(m-1,2)$ component of $\calA$.  The differential
  $\partial_1$ is given by the generators of $\calI_{n}$.

  The matrices of differentials are:
  \begin{gather*}
    \partial_3 \deq
    \begin{pmatrix}
      u_{11}&u_{12}&u_{13}&\ldots&u_{1\,2m+1}\\
      u_{21}&u_{22}&u_{23}&\ldots&u_{2\,2m+1}\\
      u_{31}&u_{32}&u_{33}&\ldots&u_{3\,2m+1}\\
      0&c_{12}&c_{13}&\ldots&c_{1\,2m+1}\\
      -c_{12}&0&c_{23}&\ldots&c_{2\,2m+1}\\
      \vdots&\vdots&\vdots&\ddots&\vdots\\
      -c_{1\,2m+1}&-c_{2\,2m+1}&-c_{3\,2m+1}&\ldots&0
    \end{pmatrix} \:,
    \\[.5ex]
    \partial_2\deq
    \begin{pmatrix}0&0&0&w_1&w_2&\ldots&w_{2m+1}\\
      x_3&x_4&0&v_{11}&v_{12}&\ldots&v_{1\,2m+1}\\
      -x_2&0&x_4&v_{21}&v_{22}&\ldots&v_{2\,2m+1}\\
      0&-x_2&-x_3&v_{31}&v_{32}&\ldots&v_{3\,2m+1}
    \end{pmatrix} \:,
  \end{gather*}
  and
  \begin{equation*}
    \partial_1 \deq
    \begin{pmatrix} x_1& x_2& x_3& x_4
    \end{pmatrix}
  \end{equation*}
  where
  \begin{gather*}
    x_1 \deq C^{m-1}\wedge u_1\wedge u_2\wedge u_3\,, \ \ x_2 \deq
    C^{m}\wedge u_1\,, \ \
    x_3 \deq C^{m}\wedge u_2\,, \ \ x_4 \deq C^{m}\wedge u_3\,,\\
    w_i \deq \pm \Pfbar[C]{i}\,, \qqand v_{\gamma i} \deq \sum_{j,k}
    \pm \Delta^{j,k}_{\alpha, \beta}\Pfbar[C]{ijk} \:.
  \end{gather*}
  Here $\Delta^{j,k}_{\alpha ,\beta}$ is a $2\times 2$ minor of the
  $3\times (2m+1)$ matrix $B$ on rows $\alpha ,\beta$ and columns
  $j,k$. Finally $\gamma$ is the complement of
  $\lbrace \alpha, \beta\rbrace$ in the set $\lbrace 1,2,3\rbrace$.

  In order to prove exactness, notice that the exterior powers $C^j$
  have a natural description in terms of Pfaffians of the matrix $A$.
  \begin{equation*}
    C^j \deq \sum_{1\le i_1<\ldots <i_{2j}\le 2m}
    \Pf[i_1\ldots i_{2j}]{A} \cdot g_{i_1}\wedge\ldots\wedge g_{i_{2j}}\;.
  \end{equation*}
  Plugging these in we see that we we get the generators of the ideals
  described in \prpref{equiv-even} from the matrix $\calU$ from 
  \prpref{even-bis}, with substitutions $c_{ij}=\tau_{i+3,j+3}$ and
  $u_{i,j}=\tau_{i,j+3}$.  Using this substitution we see our complex is
  just the complex described in  \prpref{even-bis}.
\end{prf*}

Notice that the representation theory dictates what the differentials
should be, as each component of $\partial_3$, $\partial_2$,
$\partial_1$ is determined by the equivariance property with respect
to $\operatorname{GL}(F)\times \operatorname{GL}(G)$ up to a non-zero
scalar.

  
\section[Test]{Schubert varieties in orthogonal Grassmannians vs.\
  almost complete intersection and Gorenstein ideals of codimension 3}
\label{sec:schu}

\noindent
In this section we discuss connections between the ideals described in
the previous sections and Schubert varieties in the isotropic
Grassmannian of even dimensional orthogonal space. We start with a
$2n$-dimensional vector space $\WW$ over a field $k$.  We denote by
$Q(\cdot,\cdot)$ a non-degenerate quadratic form on $\WW$ that admits
a hyperbolic basis
$\lbrace e_1,e_2,\ldots ,e_n, {\bar e}_n,\ldots ,{\bar e}_2,{\bar
  e}_1\rbrace$.  We deal with the special orthogonal group $\SO{\WW}$
of isometries of $\WW$ of determinant $1$, and its double cover
$\Spin{\WW}$.  The maximal torus $T \is (k^*)^n$ is contained in
$\SO{\WW}$ as the diagonal matrices $\underline t$ acting on $\WW$ as
follows
\begin{equation*} {\underline t}(e_i)\deq t_i e_i \qand {\underline
    t}({\bar e}_i) \deq t_i^{-1}{\bar e}_i \qtext{for} 1 \le i \le n
  \:.
\end{equation*}
We consider the lattice of integral weights for $T$, which is a free
$\ZZ$-module with coordinate basis
$\set{\epsilon_1,\ldots,\epsilon_n}$. We identify $\epsilon_i$ with
the weight of $e_i$ under this action; the weight of ${\bar e}_i$ is
$-\epsilon_i$.

There is an associated root system of type $D_n$ with roots
\begin{equation*}
  \setof{\pm\epsilon_i\pm\epsilon_j}{1\le i<j\le n} \:.
\end{equation*}
Simple roots are $\alpha_i=\epsilon_i-\epsilon_{i+1}$ for
$1\le i\le n-1$ and $\alpha_n=\epsilon_{n-1}+\epsilon_n$.  If $R(D_n)$
is a $\ZZ$-submodule of $\ZZ^n$ generated by roots, then the
\emph{fundamental weights} $\omega_i$ in the dual $\ZZ$-module
$(\ZZ^n)^*$ are generators of the dual lattice $\Lambda$, called the
\emph{weight lattice}, defined by $\omega_i(\alpha_j) = \delta_{i,j}$.
We see that
\begin{align*}
  \omega_i &\deq \epsilon_1 + \cdots + \epsilon_i\quad
             \text{for }1\le i\le n-2 \:,\\
  \omega_{n-1} &\deq \textstyle  \frac{1}{2} \sum_{i=1}^n\epsilon_i \:,
                 \text{ and} \\
  \omega_n &\deq \textstyle \frac{1}{2} \sum_{i=1}^{n-1}\epsilon_i - \frac{1}{2}\epsilon_n \:.
\end{align*}

\subsection*{The action of the Weyl group}
The Weyl group $\W{D_n}$ acts on $\Lambda$ by linear maps that permute
the roots. It is a subgroup of index 2 in a hyperoctahedral group.
$\W{D_n}$ is generated by simple reflections $s_1,s_2,\ldots
,s_{n}$. For $1\le i\le n-1$ the reflection $s_i$ simply permutes
$\epsilon_i$ and $\epsilon_{i+1}$, and $s_n$ acts as follows:
$s_n(\epsilon_i)=\epsilon_i$ for $1\le i\le n-2$,
$s_n (\epsilon_{n-1})=-\epsilon_n$, and
$s_n(\epsilon_n)=-\epsilon_{n-1}$.  It contains the permutation group
$\W{A_{n-1}}$ on $n$ elements generated by simple reflections
$s_1,\ldots ,s_{n-2}, s_n$.

Over the field of complex numbers, one can classify representations of
the group $\SO{\WW}$ and its double cover $\Spin{\WW}$. First, the
category of representations of the Spin group is semi-simple, so every
representation is a direct sum of irreducible ones. The irreducible
representations are so-called highest weight representations
$V(\lambda )$, where $\lambda=\sum_{i=1}^n\lambda_i\omega_i$ is an
integral linear combination of fundamental weights with non-negative
coefficients $\lambda_i$. The representations of $\SO{\WW}$ are direct
sums of only those irreducibles $V(\lambda)$ for which $\lambda$
written in terms of $\epsilon_i$'s belongs to $\Lambda$.  Over other
fields, and over $\ZZ$, one can define appropriate analogues of
highest weight representations.

We are interested in two particular representations: the half-spinor
representations $V(\omega_{n-1})$ and $V(\omega_n)$. They are closely
connected, as we will show, to the space of skew symmetric matrices.
To that end we recall some generalities about homogeneous spaces.

Let us work over an algebraically closed field $k$.  Let $G$ be a
reductive algebraic group and let $P_{i} \subset G$ be a parabolic
subgroup stabilizing a fundamental weight $\omega_{i} \in \Lambda$.
It is well-known that there is a canonical embedding of $G/P_{i}$ into
$\mathbb{P}(V({\omega_{i}}))$. To describe this embedding, consider
the Weyl group $W$, which naturally acts on $\Lambda$, and in it the
stabilizer $W_{\omega_{i}}$ of the $i^{\mathrm{th}}$ fundamental
weight. For each $w \in W/W_{\omega_{i}}$, let $\overset{.}{w} \in W$
be the unique minimal length representative. There is a cell
decomposition
\begin{equation*}
  G/P_{i} \deq \underset{w \in W/W_{\omega_{i}}}{\bigsqcup}
  B\overset{.}{w}P_{i}
\end{equation*}
called the Bruhat decomposition, where $B$ is the Borel subgroup
contained in $P_{i}$. The embedding
$G/P_{i} \hookrightarrow \mathbb{P}(V({\omega_{i}}))$ is given by
$b\overset{.}{w} \mapsto [b\overset{.}{w} \omega_{i}]$. In fact, we
know that $G/P = \overline{G\cdot v_{\omega_{i}}}$, where
$v_{\omega_i}$ is the highest weight vector in $V(\omega_i)$.

The cardinality of $W/W_{\omega_{i}}$ is the same as the cardinality
of the orbit $W \cdot \omega_{i}$.  Now, if the fundamental weight
$\omega_{i}$ is \textit{minuscule}, then this number coincides with
the dimension $\dim_kV({\omega_{i}})$ of the fundamental
representation. This implies that the Bruhat graph of the Bruhat
interval in the Coxeter group $(W,S)$ corresponding to the minimal
length representatives of the elements in $W/W_{\omega_{i}}$ coincides
with the \textit{crystal graph} associated to the representation
$V(\omega_{i})$.

Throughout the rest of this section we are interested in the case of a
root system of type $D_n$ and the parabolic subgroup $P_{n-1}$, the
homogeneous space $\Spin{2n}/P_{n-1}$ is one of the two connected
components of the isotropic Grassmannian $\OGr{n,2n}$.

It is well-known, see for example Laksmibai and Raghavan
\seccite[3.3]{VLkKNR08}, that the homogeneous coordinate ring of the
connected component of $\OGr{n,2n}$, considered as a projective
subvariety of the projective space $\PP(V(\omega_{n-1}))$, has a
decomposition
\begin{equation*}
  k[\OGr{n,2n}] \deq \bigoplus_{d\ge 0} V(d\omega_{n-1}) \:.
\end{equation*}
into irreducible representations of $\Spin{\WW}$, so each graded
component of this ring is irreducible.  The half-spinor representation
$V(\omega_{n-1})$ is a representation of dimension $2^{n-1}$ whose
weights with respect to the Cartan subalgebra are
$(\pm{\frac{1}{2}},\ldots ,\pm{\frac{1}{2}})$ with an even number of
minuses. It has a twin representation $V(\omega_n)$ of dimension
$2^{n-1}$ whose weights with respect to the Cartan subalgebra are
$(\pm{\frac{1}{2}},\ldots ,\pm{\frac{1}{2}})$ with an odd number of
minuses. Both half-spinor representations are constructed from the
Clifford algebra of the quadratic form $Q$.

It is also known---Kostant's Theorem, see Garfinkle's dissertation
\cite{DGr82}---that as a factor of $\Sym{V(\omega_{n-1})}$ the
coordinate ring $k[\OGr{n,2n}]$ is generated by quadratic equations.
The generators of $k[\OGr{n,2n}]$ are the spinor coordinates; they can
be indexed by the cosets $W(D_n)/W(A_{n-1})$. We denote by $q_w$ the
spinor coordinate corresponding to $w\in W(D_n)/W(A_{n-1})$; the
Schubert varieties in $\OGr{n,2n}$ are also indexed by
$W(D_n)/W(A_{n-1})$. There is a natural partial order on these
coordinates, which in the case of Schubert varieties corresponds to
the inclusion order. This partially ordered set has two combinatorial
interpretations; it is a set of $2^{n-1}$ elements.

The first interpretation of $W(D_n)/W(A_{n-1})$ is as the set
$\ESUB{n}$ of even cardinality subsets of $\set{1,2,\ldots ,n}$. The
Weyl group $W(D_n)$ acts on this set as follows. For $1\le i\le n-1$
the simple reflection $s_i$ acts by switching the numbers $i$ and
$i+1$. This means the subset is fixed by $s_i$ if it contains both or
none of the numbers $i$ and $i+1$. The reflection $s_n$ acts
non-trivially only on subsets either containing or not intersecting
the subset $\lbrace n-1,n\rbrace$. It either adds the numbers $n-1$
and $n$ or takes them away.  For a subset $I\in\ESUB{n}$ let $\ell(I)$
denote the length of a minimal representative of the corresponding
coset in $W(D_n)/W(A_{n-1})$. For a reflection $s_i$ such that
$s_i(I)\ne I$ one can prove that $\ell(I) = \ell(s_i(I))\pm 1$. The
partial order is generated by comparing $I$ and $s_i(I)$ according to
the length. In the case at hand, there is a concrete description: The
induced partial order $\ESUB{n}$ compares subsets of a given
cardinality as usual by setting
\begin{equation*}
  \set{i_1,\ldots ,i_r} \dle \set{i_1,\ldots ,i_{j-1}, i_j+1,
    i_{j+1} ,\ldots,i_r}
\end{equation*}
for $1\le i_1<i_2<\ldots <i_r$ and $i_j+1<i_{j+1}$.  The partial order
is generated by these inequalities together with the inequalities
$\lbrace i_1,\ldots,i_r\rbrace\le \lbrace i_1,\ldots
,i_r,n-1,n\rbrace$ for $1\le i_1<\ldots<i_r<n-1$; this includes
$\emptyset \le \set{n-1,n}$.

\begin{exa}
  The induced partial order on $\ESUB{4}$ is illustrated below where
  the arrows are directed such that $s_i(I) \le I$ holds.
  
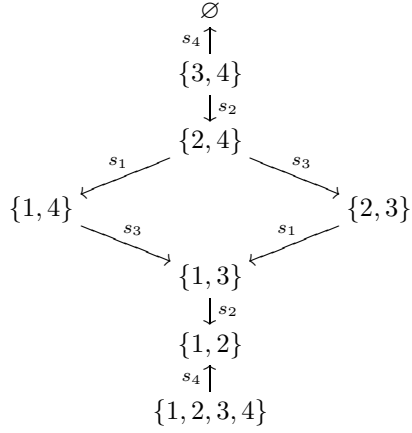
\begin{equation*}
    \xymatrix@R=.8pc{
      &\emptyset \\
      &\{3,4\} \ar[d]^-{s_2} \ar[u]^-{s_4}\\
      &\{2,4\} \ar[dl]_-{s_1} \ar[dr]^-{s_3} \\
      \{1,4\} \ar[dr]^-{s_3} && \{2,3\} \ar[dl]_-{s_1} \\
      &\{1,3\}  \ar[d]^-{s_2}\\
      &\{1,2\}  \\
      &\{1,2,3,4\} \ar[u]^-{s_4}
    }
  \end{equation*}
\end{exa}

In the second interpretation one views $W(D_n)/W(A_{n-1})$ as a
$W(D_n)$-orbit of the weight---thought of as assigning an integer to
each node of the Dynkin diagram $D_n$---under the natural action of
$W(D_n)$ on these weights.  The action of the simple reflection $s_i$
on a weight
\begin{equation*}
  \begin{matrix}
    &&&&&a_{n-1}\\
    a_1&a_2&\cdots&a_{n-3}&a_{n-2}\\
    &&&&&a_n
  \end{matrix}
\end{equation*}
changes $a_i$ to $-a_i$ and adds $a_i$ to the value at all neighboring
nodes. The partial order is generated by setting $s_i(w)>w$ if and
only if $s_i(w)\ne w$ and the node $w(i)$ is positive.

\begin{exa}
  The bijection between the set $\ESUB{4}$ and the $W(D_4)$-orbit of
  the weight
  \begin{equation*}
    w \deq \begin{smallmatrix} &&0\\0&0&\\&&1\end{smallmatrix}
  \end{equation*}
  is as follows
  \begin{alignat*}{2}
    \emptyset & \ \leftrightarrow \ w \deq \begin{smallmatrix}
      &&0\\0&0&\\&&1\end{smallmatrix} & \{3,4\} & \ \leftrightarrow
    s_4(w) \deq \begin{smallmatrix}
      &&0\\0&1&\\&&-1\end{smallmatrix} \\
    \{2,4\} & \ \leftrightarrow \ s_2s_4(w) \deq \begin{smallmatrix}
      &&1\\1&-1&\\&&0\end{smallmatrix} & \{1,4\} & \ \leftrightarrow \
    s_1s_2s_4(w) \deq \begin{smallmatrix}
      &&1\\-1&0&\\&&0\end{smallmatrix} \\
    \{2,3\} & \ \leftrightarrow \ s_3s_2s_4(w)
    \deq \begin{smallmatrix} &&-1\\1&0&\\&&0\end{smallmatrix} &
    \{1,3\} & \ \leftrightarrow \ \begin{gathered}s_3s_1s_2s_4(w)
      \\[-.5ex] s_1s_3s_2s_4(w)\end{gathered} \deq \begin{smallmatrix}
      &&-1\\-1&1&\\&&0\end{smallmatrix} \\
    \{1,2\} & \ \leftrightarrow \ s_2s_3s_1s_2s_4(w)
    \deq \begin{smallmatrix} &&0\\0&-1&\\&&1\end{smallmatrix}\qquad &
    \{1,2,3,4\} & \ \leftrightarrow \ s_4s_2s_3s_1s_2s_4(w)
    \deq \begin{smallmatrix}
      &&0\\0&0&\\&&-1\end{smallmatrix} \\
  \end{alignat*}
  Notice that this bijection commutes with the Weyl group action and
  preserves the associated partial order.
\end{exa}

\subsection*{Schubert varieties}
It is known, see for example \seccite[3.3]{VLkKNR08}, that the
defining ideal in $k[\OGr{n,2n}]$ of every Schubert variety $\Omega_w$
is generated by spinor coordinates $q_v$ for $v \not\le w$ in the
associated partial order. In our case this translates as follows.
Consider the big open cell $Y$ in $\OGr{n,2n}$ consisting of points
with Pl\"ucker coordinate $p_{\mathrm{id}}\ne 0$.  Recall the
hyperbolic basis
$\lbrace e_1,\ldots ,e_n, {\bar e}_n,\ldots , {\bar e}_1\rbrace$ of
$\WW$. To every subspace $V\in \OGr{n, 2n}$ and every basis
$\set{v_1,\ldots ,v_n}$ of $V$ we associate an $n\times 2n$ matrix $M$
whose $i^{\mathrm{th}}$ row consists of the coordinates of the vector
$v_i$ written in the basis
$\lbrace e_1,\ldots ,e_n, {\bar e}_n,\ldots , {\bar e}_1\rbrace$.  The
big open cell $Y$ in $\OGr{n, 2n}$ discussed above consists of
subspaces $V$ such that for every basis
$\lbrace v_1,\ldots ,v_n\rbrace$ of $V$ the corresponding matrix $M$
has a minor corresponding to columns $e_1,\ldots ,e_n$ not equal to
zero. The set $Y$ is an affine space of dimension $\binom{n}{2}$ as
for $V\in Y$ we can find a unique basis of $V$ such that the
corresponding matrix has a form
\begin{equation*}
  M \deq \left(
    \begin{array}{ccccc|ccccc}
      0&0&\cdots&0&1 & 0&x_{1\,2}&\cdots &x_{1\,n-1}&x_{1\,n} \\
      0&0&\cdots&1&0 & -x_{1\,2}&0&\cdots&x_{2\,n-1}&x_{2\,n}\\
      \vdots&\vdots&\ddots&\vdots&\vdots & \vdots&\vdots&\ddots&\vdots&\vdots \\
      0&1&\cdots&0&0 & -x_{1\,n-1}&-x_{2\,n-1}&\cdots &0 &x_{n-1\,n} \\
      1&0&\ldots&0&0 & -x_{1\,n}&-x_{2\,n}&\ldots&-x_{n-1\,n}&0
    \end{array}
  \right) \:.
\end{equation*}
We refer to the skew symmetric $n \times n$ block as $X$.  The
restrictions to $Y$ of the spinor coordinates correspond to
sub-Pfaffians of $X$ of all possible sizes; see for example Manivel
\cite{LMn09}. More precisely, the weights of the half-spinor
representation correspond to the subsets $I$ of the set
$\set{1,\ldots,n}$ of even cardinality. For a given $I$, the
corresponding weight is a vector
$w_I=(\pm{\frac{1}{2}},\ldots ,\pm{\frac{1}{2}})$ with $n$ coordinates
and an even number of minuses occurring in the positions determined by
$I$. The corresponding spinor coordinate $q_I$ restricts to $Y$ as the
Pfaffian of the skew symmetric matrix obtained by picking from a
generic $n\times n$ skew symmetric matrix the rows and columns
determined by $I$.  Thus, the quadratic equations generating the
defining ideal of the homogeneous coordinate ring $k[\OGr{n, 2n}]$ are
just the quadratic equations in Pfaffians of all sizes of a generic
skew symmetric matrix.

There are more facts that are known about Schubert varieties, the
reference for this is \chpcite[7]{VLkKNR08}. The half-spinor
representation $V(\omega_{n-1})$ is an example of so-called
\emph{minuscule} representation. This means all its weight vectors are
in one $W(D_n)$-orbit.  This implies that the defining ideals of
Schubert varieties and their unions behave in the optimal way
described below. For each cofinal subset $U$ of the partially ordered
set $\ESUB{n}$ we consider the ideal $J_U$ in $k[Y]$, the coordinate
ring of $Y$, generated by the spinor coordinates from that
subset. This set of ideals forms a distributive lattice $L_1$ with the
join and meet operations given by $+$ and $\cap$.  On the other hand
we can form a lattice $L_2$ of the cofinal subsets in $\ESUB{n}$ with
operations of join and meet given by $\cup$ and $\cap$.  The first
part of the next statement follows from of \seccite[7.2]{VLkKNR08};
the assertion about $J_U$ being radical follows from Brion and Kumar
\corcite[2.3.3]{MBrSKm05}.

\begin{prp}
  The lattices $L_1$ and $L_2$ are isomorphic. Moreover, the ideal
  $J_U$ is the defining ideal of the union of the Schubert varieties
  it defines set-theoretically. Thus all ideals $J_U$ are radical.
\end{prp}

Notice also that if we change the half-spinor representation
$V(\omega_{n-1})$ to the other one, i.e.\ $V(\omega_n)$ then the
lattice of Schubert varieties will change to the poset
$\mathcal{PO}_n$ of odd sized subsets of $\set{1,\ldots,n}$. The
action of the Weyl group $W(D_n)$ and the poset structure are similar
to those on $\ESUB{n}$.  We refer the reader to
\seccite[7.2]{VLkKNR08}.

To give an interpretation of the Schubert varieties of codimension 3
in concrete terms, we adopt the notation from the notes by Coskun
\cite[Lecture 5]{ICs}.  Let $\OGr{n, 2n}$ be one of the two connected
components of the orthogonal Grassmannian of $n$-dimensional isotropic
subspaces in a $2n$-dimensional vector space $\WW$.  As above denote
by $Q$ a non-degenerate quadratic form on $\WW$ that admits a
hyperbolic basis.  Fix an isotropic flag
\begin{equation*}
  F_\bullet: 0\subset F_1\subset F_2\subset \cdots \subset F_n = F_n^\perp\subset F_{n-1}^\perp
  \subset \cdots \subset F_1^\perp \subset \WW \:.
\end{equation*}
Here $F_n$ is isotropic and $F_i^\perp$ denotes the orthogonal
complement of $F_i$.  The Schubert varieties in $\OGr{n, 2n}$ are
parameterized by sequences $\lambda$,
\begin{equation*}
  n-1\ge \lambda_1>\lambda_2> \cdots >\lambda_s\ge 0 \:,
\end{equation*}
of strictly decreasing integers where $s$ has the same parity as $n$;
notice that $s\le n$ holds.  The sequence $\lambda$ determines a
unique sequence $\tilde\lambda$ of strictly decreasing integers
\begin{equation*}
  n-1\ge \tilde\lambda_{s+1}>\tilde\lambda_{s+2}>\cdots >\tilde\lambda_n\ge 0
\end{equation*}
satisfying the condition that there is no $i,j$ such that
$\lambda_i+\tilde\lambda_j=n-1$.  In other words, we obtain
$\tilde\lambda$ by removing from the sequence $n-1, n-2,\ldots ,0$ the
numbers $n-1-\lambda_1,\ldots ,n-1-\lambda_s$.

The Schubert variety $\Omega_\lambda = \Omega_\lambda (F_\bullet)$ is
defined as the closure of the locus
\begin{equation*}
  \Omega_\lambda^{(0)}(F_\bullet) \deq \left\{V\in \OGr{n,2n} \left| \;
      \begin{aligned}
        \dim_k(V\cap F_{n-\lambda_i}) &\deq i,\ \text{for}\ 1\le i\le s \\
        \dim_k(V\cap F_{\tilde\lambda_j}^\perp) &\deq j,\ \text{for}\
        s<j\le n
      \end{aligned}
    \right. \right\}\:.
\end{equation*}
The codimension of $\Omega_\lambda$ is $|\lambda | = \sum_i\lambda_i$.
The cells $\Omega_\lambda^{(0)}(F_\bullet)$ are exactly the orbits of
the Borel subgroup $B$ of the spin group $\Spin{2n}$ acting on the
connected component of $\OGr{n,2n}$.

\begin{rmk}
  In order to connect with the previous description, let us indicate
  how the partitions $\lambda$ translate to the subsets
  $\mathcal{PE}_n$ and $\mathcal {PO}_n$.  The partition
  $(\lambda_1,\ldots ,\lambda_s)$ such that
$$n-1\ge \lambda_1>\ldots >\lambda_s\ge 0$$
corresponds to the set $w(\lambda)$ of $s$ minuses in places
$n-\lambda_1, n-\lambda_2,\ldots ,n-\lambda_s$.  This set is either in
$\mathcal{PE}_n$ or $\mathcal {PO}_n$ depending on parity of $n$.
\end{rmk}

The variety $\Omega_\lambda (F_\bullet)$ in Coskun's notation is then
is then equal to the variety $\Omega_{w(\lambda)}$ in the notation of
this section.

It is well-known---see for example the works
\cite{VBMARm85,VBMVSr87,SRmARm85,ARm85} by Mehta, Ramanan, Ramanathan,
and Srinivas---that the Schubert varieties are defined over $\ZZ$ and
are normal and arithmetically Cohen-Macaulay and so are the affine
varieties $Y_\lambda = \Omega_\lambda\cap Y$. Our goal in this section
is to explicitly describe the varieties $Y_\lambda$ of codimension 3
in the affine space $Y$, i.e.\ the subvarieties $Y_\lambda$ such that
$|\lambda |=3$.

\subsection*{Spinor coordinates} Finally we describe the bijection
between the spinor coordinates and the Pfaffians of the matrix $X$.

For $n$ even and $I \in \ESUB{n}$ the corresponding spinor coordinate
$q_I$ is a square root of the minor of the matrix $M$ with columns
corresponding to $e_i$ with $i\in I$ and ${\bar e}_I$ with
$i\notin I$. This is the Pfaffian of the submatrix of $X$ obtained by
removing the rows and columns with indices $n+1-i$ for $i\in I$. The
first spinor coordinate corresponds to the subset $\emptyset$, and it
is the Pfaffian of $X$.

For $n$ odd we consider the elements $I\in \mathcal{PO}_n$. The
corresponding spinor coordinate $q_I$ is the Pfaffian of the submatrix
of $X$ obtained by removing the rows and columns with indices $n+1-i$
for $i\in I$. The first spinor coordinate corresponds to the subset
$\set{n}$, and it is the Pfaffian of the submatrix obtained from $X$
by removing the first row and column, i.e.\ the matrix
$X[\bar{1};\bar{1}]$ in the notation from \ref{notation}.

\subsection*{The case of even $n$}

  \noindent
  Let us intersect our Schubert varieties with the open cell $Y$.
  
  For $n=2$ there are evidently no Schubert varieties of codimension
  $3$, but for $n \ge 4$ there are precisely two of them, namely
  $\Omega_{(3,0)}$ and $\Omega_{(2,1)}$.  For $n=4$ it is easy to see
  they are both complete intersections, one given by
  $x_{1\,2}=x_{1\,3}=x_{1\,4}=0$ and the other by
  $x_{1\,2}=x_{1\,3}=x_{2\,3}=0$. We now assume $n\ge 6$.

  We start with the intersection $Y_{(3,0)}$.  The rank conditions
  related to the flags $F_j^\perp$ are easily seen to be empty because
  our subspace is isotropic.  The condition
  $\dim_k(V\cap F_{n-3})\ge 1$ means exactly that the rank of the
  submatrix $X[1\ldots n;4\ldots n]$ has to be less than $n-3$.  This
  is the matrix of the third differential of the almost complete
  intersection ideal of format $(1,4,n,n-3)$ described in
  \thmref{odd}.  Note that the other condition
  $\dim_k(V\cap F_n)\ge 2$ holds as the matrix $X$ has to be singular
  and therefore of rank at most $n-2$.  So the defining ideal of the
  Schubert variety $Y_{(3,0)}$ is almost complete intersection of odd
  type.

  We turn to the intersection $Y_{(2,1)}$. Again the rank conditions
  related to the flags $F_j^\perp$ are empty because our subspace is
  isotropic, so we get the conditions $\dim_k(V\cap F_{n-2})\ge 1$ and
  $\dim_k(V\cap F_{n-1})\ge 2$. The first condition is now
  redundant. The second condition means that the rank of the submatrix
  of $X[1\ldots n;2\ldots n]$ has to be at most $n-3$. This means that
  the submaximal Pfaffians $\Pfbar[X]{1i}$ of the matrix
  $X[\bar{1};\bar{1}]$ vanish for $2 \le i \le n$. It follows that the
  rank of this matrix is at most $n-4$, so the rank of
  $X[1\ldots n;2\ldots n]$ is at most $n-3$. We conclude that
  $Y_{(2,1)}$ is the subvariety given by vanishing of Pfaffians
  $\Pfbar[X]{1i}$; in other words, the defining ideal is a generic
  Gorenstein ideal of codimension $3$.

  \subsection*{The case of odd $n$}
  
  For $n\le 3$ there are evidently no Schubert varieties of
  codimension $3$, but for $n \ge 5$ there are precisely two of them,
  namely $\Omega_{(3)}$ and $\Omega_{(2,1,0)}$.  Let us intersect them
  with the open cell $Y$.

  We start with the intersection $Y_{(3)}$.  The rank conditions
  related to the flags $F_j^\perp$ are easily seen to be empty because
  our subspace is isotropic.  The condition
  $\dim_k(V\cap F_{n-3})\ge 1$ means that the rank of the submatrix
  $X[1\ldots n;4\ldots n]$ has to be less than $n-3$.  But this is the
  matrix of third differential of the almost complete intersection
  ideal of format $(1,4,n,n-3)$ described in \thmref{even}.  So the
  defining ideal of the Schubert variety $Y_{(3)}$ is an almost
  complete intersection of even type.
  
  We turn to the intersection $Y_{(2,1,0)}$. Again the rank conditions
  related to the flags $F_j^\perp$ are empty because our subspace is
  isotropic.  This leaves us with the conditions
  $\dim_k(V \cap F_{n-2})\ge 1$ and $\dim_k(V \cap F_{n-1})\ge
  2$. This just means that submaximal Pfaffians of the matrix $X$ are
  zero, so we get a generic Gorenstein ideal of codimension $3$.

  \subsection*{Minimal free resolutions}
  
  Let us look at the minimal free resolutions of the coordinate rings
  of the codimension $3$ varieties $Y_\lambda$ from the point of view
  of Schubert varieties. The defining equations of Schubert varieties
  have a general description in terms of ideals $J_U$; we recall its
  meaning in our case.  The Schubert varieties of codimension $3$
  correspond to elements
  \begin{equation*}
    w' \deq s_n s_{n-2}s_{n-1} \qqand w'' \deq s_{n-3}s_{n-2}s_{n-1} \:,
  \end{equation*}
  as these are the only two elements of length 3 in
  $W(D_n)/W(A_{n-1})$.  The generators of the corresponding ideals
  are:
  \begin{equation*}
    (q_{\mathrm{id}}, q_{s_{n-1}}, q_{s_{n-2} s_{n-1}}, q_{s_{n-3} s_{n-2} s_{n-1}},\ldots ,  q_{s_{n-i}\ldots s_{n-3} s_{n-2} s_{n-1}},\ldots ,  q_{s_{1}\ldots s_{n-3} s_{n-2} s_{n-1}}) \:,
  \end{equation*}
  where there are $n$ generators in total, and
  \begin{equation*}
    (q_{\mathrm{id}}, q_{s_{n-1}}, q_{s_{n-2} s_{n-1}}, q_{s_n s_{n-2} s_{n-1}}) \:.
  \end{equation*}
  Identifying these generators with the corresponding Pfaffians, we
  see that for $n$ odd the first ideal generated by the submaximal
  Pfaffians of $X$. The second ideal gives the almost complete
  intersection ideal described in \thmref{odd}. For even $n$, the
  first ideal gives the Pfaffian of $X$ and submaximal Pfaffians of
  $X[\bar{1};\bar{1}]$; in this case the first generator is redundant.
  The four generator ideal gives the almost complete intersection
  ideal described in \thmref{even}.
 
  There is one more statement one can make which plays an important
  role. It is be proved in terms of commutative algebra
  \lemref[Lemmas~]{regseq-odd} and \lemref[]{regseq-even}. Here we
  give a geometric reasoning proving the statement.
 
\begin{prp}
  The ideal generated by the first three spinor coordinates,
  \begin{equation*}
    (q_{\mathrm{id}}, q_{s_{n-1}}, q_{s_{n-2} s_{n-1}}) \:,
  \end{equation*}
  in the partial order on the homogeneous coordinate ring
  $k[\OGr{n,2n}]$ of the orthogonal Grassmannian is generated by a
  regular sequence. Therefore, these coordinates restricted to the
  open cell $Y$ also generate an ideal generated by a regular sequence
  in $k[Y]$. Moreover, these elements generate a radical ideal.
\end{prp}

\begin{prf*}
  Let $B$ be a Borel subgroup of the group $\Spin{2n}$. The almost
  complete intersection ideal
  $(q_{\mathrm{id}}, q_{s_{n-1}}, q_{s_{n-2} s_{n-1}})$ in
  $k[\OGr{n,2n}]$ is $B$-equivariant. This means that its vanishing
  locus is a union of Schubert cells. It follows that this vanishing
  set consists of the closure of the union of two Schubert varieties
  of codimension $3$. It is therefore an ideal of depth three
  generated by three elements. Such ideal is then generated by a
  regular sequence, as the ring $k[\OGr{n,2n}]$ is Cohen-Macaulay.
  The ideal is radical by of \corcite[2.3.3]{MBrSKm05}.  The result
  for $k[Y]$ follows by localization.
\end{prf*}
 
This result means we have an occurrence of the situation described by
Ulrich \cite{BUl90}. The ideal
$(q_{\mathrm{id}}, q_{s_{n-1}}, q_{s_{n-2} s_{n-1}})$ is the
intersection $I_{w'}\cap I_{w''}$, and thus the ideals $I_{w'}$ and
$I_{w''}$ are linked via the regular sequence
$(q_{id}, q_{s_{n-1}}, q_{s_{n-2} s_{n-1}})$. This is exactly the
procedure described in \cite{AEB87}.  By this argument, we can
describe the format of the resolutions of our almost complete
intersections.

Set $R=k[Y]$ and $n=2m+2$ for some natural number $m$. The resolution
of the Gorenstein ideal $I_{w'}$ has format
\begin{equation*}
  0 \lra R(-2m-1)\lra R^{2m+1}(-m-1)\lra R^{2m+1}(-m)\lra R \:;
\end{equation*}
we link by a regular sequence of elements of degrees $m$, $m$, and
$m+1$.  Looking at the mapping cone
\begin{equation*}
  \xymatrix@C=1pc{
    0 \ar[r] & R(-2m-1) \ar[r] & R^{2m+1}(-m-1) \ar[r] & R^{2m+1}(-m) \ar[r] & R\\
    0 \ar[r] & R(-3m-1) \ar[r] \ar[u] & R^{2}(-2m-1)\oplus
    R(-2m) \ar[r] \ar[u] & R(-m-1)\oplus R^{2}(-m) \ar[r] \ar[u] & R \ar[u]
  }
\end{equation*}
we deduce that the other ideal has a resolution with the format
\begin{equation*}
  0\lra R^{2m-1}(-2m-1)\lra R^{2m+2}(-2m)\lra R(-m-1)\oplus R^{3}(-m)\lra R
\end{equation*}

which is exactly the format of the resolution from section
\ref{sec:generic}.

Let us do this calculation for odd $n=2m+3$.  The resolution of
Gorenstein ideal of codimension 3 has format
\begin{equation*}
  0\lra R(-2m-3)\lra R^{2m+3}(-m-2)\lra R^{2m+3}(-m-1)\lra R
\end{equation*}
and we link by a regular sequence with degrees $m+1,m+1,m+1$.  Looking
at the mapping cone
\begin{equation*}
  \xymatrix{
    0 \ar[r] & R(-2m-3) \ar[r] & R^{2m+3}(-m-2) \ar[r] & R^{2m+3}(-m-1) \ar[r] & R\\
    0 \ar[r] & R(-3m-3) \ar[r] \ar[u] & R^{3}(-2m-2) \ar[r] \ar[u] &
    R^{3}(-m-1) \ar[r] \ar[u] & R \ar[u]
  }
\end{equation*}
we deduce that the other ideal has a resolution with the format
\begin{equation*}
  0\lra R^{2m}(-2m-2)\lra R^{2m+3}(-2m-1)\lra R(-m)\oplus R^{3}(-m-1)\lra R
\end{equation*}
which is exactly the format of the resolution from section
\secref{generic}.
 
Next we interpret the matrices of the free resolutions of almost
complete intersections in terms of spinor coordinates.  Before we
start, let us comment on the defining ideals of the coordinate rings
$k[\OGr{n,2n}]$ thought of as factors of the symmetric algebra on the
half-spinor representation.  By Kostant's Theorem \cite{DGr82} these
ideals are defined by quadratic equations, therefore, they are
generated by the kernel of the map
\begin{equation*}
  S_2(V(\omega_{n-1}))\lra V(2\omega_{n-1}) \:.
\end{equation*}
One has the following formula, see Adams \cite[p.~25]{JFA96},
\begin{equation*}
  S_2(V(\omega_{n-1}) \deq
  V(2\omega_{n-1})\oplus \bigoplus_{i\ge 1} V(\omega_{n-4i}) \:.
\end{equation*}
We use the notation from \secref{equiv} for the differentials in our
complexes.

Let us start with the case of odd $n=2m+3$. The generators of our
ideal in terms of Pl\"ucker coordinates are
\begin{align*}
  x_1 &\deq p_{{\bar{1}},{\bar{2}},\ldots, \widebar{2m}, 2m+1, 2m+2, 2m+3}\\
  x_2 &\deq p_{{\bar{1}},{\bar{2}},\ldots, \widebar{2m}, 2m+1, \widebar{2m+2},
        \widebar{2m+3}}\\
  x_3 &\deq p_{\bar{1},\bar{2},\ldots, \widebar{2m},\widebar{2m+1}, 2m+2, \widebar{2m+3}}\\
  x_4 &\deq p_{\bar{1},\bar{2},\ldots, \widebar{2m}, \widebar{2m+1}, \widebar{2m+2}, 2m+3} \:.
\end{align*}
The entries of the second differential $\partial_2$ are as follows.
The element $w_i$ is the Pl\"ucker coordinate with $2m+2$ bars, the
only number without bar is $2m+1-i$.  The element
$v_{\{\alpha ,\beta\}\,i}$ is a Pl\"ucker coordinate with $2m$
bars. The numbers without bars are $2m+1-i$, $2m+\alpha $,
$2m+\beta $.

The entries of the matrix $\partial_3$ are also Pl\"ucker coordinates.
The element $u_{\alpha i}$ is a Pl\"ucker coordinate with two bars, at
numbers $2m+1-i$ and $2m+4-\alpha$.  The element $c_{ij}$ is a
Pl\"ucker coordinate with two bars, at numbers $2m+1-i$ and $2m+1-j$.

The gradings of the basis vectors in the modules of the complex are:
The basis element in $F_0=R$ has weight $(0^{2m+3})$. In the following
we use $1^m$ to denote $(1,1,\ldots ,1)$ with $m$ coordinates. The
basis elements in $F_1=R\oplus R^3$ have weights
\begin{equation*}
  \textstyle
  ((\frac{1}{2})^{2m}, {-\frac{1}{2}}, {-\frac{1}{2}}, {-\frac{1}{2}}),\quad
  (({\frac{1}{2}})^{2m}, {-\frac{1}{2}}, {\frac{1}{2}}, {\frac{1}{2}}), \quad
  (({\frac{1}{2}})^{2m}, \frac{1}{2}, {-\frac{1}{2}}, {\frac{1}{2}}), \quad
  (({\frac{1}{2}})^{2m}, \frac{1}{2}, {\frac{1}{2}}, {-\frac{1}{2}}) \:.
\end{equation*}
The basis elements in $F_2=R^{2m}\oplus R^3$ have weights
\begin{gather*}
  (1^{2m-1}, 0, 0^3), \quad (1^{2m-2},0,1,0^3),\ldots ,(0, 1^{2m-1}, 0^3), \\
  (1^{2m}, 0, 0, -1), \quad (1^{2m}, 0, -1, 0), \quad (1^{2m}, -1, 0,
  0) \:.
\end{gather*}
Finally the basis vectors in $F_3=R^{2m}$ have weights
$({\frac{1}{2}}, \ldots ,{\frac{1}{2}}, \frac{3}{2}, {\frac{1}{2}},
\ldots ,{\frac{1}{2}}, ({-\frac{1}{2}})^3)$.

The composite $\partial_1\partial_2$ is easily explained. It is a
$1 \times (2m+3)$ matrix. Its entry in the $i^{\mathrm{th}}$ row is
the Pl\"ucker coordinate with the weight with $2m-1$ entries of $-1$'s
and $4$ zeros in positions $2m+1-i, 2m+1, 2m+2, 2m+3$.  This entry is
zero because it corresponds to the extremal weight vector in the
representation $V(\omega_{2m-1})$ which occurs in the
$2^{\mathrm{nd}}$ symmetric power of $V(\omega_{2m+2})$.

The composite $\partial_2 \partial_3$ is a $4\times 2m$ matrix with
the weights in the first row being
$(0,0,\ldots ,0,-1,0,\ldots,0,0,0,0)$, where $-1$ appears in positions
$1,\ldots ,2m$; the entries in the second row are
$(0,0,\ldots ,0,-1,0,\ldots,0,0,1,1)$, where $-1$ appears in positions
$1,\ldots ,2m$; the entries in the third row are
$(0,0,\ldots ,0,-1,0,\ldots,0,1,0,1)$, where $-1$ appears in positions
$1,\ldots ,2m$, finally the entries in the fourth row are
$(0,0,\ldots ,0,-1,0,\ldots,0,1,1,0)$, where $-1$ appears in positions
$1,\ldots ,2m$.

The interpretation of the identity $\partial_2 \partial_3=0$ from this
point of view requires further analysis.

Similarly we treat the even case $n=2m+4$. The generators of the ideal
$I_{w''}$ in terms of Pl\"ucker coordinates are
\begin{align*}
  x_1 &\deq p_{\bar{1},\bar{2},\ldots, \widebar{2m+4}}\\
  x_2 &\deq p_{\bar{1},\bar{2},\ldots, \widebar{2m}, \widebar{2m+1},
        \widebar{2m+2}, 2m+3, 2m+4}\\
  x_3 &\deq p_{\bar{1},\bar{2},\ldots{\widebar{2m}}, \widebar {2m+1}, 2m+2,
        \widebar{2m+3}, 2m+4}\\
  x_4 &\deq p_{\bar{1},\bar{2},\ldots, \widebar{2m}, \widebar{2m+1},
        2m+2, 2m+3,\widebar{2m+4}} \:.
\end{align*}
The entries of the second differential $\partial_2$ are as follows.
The element $w_i$ is the Pl\"ucker coordinate with $2m$ bars, the only
numbers without bar are $2m+2-i$ and $2m+2, 2m+3, 2m+4$.  The element
$v_{\alpha i}$ is a Pl\"ucker coordinate with $2m+2$ bars. The numbers
without bars are $2m+2-i$, $2m+1+\alpha $.

The entries of the matrix $\partial_3$ are also Pl\"ucker coordinates.
The element $u_{\alpha i}$ is a Pl\"ucker coordinate with two bars, at
numbers $2m+1-i$ and $2m+4-\alpha$.  The element $c_{ij}$ is a
Pl\"ucker coordinate with two bars, at numbers $2m+1-i$ and $2m+1-j$.

The gradings of the basis vectors in the modules of the complex are:
The basis element in $F_0=R$ has weight $(0^{2m+4})$.  The basis
elements in $F_1=R\oplus R^3$ have weights
\begin{equation*}
  \textstyle
  (({\frac{1}{2}})^{2m+4}), \ \  (({\frac{1}{2}})^{2m+1}, {\frac{1}{2}}, {\frac{-1}{2}}, {\frac{-1}{2}}), \ \ 
  (({\frac{1}{2}})^{2m+1}, -\frac{1}{2}, \frac{1}{2}, -\frac{1}{2}),\ \ 
  ((\frac{1}{2})^{2m+1}, -\frac{1}{2}, -\frac{1}{2}, \frac{1}{2}) \:.
\end{equation*}
The basis elements in $F_2=R^{2m}\oplus R^3$ have weights
\begin{gather*}
  (1^{2m}, 0, 0^3), (1^{2m-1},0,1,0^3),\ldots ,(0, 1^{2m}, 0^3), \\
  (1^{2m+1}, 0, 0, -1), \quad (1^{2m+1}, 0, -1, 0), \quad (1^{2m+1},
  -1, 0, 0) \:.
\end{gather*}
Finally the basis vectors in $F_3=R^{2m}$ have weights
$({\frac{1}{2}}, \ldots ,{\frac{1}{2}}, \frac{3}{2}, {\frac{1}{2}},
\ldots ,{\frac{1}{2}}, (-\frac{1}{2})^3)$.

The composite $\partial_1 \partial_2$ is easily explained. It is a
$1 \times (2m+4)$ matrix.  Its entry in the $i^{\mathrm{th}}$ row is
the Pl\"ucker coordinate with the weight with $2m$ entries of $-1$'s
and $4$ zeros in positions $2m+2-i, 2m+2, 2m+3, 2m+4$.  This entry is
zero because it corresponds to the extremal weight vector in the
representation $V(\omega_{2m})$ which occurs in the $2^{\rm nd}$
symmetric power of $V(\omega_{2m+3})$.

The composition $\partial_2 \partial_3$ is a $4\times (2m+1)$ matrix
with the weights in the first row being
$(0,0,\ldots ,0,-1,0,\ldots,0,1,1,1)$, where $-1$ appears in positions
$1,\ldots ,2m+1$; the entries in the second row are
$(0,0,\ldots ,0,-1,0,\ldots,0,1,0,0)$, where $-1$ appears in positions
$1,\ldots ,2m+1$; the third row entries are
$(0,0,\ldots ,0,-1,0,\ldots,0,0,1,0)$, where $-1$ appears in positions
$1,\ldots ,2m+1$; finally the entries in the fourth row are
$(0,0,\ldots ,0,-1,0,\ldots,0,0,0,1)$, where $-1$ appears in positions
$1,\ldots ,2m+1$.

The interpretation of the identity $\partial_2 \partial_3=0$ from this
point of view requires further analysis.

\appendix
\section{Pfaffian identities following Knuth}
\label{sec:knuth}
\noindent
For the benefit of the reader, we quote from \cite{CVW-20a} a short
introduction to Knuth's \cite{DEK96} combinatorial approach to
Pfaffians.

Let $T = (t_{ij})$ be an $n \times n$ skew symmetric matrix with
entries in a commutative ring. Assume that $T$ has zeros on the
diagonal; this is, of course, automatic if the characteristic of the
ring is not $2$. Set $\pff{ij} = t_{ij}$ for
$i,j \in \set{1,\ldots,n}$ and extend $\mathcal{P}$ to a function on
words in letters from $\set{1,\ldots,n}$ as follows:
\begin{equation*}
  \pff{i_1\ldots i_{m}} \deq
  \begin{cases}
    0 & \text{ if $m$ is odd} \\
    \sum \sgn{i_1\ldots i_{2k}}{j_1 \ldots j_{2k}}\pff{j_1j_2} \cdots
    \pff{j_{2k-1}j_{2k}} & \text{ if $m = 2k$ is even}
  \end{cases}
\end{equation*}
where the sum is over all partitions of $\set{i_1,\ldots,i_{2k}}$ in
$k$ subsets of cardinality $2$. The order of elements in each subset
is irrelevant as the difference in sign $\pff{jj'} = - \pff{j'j}$ is
offset by a change of sign of the permutation; see
\seccite[0]{DEK96}. The value of $\mathcal{P}$ on the empty word is by
convention $1$, and the value of $\mathcal{P}$ on a word with a
repeated letter is $0$. The latter is a convention in characteristic
$2$ and otherwise automatic.

The function $\mathcal{P}$ computes the Pfaffians of submatrices of
$T$. Indeed, for a subset
$\set{i_1,\ldots,i_k} \subseteq \set{1,\ldots,n}$ with elements
$i_1 < \cdots < i_k$ one has
\begin{equation*}
  \label{eq:Pff}
  \pf{i_1\ldots i_k} \deq \pff{i_1\ldots i_{k}} \:,
\end{equation*}
in the notation introduced in \ref{notation} and \eqref{notation}.

 \begin{bfhpg}[Overlapping Pfaffians]
   Let $\gra$, $\grb$, and $\grg$ be disjoint words in letters from
   $\set{1,\ldots,n}$. For $b$ a letter in $\grb$, the formula
   \cite[(5.0)]{DEK96} reads
   \begin{equation}
     \label{eq:DEK}
     \begin{aligned}
       \pff{\gra\grb}\pff{\gra\grg} &\deq \sum_{i\in\grb}
       \sgn{\grb}{bi(\grb\setminus bi)}
       \pff{\gra\grb\setminus bi}\pff{\gra\grg bi} \\
       &\hspace{3pc} + \sum_{j\in\grg} \sgn{\grb}{b(\grb\setminus
         b)}\sgn{\grg}{j(\grg\setminus j)} \pff{\gra j\grb \setminus
         b} \pff{\gra b \grg \setminus j}\;.
     \end{aligned}
   \end{equation}
   We record a number of special cases of this formula.
 
   For $\grb = b$ the formula \eqref{DEK} reduces to
   \begin{equation}
     \label{eq:DEK-b}
     \pff{\gra b}\pff{\gra\grg}
     \deq \sum_{j\in\grg}
     \sgn{\grg}{j(\grg\setminus j)}
     \pff{\gra j}
     \pff{\gra b \grg \setminus j}\;.
   \end{equation}

   For $\grg = c$ the formula \eqref{DEK} reduces to
   \begin{equation}
     \label{eq:DEK-c}
     \begin{aligned}
       \pff{\gra\grb}\pff{\gra c} &\deq \sum_{i\in\grb}
       \sgn{\grb}{bi(\grb\setminus bi)}
       \pff{\gra\grb\setminus bi}\pff{\gra cbi} \\
       &\hspace{3pc} + \sgn{\grb}{b(\grb\setminus b)} \pff{\gra c\grb
         \setminus b} \pff{\gra b} \:.
     \end{aligned}
   \end{equation}

   With $\grg$ empty the formula \eqref{DEK} reduces to
   \begin{equation}
     \label{eq:g=0}
     \pff{\gra\grb}\pff{\gra} \deq \sum_{i\in\grb}
     \sgn{\grb}{bi(\grb\setminus bi)}
     \pff{\gra\grb\setminus bi}\pff{\gra bi} \:.
   \end{equation}

   With $\gra$ and $\grg$ empty the formula \eqref{DEK} reduces to
   \begin{equation}
     \label{eq:DEK-ab=0}
     \pff{\grb} \deq \sum_{i\in\grb}
     \sgn{\grb}{bi(\grb\setminus bi)}
     \pff{\grb\setminus bi}\pff{bi} \:.
   \end{equation}
 \end{bfhpg}

 In this first appendix we derive some consequences of \eqref{DEK}
 that facilitate the computations in Appendices \secref[]{brill} and
 \secref[]{odd}. The first lemma is just the classic Laplacian
 expansion of the Pfaffian of a skew symmetric submatrix of $T$.
 
 \begin{lem}
   \label{lem:exp}
   For integers $1 \le u_1 < \cdots < u_k \le n$ and every integer
   $\ell$ with $1 \le \ell \le k$ one has
   \begin{align*}
     \Sign{\ell-1}\pf{u_1\ldots u_k} 
     &\deq \sum_{i =1}^{\ell-1} \Sign{i} t_{u_iu_\ell}\pf{u_1\ldots u_k\setminus u_iu_\ell} \\
     &\hspace{3pc} {} + \sum_{i =\ell+1}^{k} \Sign{i}
       t_{u_\ell u_i}\pf{u_1\ldots u_k\setminus u_iu_\ell} \:.
   \end{align*}
 \end{lem}

\begin{prf*}
  With $\grb = u_1\ldots u_k$ and $b=u_\ell$ the formula
  \eqref{DEK-ab=0} yields
  \begin{align*}
    \pff{\grb} & \deq \sum_{i=1}^{\ell-1}\Sign{i+\ell}
                 \pff{\grb\setminus u_\ell u_i}\pff{u_\ell u_i}
                 + \sum_{i=\ell+1}^{k}\Sign{i+\ell-1}
                 \pff{\grb\setminus u_\ell u_i}\pff{u_\ell u_i}\\
               & \deq \sum_{i=1}^{\ell-1}\Sign{i+\ell-1}
                 \pff{\grb\setminus u_\ell u_i}\pff{u_iu_\ell}
                 + \sum_{i=\ell+1}^{k}\Sign{i+\ell-1}
                 \pff{\grb\setminus u_\ell u_i}\pff{u_\ell u_i} \:. \qedhere
  \end{align*}
\end{prf*}

\begin{lem}
  \label{lem:pf1-2}
  For integers $1 \le u_1 < \cdots < u_k \le n$ and for every integer
  $\ell$ with $1 \le \ell \le k$ one has
  \begin{align*}
    \sum_{i = 1}^{\ell -1} \Sign{i} t_{ u_i u_\ell}
    &\pf{u_1\ldots u_{i-1}u_{i+1}\ldots u_k} \\
    & \deq
      \sum_{i = \ell + 1}^{n} \Sign{i} t_{ u_\ell u_i}
      \pf{u_1\ldots u_{i-1}u_{i+1}\ldots u_k} \:.
  \end{align*}
\end{lem}

\begin{prf*}
  First assume that $\ell \ge 2$ holds.  With $\gra = u_\ell$,
  $b = u_1$, and $\grg = u_2\ldots u_{k}\setminus u_\ell$ the equation
  \eqref{DEK-b} yields
  \begin{align*}
    \pff{u_\ell u_1}&\pff{u_\ell\grg} \\
    \deq & \sum_{j = 2}^{\ell-1}\Sign{j}
           \pff{u_\ell u_j}
           \pff{u_\ell u_1\grg \setminus u_j}
           + \sum_{j = \ell+1}^{k}\Sign{j-1}
           \pff{u_\ell u_j}
           \pff{u_\ell u_1\grg \setminus u_j} \:,
  \end{align*}
  which after reordering and multiplication by a sign becomes
  \begin{align*}
    &\pff{u_1u_\ell}\pff{u_2\ldots u_k} \\
    & \quad \deq \sum_{j = 2}^{\ell -1}\Sign{j}
      \pff{u_ju_\ell} \pff{u_1\ldots u_k \setminus u_j} 
      +  \sum_{j = \ell +1}^{n}\Sign{j-1}
      \pff{u_\ell u_j} \pff{u_1\ldots u_k \setminus u_j} \:,
  \end{align*}
  and that can be rewritten as
  \begin{equation*}
    \sum_{j = 1}^{\ell -1}\Sign{j}
    \pff{u_ju_\ell}    \pff{u_1\ldots u_k \setminus u_j} 
    \deq
    \sum_{j = \ell+1}^{k} \Sign{j}
    \pff{u_\ell u_j}    \pff{u_1\ldots u_k \setminus u_j} \:. 
  \end{equation*}

  Next assume that $\ell = 1$ holds.  With $\gra = u_1$, $b = u_2$,
  and $\grg = u_3\ldots u_k$ the equation \eqref{DEK-b} yields
  \begin{equation*}
    \pff{u_1 u_2}\pff{u_1\grg} \deq \sum_{j = 3}^{k}\Sign{j-1}
    \pff{u_1 u_j}
    \pff{u_1 u_2\grg \setminus u_j}
  \end{equation*}
  which can be rewritten as
  \begin{equation*}
    \sum_{j = 2}^{k}\Sign{j}
    \pff{u_1u_j}    \pff{u_1\ldots u_k \setminus u_j} 
    \deq 0 \:. \qedhere
  \end{equation*}
\end{prf*}

\begin{lem}
  \label{lem:pf3-1}
  For integers $1 \le u_1 < \cdots < u_k \le n$ and every integer
  $\ell$ with $1 \le \ell \le k$ one has,
  \begin{align*}
    \Sign{\ell-1}\pfT[T] \pfbar[T]{u_1\ldots u_k} \deq \sum_{i=1}^k \Sign{i}
    \pfbar[T]{u_iu_\ell}\pfbar[T]{u_1\ldots u_k\setminus u_iu_\ell} \:.
  \end{align*}
\end{lem}

\begin{prf*}
  With $\gra = 1\ldots n \setminus u_1\ldots u_k$,
  $\grb = u_1\ldots u_k$, and $b = u_\ell$ the formula \eqref{g=0}
  yields
  \begin{align*}
    \pff{\gra\grb}\pff{\gra}
    &\deq \sum_{i=1}^{\ell-1} \Sign{i+\ell}
      \pff{\gra\grb\setminus u_\ell u_i}\pff{\gra u_\ell u_i}  \\
    & \hspace{3pc} {} + \sum_{i=\ell+1}^{k} \Sign{i+\ell-1}
      \pff{\gra\grb\setminus u_\ell u_i}\pff{\gra u_\ell u_i} \:,
  \end{align*}
  which after reordering and multiplying by a sign becomes
  \begin{align*}
    \pff{1\ldots n}
    &\pff{1\ldots n \setminus u_1\ldots u_k}\\
    &\deq \sum_{i=1}^{k} \Sign{i+\ell-1}
      \pff{1\ldots n\setminus u_\ell u_i}
      \pff{1\ldots n \setminus (u_1\ldots u_k\setminus u_\ell u_i)} \:.
      \qedhere
  \end{align*}
\end{prf*}

\begin{lem}
  \label{lem:421-even}
  For integers $1 \le u_1 < \cdots < u_k \le n$ and every integer
  $\ell$ with $1 \le \ell \le k-1$ one has
  \begin{equation*}
    \sum_{i=1}^{\ell-1} \Sign{i}  \pfbar[T]{u_1\ldots u_{k}\setminus u_i}
    \pfbar[T]{u_iu_\ell}  \deq
    \sum_{i=\ell+1}^{k} \Sign{i} \pfbar[T]{u_1\ldots u_k\setminus u_i}
    \pfbar[T]{u_\ell u_i} \:.
  \end{equation*}
\end{lem}

\begin{prf*}
  With $\gra = 1\ldots n\setminus u_1\ldots u_k$, $b = u_k$, and
  $\grg = u_1\ldots u_{k-1}\setminus u_\ell$ formula \eqref{DEK-b}
  yields
  \begin{align*}
    \pff{\gra u_k}\pff{\gra \grg} %
    &\deq
      \sum_{j=1}^{\ell-1} \Sign{j-1} \pff{\gra u_j}\pff{\gra u_k \grg \setminus u_j} \\
    & \hspace{3pc} + \sum_{j=\ell+1}^{k-1} \Sign{j} \pff{\gra u_j}
      \pff{\gra u_k \grg \setminus u_j} \:,
  \end{align*}
  which after reordering and multiplication by a sign can be rewritten
  as
  \begin{align*}
    \Sign{k-1}\pff{1&\ldots n\setminus u_1\ldots u_{k-1}}
                      \pff{1\ldots n\setminus u_\ell u_k} \\
                    &\deq
                      \sum_{j=1}^{\ell-1} \Sign{j-1}
                      \pff{1\ldots n\setminus(u_1\ldots u_{k}\setminus u_j)}
                      \pff{1\ldots n\setminus u_ju_\ell}  \\
                    & \hspace{3pc} + \sum_{j=\ell+1}^{k-1} \Sign{j}
                      \pff{1\ldots n\setminus(u_1\ldots u_k\setminus u_j)}
                      \pff{1\ldots n\setminus u_\ell u_j} \:.
  \end{align*}
  This can also be written
  \begin{align*}
    \sum_{j=1}^{\ell-1} \Sign{j}
    & \pff{1\ldots n\setminus(u_1\ldots u_{k}\setminus u_j)}
      \pff{1\ldots n\setminus u_ju_\ell} \\
    &\deq
      \sum_{j=\ell+1}^{k} \Sign{j} \pff{1\ldots n\setminus(u_1\ldots u_k\setminus u_j)}
      \pff{1\ldots n\setminus u_\ell u_j} \:. \qedhere
  \end{align*}
\end{prf*}

\begin{lem}
  \label{lem:pf1-1}
  For integers $1 \le u_1 < \cdots < u_k \le n$ one has
  \begin{equation*}
    \sum_{i =1}^k \Sign{i} \pfbar[T]{u_i}
    \pfbar[T]{u_1\ldots u_{i-1}u_{i+1}\ldots u_k} \deq 0 \:.
  \end{equation*}
\end{lem}

\begin{prf*}  
  With $\gra = 1 \ldots n \setminus u_1\ldots u_k$, $b = u_1$, and
  $\grg = u_2\ldots u_k$ equation \eqref{DEK-b} yields
  \begin{equation*}
    \pff{\gra u_1}\pff{\gra\grg}
    \deq \sum_{j =2}^k \Sign{j}\pff{\gra u_j} \pff{\gra u_1\grg\setminus u_j} \:,
  \end{equation*}
  which after reordering and multiplication by a sign becomes
  \begin{equation*}
    \sum_{j =1}^k \Sign{j} 
    \pff{1\ldots n \setminus u_1\ldots u_{j-1}u_{j+1}\ldots u_k}\pff{1\ldots n\setminus u_j} \deq 0 \:.
    \qedhere
  \end{equation*}
\end{prf*}

\begin{lem}
  \label{lem:531-odd}
  For integers $1 \le u < v < w < x < y < z\le n$ one has
  \begin{align*}
    \pfbar[T]{y} %
    & \pfbar[T]{uvwxz} - \pfbar[T]{z}\pfbar[T]{uvwxy}  \\
    & \deq \pfbar[T]{uyz}\pfbar[T]{vwx} - \pfbar[T]{vyz}\pfbar[T]{uwx} \\
    & \hspace{3pc} {}  + \pfbar[T]{wyz}\pfbar[T]{uvx} - \pfbar[T]{xyz}\pfbar[T]{uvw}
      \:.
  \end{align*}
\end{lem}

\begin{prf*}  
  With $\gra = 1 \ldots n \setminus uvwxyz$, $\grb = uvwxy$, $b=y$,
  and $c = z$ equation \eqref{DEK-c} yields
  \begin{equation*}
    \pff{\gra\grb}\pff{\gra z} - \pff{\gra z\grb\setminus y} \pff{\gra y}
    \deq \sum_{i \in \grb}\sgn{\grb}{yi(\grb\setminus yi)}
    \pff{\gra\grb\setminus yi} \pff{\gra zyi}   \:,
  \end{equation*}
  which expands into
  \begin{align*}
    \pff{\gra uvwxy}\pff{\gra z} - \pff{\gra zuvwx} \pff{\gra y}
    & \deq \pff{\gra vwx} \pff{\gra zyu} 
      - \pff{\gra uwx} \pff{\gra zyv} \\
    & \hspace{1pc} {} +  \pff{\gra uvx} \pff{\gra zyw}
      - \pff{\gra uvw} \pff{\gra zyx} \:.
  \end{align*}
  After reordering and multiplication by $\Sign{u+v+w+x+y+z}$ it
  becomes
  \begin{align*}
    -\pff{1\ldots n \setminus z} %
    &\pff{1\ldots n \setminus uvwxy} 
      + \pff{1\ldots n \setminus y}\pff{1\ldots n\setminus uvwxz}
    \\
    & \deq \pff{1\ldots n \setminus uyz} \pff{1\ldots n \setminus vwx} \\
    & \hspace{3pc} - \pff{1\ldots n \setminus vyz} \pff{1\ldots n \setminus uwx} \\
    & \hspace{6pc} + \pff{1\ldots n \setminus wyz}
      \pff{1\ldots n\setminus uvx} \\
    & \hspace{9pc} + \pff{1\ldots n \setminus xyz}
      \pff{1\ldots n\setminus uvw}  \:. \qedhere
  \end{align*}
\end{prf*}

\begin{lem}
  \label{lem:42-even}
  For integers $1 \le u < v < w < x < y < z\le n$ one has
  \begin{align*}
    \pfbar[T]{xy}&\pfbar[T]{uvwz} - \pfbar[T]{xz}\pfbar[T]{uvwy}  
                   + \pfbar[T]{yz}\pfbar[T]{uvwx} \\
                 & \deq \pfbar[T]{uv}\pfbar[T]{wxyz} - \pfbar[T]{uw}\pfbar[T]{vxyz}
                   + \pfbar[T]{vw}\pfbar[T]{uxyz}                \:.
  \end{align*}
\end{lem}

\begin{prf*}  
  With $\gra = 1 \ldots n \setminus uvwxyz$, $\grb = uvwx$, $b=x$, and
  $\grg = yz$ equation \eqref{DEK} yields
  \begin{align*}
    \pff{\gra uvwx}\pff{\gra yz} %
    & \deq -\pff{\gra vw}\pff{\gra yzxu} + \pff{\gra uw}\pff{\gra yzxv}
      - \pff{\gra uv}\pff{\gra yzxw} \\
    & \hspace{3pc} {}  -  \pff{\gra yuvw}\pff{\gra xz} +  \pff{\gra zuvw}\pff{\gra xy}
    \\
    & \deq \pff{\gra vw}\pff{\gra uxyz} - \pff{\gra uw}\pff{\gra vxyz}
      + \pff{\gra uv}\pff{\gra wxyz} \\
    & \hspace{3pc} {}  +  \pff{\gra uvwy}\pff{\gra xz} -
      \pff{\gra uvwz}\pff{\gra xy} \:,
  \end{align*}
  which after reordering and multiplication by a sign becomes
  \begin{align*}
    \pff{1\ldots n
    &\setminus yz} %
      \pff{1\ldots n \setminus uvwx} 
      - \pff{1\ldots n \setminus xz}\pff{1\ldots n \setminus uvwy} \\
    & \hspace{3pc} {}  + \pff{1\ldots n \setminus xy}\pff{1\ldots n \setminus uvwz}
    \\
    & \deq
      \pff{1\ldots n \setminus vw}\pff{1\ldots n \setminus uxyz}                                   
      - \pff{1\ldots n \setminus uw}\pff{1\ldots n \setminus vxyz} \\
    & \hspace{3pc} {} + \pff{1\ldots n \setminus uv}\pff{1\ldots n \setminus wxyz} \:.                                       \qedhere
  \end{align*}
\end{prf*}

\begin{lem}
  \label{lem:51-odd}
  For integers $1 \le u < x < y \le n$ and $1 \le v < w < x$ one has
  \begin{align*}
    \pfbar[T]{uxy}\pfbar[T]{uvw} 
    - \pfbar[T]{u}&\pfbar[T]{uvwxy} \\
                  & \deq \pfbar[T]{uvx}\pfbar[T]{uwy}
                    -   \pfbar[T]{uwx}\pfbar[T]{uvy} \:.
  \end{align*}
\end{lem}

\begin{prf*}  
  With $\gra = 1 \ldots n \setminus uvwxy$, $\grb = vwxy$, and $b = x$
  equation \eqref{g=0} yields
  \begin{equation*}
    \pff{\gra\grb}\pff{\gra}
    \deq \sum_{i \in \grb}\sgn{\grb}{xi(\grb\setminus xi)}
    \pff{\gra\grb\setminus xi} \pff{\gra xi} \:,
  \end{equation*}
  which expands into
  \begin{equation*}
    \pff{\gra vwxy}\pff{\gra}
    \deq \pff{\gra wy} \pff{\gra xv} 
    - \pff{\gra vy} \pff{\gra xw} + \pff{\gra vw} \pff{\gra xy} \:.
  \end{equation*}
  After reordering and multiplication by $\Sign{v+w+x+y}$ this
  expression becomes
  \begin{align*}
    \pff{1\ldots n \setminus u}\pff{\gra} 
    & \deq -\pff{1\ldots n \setminus uvx} \pff{1\ldots n \setminus uwx} \\
    & \hspace{3pc} + \pff{1\ldots n \setminus uwx} \pff{1\ldots n \setminus uvy} \\
    & \hspace{6pc} + \pff{1\ldots n \setminus uxy}
      \pff{1\ldots n\setminus uvw} \:. \qedhere
  \end{align*}
\end{prf*}

\section{Minors via Pfaffians following Brill}
\label{sec:brill}

\noindent
The formula in the next theorem was first discovered by Brill
\cite{JBr04}; the theorem stated here is \thmcite[2.1]{CVW-20a}.

\begin{thm}
  \label{thm:brill}
  Let $T$ be an $n \times n$ skew symmetric matrix. Let
  $\set{i_1,\ldots,i_m}$ and $\set{j_1,\ldots,j_m}$ be subsets of
  $\set{1,\ldots,n}$ with $i_1<\cdots <i_m$ and $j_1<\cdots <j_m$, and
  set $\grr = i_1\ldots i_m$ and $\grs = j_1\ldots j_m$. The following
  equality holds:
  \begin{align*}
    \det(T[i_1&\ldots i_m;j_1\ldots j_m]) \\
              &\deq \Sign{\left\lfloor \frac{m}{2} \right\rfloor}
                \sum_{0 \le k \le \left\lfloor
                \frac{m}{2} \right\rfloor}
                \Sign{k}\sum_{\substack{\scriptscriptstyle |\omega| = 2k \\
    \scriptscriptstyle \gro \subseteq \grr}}
    \sgn{\grr}{\gro(\grr\setminus\gro)}\pff{\omega}
    \pff{(\grr\setminus\omega)\grs}\:.
  \end{align*}
\end{thm}

\noindent Notice that only subwords $\omega$ of $\rho$ that contain
$\grr \cap \grs$ contribute to the sum above.

The two lemmas proved below are applied in Appendix \secref[]{odd} to
calculate the maximal minors of the matrices $\partial_3$ from
\thmref[Theorems~]{odd} and \thmref[]{even}.

\begin{lem}
  \label{lem:d3minors-odd}
  Let $n \ge 5$ be an odd number. For integers
  $1 \le r_1 < r_2 <r_3 \le n$ one has
  \begin{align*}
    \det&\bigr(T[\overline{r_1r_2r_3};\overline{123}]\bigl)
    \\
        & \deq
          \begin{cases}
            \pfbar[T]{r_1r_2r_3}\pfbar[T]{123} & \text{if\,
              $r_2 \le 3$}
            \\[.3ex]
            \pfbar[T]{r_1r_2r_3} \pfbar[T]{123} - \pfbar[T]{123r_2r_3}
            \pfbar[T]{r_1} & \text{if\, $r_1 \le 3 < r_2$}
            \\[.3ex]
            \pfbar[T]{r_1r_2r_3} \pfbar[T]{123} -
            \pfbar[T]{23r_1r_2r_3}
            \pfbar[T]{1} \\
            \qquad {} + \pfbar[T]{13r_1r_2r_3}\pfbar[T]{2} -
            \pfbar[T]{12r_1r_2r_3} \pfbar[T]{3} & \text{if\,
              $3 < r_1 \:.$}
          \end{cases}
  \end{align*}
\end{lem}

\begin{prf*}
  Consider the words
  \begin{equation*}
    \grr \deq 1\ldots n \setminus r_1r_2r_3 \qqand \grs \deq 4\ldots n
  \end{equation*}
  of length $n-3$. One has $\grr \cap \grs = \grs\setminus r_1r_2r_3$,
  and \thmref{brill} yields
  \begin{equation}
    \label{eq:d3-brill}
    \begin{aligned}
      \det(T[&\overline{r_1r_2r_3};\overline{123}])  \\
      & \deq \Sign{\frac{n-3}{2}} \mspace{-18mu} \sum_{k = \left\lceil
          \frac{|\grs\setminus r_1r_2r_3|}{2} \right\rceil}^{
        \frac{n-3}{2}}
      \Sign{k} \mspace{-24mu} \sum_{\substack{\scriptscriptstyle |\omega| = 2k \\
          \scriptscriptstyle \grs \setminus r_1r_2r_3 \subseteq \gro
          \subseteq \grr}} \mspace{-18mu}
      \sgn{\grr}{\gro(\grr\setminus\gro)}\pff{\omega}
      \pff{(\grr\setminus\omega)\grs} \:.
    \end{aligned}
  \end{equation}

  If $r_2 \le 3$ holds, then one has $|\grr \cap \grs| = n-3$ or
  $|\grr \cap \grs| = n-4$. In either case the shortest word $\gro$
  contributing to the sum \eqref{d3-brill} has length $n-3$. Thus,
  $\gro = \grr$ is the only contributing word and one gets
  \begin{align*}
    \det(T[\overline{r_1r_2r_3};\overline{123}])
    \deq \pff{\grr}\pff{\grs} \deq \pfbar[T]{r_1r_2r_3}\pfbar[T]{123}\:.
  \end{align*}

  If $r_1 \le 3 < r_2$ holds, then one has $|\grr\cap\grs| = n-5$. As
  $n-5$ is even, the shortest subwords $\gro$ of
  \begin{equation*}
    \grr \deq (123\setminus r_1)(\grs\setminus r_2r_3)
  \end{equation*}
  that contribute to the sum \eqref{d3-brill} have length $n-5$, so
  $\gro = \grs\setminus r_2r_3$ is the only one. Now one has
  \begin{align*}
    \det(T&[\overline{r_1r_2r_3};\overline{123}])
    \\
          &  \deq \Sign{\frac{n-3}{2}}
            \sum_{k = \frac{n-5}{2}}^{\frac{n-3}{2}}
            \Sign{k} \mspace{-18mu} \sum_{\substack{\scriptscriptstyle |\omega| = 2k \\
    \scriptscriptstyle \grs \setminus r_2r_3
    \subseteq \gro \subseteq \grr}}
    \mspace{-18mu} \sgn{\grr}{\gro(\grr\setminus\gro)}\pff{\omega}
    \pff{(\grr\setminus\omega)\grs}
    \\
          & \deq \Sign{\frac{n-3}{2}} \Bigl(
            \Sign{\frac{n-5}{2}}
            \sgn{\grr}{\grs\setminus r_2r_3(123\setminus r_1)}
            \pff{\grs\setminus r_2r_3}
            \pff{(123\setminus r_1)\grs}\\
          & \hspace{6pc}  {} + \Sign{\frac{n-3}{2}}
            \pff{\grr}\pff{\grs}
            \Bigr)
    \\
          & \deq 
            -\sgn{\grr}{\grs\setminus r_2r_3(123\setminus r_1)}
            \pff{\grs\setminus r_2r_3} 
            \pff{(123\setminus r_1)\grs}  +    \pff{\grr}\pff{\grs}
    \\
          & \deq 
            -\pff{\grs\setminus r_2r_3} 
            \pff{(123\setminus r_1)\grs}  +    \pff{\grr}\pff{\grs}
    \\
          & \deq -\pfbar[T]{123r_2r_3} \pfbar[T]{r_1}
            +    \pfbar[T]{r_1r_2r_3} \pfbar[T]{123}      \:.
  \end{align*}

  If $3 < r_1$ holds, then one has $|\grr\cap\grs| = n-6$. As $n-6$ is
  odd, the shortest subwords $\gro$ of
  \begin{equation*}
    \grr \deq 123(\grs\setminus r_1r_2r_3)
  \end{equation*}
  that contribute to the sum \eqref{d3-brill} have length $n-5$. Now
  one has
  \begin{align*}
    \det&(T[\overline{r_1r_2r_3};\overline{123}])
    \\
        &  \deq \Sign{\frac{n-3}{2}}
          \sum_{k = \frac{n-5}{2}}^{\frac{n-3}{2}}
          \Sign{k} \mspace{-24mu} \sum_{\substack{\scriptscriptstyle |\omega| = 2k \\
    \scriptscriptstyle \grs \setminus r_1r_2r_3
    \subseteq \gro \subseteq \grr}}
    \mspace{-18mu}\sgn{\grr}{\gro(\grr\setminus\gro)}\pff{\omega}
    \pff{(\grr\setminus\omega)\grs}
    \\
        & \deq \Sign{\frac{n-3}{2}} \Bigl(
          \Sign{\frac{n-5}{2}}\mspace{-18mu}\sum_{\substack{\scriptscriptstyle |\omega| = n-5 \\
    \scriptscriptstyle \grs \setminus r_1r_2r_3
    \subseteq \gro \subseteq \grr}}
    \mspace{-18mu}\sgn{\grr}{\gro(\grr\setminus\gro)}\pff{\omega}
    \pff{(\grr\setminus\omega)\grs} 
    \\
        & \hspace{6pc} {} + \Sign{\frac{n-3}{2}}
          \pff{\grr}\pff{\grs}
          \Bigr)
    \\
        & \deq -\Bigl(
          \sgn{\grr}{1(\grs\setminus r_1r_2r_3)23}
          \pff{1\grs\setminus r_1r_2r_3} 
          \pff{23\grs} \\
        & \hspace{3pc} + \sgn{\grr}{2(\grs\setminus r_1r_2r_3)13}
          \pff{2\grs\setminus r_1r_2r_3} 
          \pff{13\grs}  \\
        & \hspace{6pc} +      \sgn{\grr}{3(\grs\setminus r_1r_2r_3)12}
          \pff{3\grs\setminus r_1r_2r_3} 
          \pff{12\grs}  
          \Bigr) + \pff{\grr}\pff{\grs}
    \\
        & \deq 
          \pff{1\grs\setminus r_1r_2r_3} 
          \pff{23\grs} -
          \pff{2\grs\setminus r_1r_2r_3} 
          \pff{13\grs} +
          \pff{3\grs\setminus r_1r_2r_3} 
          \pff{12\grs}  
    \\
        & \hspace{6pc} 
          {} + \pff{\grr}\pff{\grs}
    \\
        & \deq -\pfbar[T]{23r_1r_2r_3}\pfbar[T]{1} 
          +      \pfbar[T]{13r_1r_2r_3}\pfbar[T]{2} \\
        & \hspace{6pc} {} - \pfbar[T]{12r_1r_2r_3}
          \pfbar[T]{3} + \pfbar[T]{r_1r_2r_3}\pfbar[T]{123} \:. \qedhere
  \end{align*}
\end{prf*}

\begin{lem}
  \label{lem:d3minors-even}
  Let $n\ge 6$ be an even number. For integers
  $1 \le r_1 < r_2 <r_3 \le n$ one has
  \begin{align*}
    \det\bigr(T&[\overline{r_1r_2r_3};\overline{123}]\bigl)
    \\
               & \deq
                 \begin{cases}
                   0 & \text{if\, $r_3 = 3$}
                   \\[.3ex]
                   \pfbar[T]{123r_3}\pfbar[T]{r_1r_2} & \text{if\,
                     $r_2 \le 3 < r_3$}
                   \\[.3ex]
                   \pfbar[T]{12r_2r_3}\pfbar[T]{13} -
                   \pfbar[T]{13r_2r_3}\pfbar[T]{12} & \text{if\,
                     $1 = r_1 \le 3 < r_2$}
                   \\[.3ex]
                   \pfbar[T]{12r_2r_3}\pfbar[T]{23} -
                   \pfbar[T]{23r_2r_3}\pfbar[T]{12} & \text{if\,
                     $2 = r_1 \le 3 < r_2$}
                   \\[.3ex]
                   \pfbar[T]{13r_2r_3}\pfbar[T]{23} -
                   \pfbar[T]{23r_2r_3}\pfbar[T]{13} & \text{if\,
                     $3 = r_1 < r_2$}
                   \\[.3ex]
                   \pfbar[T]{1r_1r_2r_3}\pfbar[T]{23} -
                   \pfbar[T]{2r_1r_2r_3}\pfbar[T]{13} & {} \\
                   \quad+\pfbar[T]{3r_1r_2r_3}\pfbar[T]{12}
                   -\pfbar[T]{123r_1r_2r_3} \pfT[T] & \text{if\,
                     $3 < r_1 \:.$}
                 \end{cases}
  \end{align*}
\end{lem}

\begin{prf*}
  Consider the words
  \begin{equation*}
    \grr \deq 1\ldots n \setminus r_1r_2r_3 \qand \grs \deq 4\ldots n
  \end{equation*}
  of length $n-3$. One has $\grr \cap \grs = \grs\setminus r_1r_2r_3$
  and \thmref{brill} yields
  \begin{equation}
    \label{eq:d3-brill-even}
    \begin{aligned}
      \det&(T[\overline{r_1r_2r_3};\overline{123}])
      \\
      & \deq \Sign{\frac{n-4}{2}} \mspace{-18mu} \sum_{k = \left\lceil
          \frac{|\grs\setminus r_1r_2r_3|}{2} \right\rceil}^{
        \frac{n-4}{2}}
      \Sign{k}\mspace{-24mu} \sum_{\substack{\scriptscriptstyle |\omega| = 2k \\
          \scriptscriptstyle \grs \setminus r_1r_2r_3 \subseteq \gro
          \subseteq \grr}}
      \mspace{-18mu}\sgn{\grr}{\gro(\grr\setminus\gro)}\pff{\omega}
      \pff{(\grr\setminus\omega)\grs}\:.
    \end{aligned}
  \end{equation}

  If $r_3 = 3$ holds, then one has $|\grr \cap \grs| = |\grs| = n-3$,
  so the sum \eqref{d3-brill-even} is empty. i.e.
  \begin{align*}
    \det(T[\overline{123};\overline{123}])
    \deq  0\:.
  \end{align*}

  If $r_2 \le 3 < r_3$ hold, then one has $|\grr \cap \grs| = n-4$, so
  the shortest subwords $\gro$ of
  \begin{equation*}
    \grr \deq (123\setminus r_1r_2)(\grs\setminus r_3)
  \end{equation*}
  contributing to the sum \eqref{d3-brill-even} have length
  $n-4$. Thus, $\gro = \grs\setminus r_3$ is the only contributing
  word, and one gets
  \begin{align*}
    \det(T[\overline{r_1r_2r_3};\overline{123}])
    \deq \pff{\grs\setminus r_3}\pff{(123\setminus r_1r_2)\grs} =
    \pfbar[T]{123r_3}\pfbar[T]{r_1r_2} \:.
  \end{align*}

  If $r_1 \le 3 < r_2$ holds, then one has $|\grr \cap \grs| = n-5$,
  which is odd. Therefore, the shortest subwords $\gro$ of
  \begin{equation*}
    \grr \deq (123\setminus r_1)(\grs\setminus r_2r_3)
  \end{equation*}
  contributing to the sum \eqref{d3-brill-even} have length
  $n-4$. Hence, one gets
  \begin{equation}
    \label{eq:d3-minors-even}
    \begin{aligned}
      \det(T&[\overline{r_1r_2r_3};\overline{123}]) \\
      & \deq \sum_{r \in 123\setminus r_1} \sgn{\grr}{(r\grs\setminus
        r_2r_3)(123\setminus rr_1)} \pff{r\grs\setminus r_2r_3}
      \pff{(123\setminus r_1r)\grs}\:.
    \end{aligned}
  \end{equation}
  For $r_1 = 1$ this specializes to
  \begin{align*}
    \det(T[\overline{1r_2r_3};\overline{123}]) 
    & \deq \sum_{r \in 23}
      \sgn{\grr}{(r\grs\setminus r_2r_3)(23\setminus r)}
      \pff{r\grs\setminus r_2r_3}
      \pff{(23\setminus r)\grs}
    \\
    &  \deq 
      -\pfbar[T]{13r_2r_3}\pfbar[T]{12} + \pfbar[T]{12r_2r_3}\pfbar[T]{13} \:.
  \end{align*}
  The specialization of \eqref{d3-minors-even} with $r_1 = 2$ is
  \begin{align*}
    \det(T[\overline{2r_2r_3};\overline{123}])
    &  \deq \sum_{r \in 13}
      \sgn{\grr}{(r\grs\setminus r_2r_3)(13\setminus r)}
      \pff{r\grs\setminus r_2r_3}
      \pff{(13\setminus r)\grs}
    \\
    & \deq -\pfbar[T]{23r_2r_3}\pfbar[T]{12} + \pfbar[T]{12r_2r_3}\pfbar[T]{23} \:.
  \end{align*}
  The specialization of \eqref{d3-minors-even} with $r_1 = 3$ is
  \begin{align*}
    \det(T[\overline{3r_2r_3};\overline{123}])
    &  \deq \sum_{r \in 12}
      \sgn{\grr}{(r\grs\setminus r_2r_3)(12\setminus r)}
      \pff{r\grs\setminus r_2r_3}
      \pff{(12\setminus r)\grs}
    \\
    & \deq
      - \pfbar[T]{23r_2r_3}\pfbar[T]{13} + \pfbar[T]{13r_2r_3}\pfbar[T]{23} \:.
  \end{align*}
  
  If $3 < r_1$ holds, then one has $|\grr \cap \grs| = n-6$, which is
  even. Therefore, the shortest subwords $\gro$ of
  \begin{equation*}
    \grr \deq 123(\grs\setminus r_1r_2r_3)
  \end{equation*}
  that contribute to the sum \eqref{d3-brill-even} have length $n-6$,
  which means that $\gro =\grs\setminus r_1r_2r_3$ is the only
  one. Thus one has
  \begin{align*}
    \det&(T[\overline{r_1r_2r_3},\overline{123}])
    \\
        &  \deq \Sign{\frac{n-4}{2}}
          \sum_{k = \frac{n-6}{2}}^{\frac{n-4}{2}}
          \Sign{k} \mspace{-24mu} \sum_{\substack{\scriptscriptstyle |\omega| = 2k \\
    \scriptscriptstyle \grs \setminus r_1r_2r_3
    \subseteq \gro \subseteq \grr}}
    \mspace{-18mu} \sgn{\grr}{\gro(\grr\setminus\gro)}\pff{\omega}
    \pff{(\grr\setminus\omega)\grs}
    \\
        & \deq - \pff{\grs\setminus r_1r_2r_3} \pff{123\grs} +
          \sum_{\substack{\scriptscriptstyle |\omega| = n-4 \\
    \scriptscriptstyle \grs \setminus r_1r_2r_3
    \subseteq \gro \subseteq \grr}}
    \mspace{-18mu} \sgn{\grr}{\gro(\grr\setminus\gro)}\pff{\omega}
    \pff{(\grr\setminus\omega)\grs}
    \\
        & \deq - \pff{\grs\setminus r_1r_2r_3} \pff{123\grs} 
          +
          \pff{23\grs\setminus r_1r_2r_3}
          \pff{1\grs} \\
        &\hspace{6pc} -
          \pff{13\grs\setminus r_1r_2r_3}
          \pff{2\grs} + 
          \pff{12\grs\setminus r_1r_2r_3}
          \pff{3\grs}
    \\
        & \deq 
          -\pfbar[T]{123r_1r_2r_3} \pfT[T] + \pfbar[T]{1r_1r_2r_3}\pfbar[T]{23} \\
        & \hspace{6pc}       -
          \pfbar[T]{2r_1r_2r_3}\pfbar[T]{13} + 
          \pfbar[T]{3r_1r_2r_3}\pfbar[T]{12}
          \:. \qedhere
  \end{align*}
\end{prf*}


\section{Generic almost complete intersections: the proofs}
\label{sec:odd}

\noindent
In this final appendix we provide the detailed computations that
underpin the theorems in \secref{generic}.

\subsection*{Quotients of even type}

\begin{lem}
  \label{lem:setup-odd}
  Let $n \ge 5$ be an odd number and adopt the setup from
  \thmref[]{odd}. The sequence
  $0 \lra \calR^{n-3} \xra{\partial_3} \calR^{n} \xra{\partial_2}
  \calR^4 \xra{\partial_1} \calR \lra 0$ is a complex.
\end{lem}

\begin{prf*}
  The product $\partial_1\partial_2$ is a $1 \times n$ matrix; the
  first three entries are evidently $0$. For $i \in \set{4,\ldots,n}$
  the $i^{\mathrm{th}}$ entry is
  \begin{equation*}
    \pm\bigl(-\pfbar{1}\pfbar{23i} + \pfbar{2}\pfbar{13i}
    - \pfbar{3}\pfbar{12i}
    + \pfbar{123}\pfbar{i}\bigr) \:,
  \end{equation*}
  which is zero by \lemref{pf1-1} applied with $u_1\ldots u_k = 123i$.

  The product $\partial_2\partial_3$ is a $4 \times (n-3)$ matrix. Let
  $i \in \set{4,\ldots, n}$; the entry in position $(1,i-3)$ is
  \begin{equation*}
    \tau_{1i}\pfbar{123} +
    \sum_{j=4}^{i -1} \Sign{j-1}\tau_{ji}\pfbar{23j} -
    \sum_{j =i+1}^{n}\Sign{j-1}\tau_{ij}\pfbar{23j} \:.
  \end{equation*}
  Applied with $u_1\ldots u_k = 14\ldots n$ and $u_\ell = i$,
  \lemref{pf1-2} shows that this quantity is zero. Similarly,
  \lemref{pf1-2} applied with $u_1\ldots u_k = 24\ldots n$ and
  $u_\ell = i$ shows that the entry in position $(2,i-3)$ is zero, and
  an application with $u_1\ldots u_k = 3\ldots n$ and $u_\ell = i$
  shows that the entry in position $(3,i-3)$ is zero.
  The entry in position $(4,i-3)$ is
  \begin{equation*}
    \sum_{j=1}^{i-1}\Sign{j-1}\tau_{ji}\pfbar{j}
    - \sum_{j=i+1}^{n}\Sign{j-1}\tau_{i j}\pfbar{j} \:.
  \end{equation*}
  Applied with $u_1\ldots u_k = 1\ldots n$ and $u_\ell = i$,
  \lemref{pf1-2} shows that also this quantity is zero.
\end{prf*}

J\'{o}zefiak and Pragacz \cite{TJzPPr79} calculate the grade of ideals
generated by Pfaffians; we combine this with a classic result of Eagon
and Northcott \cite{JAEDGN67} to obtain the next lemma and
\lemref{regseq-even}, which deals with the case of even $n$.

\begin{lem}
  \label{lem:regseq-odd}
  Let $n\ge$ 5 be an odd number and adopt the setup from
  \thmref[]{odd}. The Pfaffians $\pfbar{1}$, $\pfbar{2}$, and
  $\pfbar{123}$ form a regular sequence in $\calR$.
\end{lem}

\begin{prf*}
  The $(n-3) \times (n-3)$ Pfaffians of the matrix
  $\calT[3\ldots n;3\ldots n]$ generate by \corcite[2.5]{TJzPPr79} an
  ideal of grade $3$ in the subring
  $\calR' = \ZZ[\grt_{ij}\mid 3 \le i < j \le n]$ of $\calR$; they are
  the Pfaffians $\pfbar{12i}$ for $3 \le i \le n$. As $\pfbar{123}$ is
  a regular element in the domain $\calR'$, the Pfaffians
  $\pfbar{12i}$ for $4 \le i \le n$ generate an ideal of grade $2$ in
  $\calS' = \calR'/\pfbar{123}$.  In $\calS = \calR/\pfbar{123}$ one
  has,
  \begin{equation*}
    \pfbar{1} \deq \sum_{i=4}^n
    \Sign{i-1}\tau_{2i}\pfbar{12i} \qand 
    \pfbar{2} \deq \sum_{i=4}^n
    \Sign{i-1}\tau_{1i}\pfbar{12i} \:.
  \end{equation*}
  Indeed, the first equality follows from \lemref{exp} applied with
  $u_1\ldots u_k = 2\ldots n$ and $\ell = 1$; the same lemma applied
  with $u_1\ldots u_k = 13\ldots n$ and $\ell =1$ yields the second
  equality.  Now it follows from \lemcite[6]{JAEDGN67} that
  $\pfbar{1}$ and $\pfbar{2}$ form a regular sequence in $\calB$.
\end{prf*}

\begin{lem}
  \label{lem:d3-minors-odd}
  Let $n \ge 5$ be an odd number and adopt the setup from
  \thmref[]{odd}. The ideal generated by the $(n-3) \times (n-3)$
  minors of the matrix $\partial_3$ contains the elements
  \begin{equation*}
    (\pfbar{1})^2\:, \quad (\pfbar{2})^2 \:,
    \quad (\pfbar{3})^2\,, \qand\, (\pfbar{123})^2 \:.
  \end{equation*}
\end{lem}

\begin{prf*}
  One has $(\pfbar{1})^2 = \det(\calT[2\ldots n; 2\ldots n])$ and
  expansion of this determinant along the first two columns, see Horn
  and Johnson \cite[0.8.9]{matrixanalysis2}, yields:
  \begin{align*}
    \det(\calT[2\ldots n; 2\ldots n]) 
    & \deq \sum_{2\le i < j \le n} \pm \det(\calT[ij;23])
      \det(\calT[\overline{1ij};\overline{123}]) \\
    & \deq \sum_{2\le i < j \le n} \pm\det(\calT[ij;23])
      \det(\partial_3[\overline{1ij};1\ldots n-3]) \:.
  \end{align*}
  Similarly, one gets
  \begin{align*}
    (\pfbar{2})^2 
    & \deq \det(\calT[13\ldots n; 13\ldots n] )
    \\
    & \deq \sum_{3\le j \le n}
      \pm\det(\calT[1j;13])
      \det(\partial_3[\overline{12j};1\ldots n-3]) \\
    & \hspace{3pc}{} + \sum_{3\le i < j \le n} \pm\det(\calT[ij;13])
      \det(\partial_3[\overline{1ij};1\ldots n-3])
  \end{align*}
  and
  \begin{align*}
    (\pfbar{3})^2 
    & \deq \det(\calT[124\ldots n; 124\ldots n] )\\
    & \deq \sum_{\substack{1 \le i < j \le n \\ i \ne 3 \ne j }}
    \pm \det(\calT[ij;12])
    \det(\partial_3[\overline{1ij};1\ldots n-3])   \:.
  \end{align*}
  Finally, one trivially has
  \begin{equation*}
    (\pfbar{123})^2 \deq \det(\calT[\overline{123};\overline{123}]) \deq
    \det(\partial_3[\overline{123};1\ldots n-3]) \:. \qedhere
  \end{equation*}
\end{prf*}

\begin{prp}
  \label{prp:minors-of-d2-odd}
  Let $n\ge 5$ be an odd number and adopt the setup from
  \thmref[]{odd}. For integers $1 \le r_1 < r_2 < r_3 \le n$ and
  $1 \le s_1 < s_2 < s_3 \le 4$ one has
  \begin{equation*}
    \det(\partial_3[\overline{r_1r_2r_3}; 1\ldots n-3])
    \det(\partial_1[1;\overline{s_1s_2s_3}])
    \deq \pm\det(\partial_2[s_1s_2s_3;r_1r_2r_3]) \:.
  \end{equation*}
\end{prp}

\begin{prf*}
  First notice that one has
  $\det(\partial_3[\overline{r_1r_2r_3}; 1\ldots n-3]) =
  \det(\calT[\overline{r_1r_2r_3};\overline{123}])$.  With the
  notation
  \begin{align*}
    \LHS &\deq \det(\calT[\overline{r_1r_2r_3};\overline{123}])
           \det(\partial_1[1;\overline{s_1s_2s_3}]) \quad\text{and} \\
    \RHS &\deq \det(\partial_2[s_1s_2s_3;r_1r_2r_3])
  \end{align*}
  the goal is to prove that $\LHS = \pm\RHS$ holds.  Set
  \begin{equation*}
    \grr \deq 1\ldots n \setminus r_1r_2r_3
    \qqand \set{s} \deq \overline{\set{s_1,s_2,s_3}} \:.
  \end{equation*}
  The possible values of $s_3$ are $3$ and $4$, and we treat these
  cases separately.

  \textbf{\em Case I.} Assuming that $s_3 = 3$ holds one has $s = 4$
  and, therefore,
  \begin{equation}
    \label{eq:d3-1}
    \det(\partial_1[1;\overline{123}]) \deq \pfbar{123} \:.
  \end{equation}
  Because the first three columns of the matrix $\partial_2$ are
  special, our argument depends on the size of the intersection
  $\set{1,2,3}\cap\set{r_1,r_2,r_3}$. We therefore consider four
  subcases determined by the (in)equalities
  \begin{equation}
    \label{eq:d3-c}
    r_3 = 3\,,\qquad r_2 \le 3 < r_3\,, \qquad r_1 \le 3 < r_2\,,
    \qqand r_1 < 3 \:.
  \end{equation}

  \textbf{\em Subcase I.a.} If $r_3 = 3$ holds, then \eqref{d3-1} and
  \lemref{d3minors-odd} yield
  \begin{equation*}
    \LHS \deq (\pfbar{123})^2\pfbar{123} \:,
  \end{equation*}
  and evidently one has $\RHS = (\pfbar{123})^3$.

  \textbf{\em Subcase I.b.} If $r_2 \le 3 < r_3$ hold, then
  \eqref{d3-1} and \lemref{d3minors-odd} yield
  \begin{equation*}
    \LHS \deq \pfbar{r_1r_2r_3}(\pfbar{123})^2 \:.
  \end{equation*}
  Expanding the determinant along the first column one has
  \begin{align*}
    \pm\RHS 
    & \deq \det
      \begin{pmatrix}
        \kd{1r_1}\pfbar{123} & 0 & \pfbar{23r_3} \\
        \kd{2r_1}\pfbar{123} & \kd{2r_2}\pfbar{123} & \pfbar{13r_3} \\
        0 & \kd{3r_2}\pfbar{123} & \pfbar{12r_3}
      \end{pmatrix}
    \\
    & \deq \kd{1r_1}\pfbar{123} \bigl(\kd{2r_2}\pfbar{123}\pfbar{12r_3}
      - \pfbar{13r_3}\kd{3r_2}\pfbar{123} \bigr)
    \\
    & \hspace{2pc} + \kd{2r_1}\pfbar{123}\kd{3r_2}\pfbar{123}\pfbar{23r_3}
    \\
    & \deq (\pfbar{123})^2\\
    & \hspace{2pc} \cdot \bigl(
      \kd{1r_1}\kd{2r_2}\pfbar{12r_3} -
      \kd{1r_1}\kd{3r_2}\pfbar{13r_3} + \kd{2r_1}\kd{3r_2}\pfbar{23r_3}
      \bigr) \:.
  \end{align*}
  For all three choices of $r_1< r_2$ in $\set{1,2,3}$ one gets
  $\RHS = \pm(\pfbar{123})^2\pfbar{r_1r_2r_3}$ as desired.

  \textbf{\em Subcase I.c.} If $r_1 \le 3 < r_2$ hold, then
  \eqref{d3-1} and \lemref{d3minors-odd} yield
  \begin{equation*}
    \LHS \deq \bigl( \pfbar{r_1r_2r_3}\pfbar{123} -\pfbar{123r_2r_3}
    \pfbar{r_1}\bigr)\pfbar{123} \:.
  \end{equation*}
  In view of \lemref{51-odd} this can be rewritten as
  \begin{align*}
    \LHS &\deq \kd{1r_1}\bigl(\pfbar{12r_2}\pfbar{13r_3} -
           \pfbar{13r_2}\pfbar{12r_3}\bigr)
           \pfbar{123}
    \\
         & \hspace{2pc} + \kd{2r_1}\bigl(\pfbar{12r_2}\pfbar{23r_3} -
           \pfbar{23r_2}\pfbar{12r_3}\bigr)\pfbar{123}
    \\
         & \hspace{4pc} + \kd{3r_1}\bigl(\pfbar{13r_2}\pfbar{23r_3} -
           \pfbar{23r_2}\pfbar{13r_3}\bigr)
           \pfbar{123} \:.
  \end{align*}
  Expansion of the determinant along the first column yields the
  matching expression
  \begin{align*}
    \pm\RHS %
    &\deq \det
      \begin{pmatrix}
        \kd{1r_1}\pfbar{123} & \pfbar{23r_2}
        & \pfbar{23r_3} \\
        \kd{2r_1}\pfbar{123} & \pfbar{13r_2} & \pfbar{13r_3} \\
        \kd{3r_1}\pfbar{123} & \pfbar{12r_2} & \pfbar{12r_3}
      \end{pmatrix}
    \\
    & \deq \kd{1r_1}\pfbar{123}\bigl(\pfbar{13r_2}\pfbar{12r_3}
      - \pfbar{13r_3}\pfbar{12r_2} \bigr)
    \\
    & \hspace{2pc} - \kd{2r_1}\pfbar{123}\bigl(\pfbar{23r_2}\pfbar{12r_3}
      - \pfbar{23r_3}\pfbar{12r_2} \bigr) \\
    & \hspace{4pc} + \kd{3r_1}\pfbar{123}\bigl(\pfbar{23r_2}\pfbar{13r_3}
      - \pfbar{23r_3}\pfbar{13r_2} \bigr) \:.            
  \end{align*}

  \textbf{\em Subcase I.d.} If $3 < r_1$ holds, then \eqref{d3-1} and
  \lemref{d3minors-odd} yield
  \begin{align*}
    \LHS & \deq \bigl( \pfbar{r_1r_2r_3}\pfbar{123}
           - \pfbar{23r_1r_2r_3} \pfbar{1} \\
         & \hspace{3pc} {}  + \pfbar{13r_1r_2r_3} \pfbar{2} 
           - \pfbar{12r_1r_2r_3} \pfbar{3} \bigr)\pfbar{123} \:.
  \end{align*}
  Expansion of the determinant along the first column yields the
  second equality in the computation below. The third equality follows
  from \lemref{51-odd} while the fifth follows from
  \lemref[Lemmas~]{pf1-1} and \lemref[]{531-odd}. Finally, the last
  equality follows from another application of \lemref{pf1-1}.
  \begin{align*}
    \pm\RHS 
    & \deq \det
      \begin{pmatrix}
        \pfbar{23r_1} & \pfbar{23r_2}
        & \pfbar{23r_3} \\
        \pfbar{13r_1} & \pfbar{13r_2} & \pfbar{13r_3} \\
        \pfbar{12r_1} & \pfbar{12r_2} & \pfbar{12r_3}
      \end{pmatrix}
    \\
    & \deq \pfbar{23r_1}\bigl( \pfbar{13r_2}\pfbar{12r_3} -
      \pfbar{13r_3}\pfbar{12r_2} \bigr) \\
    & \hspace{2pc} - \pfbar{13r_1}\bigl( \pfbar{23r_2}\pfbar{12r_3} -
      \pfbar{23r_3}\pfbar{12r_2} \bigr) \\
    & \hspace{4pc} + \pfbar{12r_1}\bigl( \pfbar{23r_2}\pfbar{13r_3} - 
      \pfbar{23r_3}\pfbar{13r_2} \bigr)
    \\
    & \deq \pfbar{23r_1}\bigl( \pfbar{123r_2r_3}\pfbar{1} -
      \pfbar{123}\pfbar{1r_2r_3} \bigr) \\
    &\hspace{2pc} - \pfbar{13r_1}\bigl(\pfbar{123r_2r_3}\pfbar{2} -
      \pfbar{123}\pfbar{2r_2r_3} \bigr) \\
    & \hspace{4pc} + \pfbar{12r_1}\bigl( \pfbar{123r_2r_3}\pfbar{3} - 
      \pfbar{123}\pfbar{3r_2r_3} \bigr)
    \\
    & \deq \pfbar{123r_2r_3}\bigl(
      \pfbar{23r_1}\pfbar{1} - \pfbar{13r_1}\pfbar{2}
      + \pfbar{12r_1}\pfbar{3}\bigr) \\
    &  \hspace{2pc} - \pfbar{123}\bigl(\pfbar{23r_1}\pfbar{1r_2r_3} \\
    &  \hspace{7pc}  - \pfbar{13r_1}\pfbar{2r_2r_3}
      + \pfbar{12r_1}\pfbar{3r_2r_3} \bigr)
    \\
    & \deq \pfbar{123r_2r_3}\pfbar{123}\pfbar{r_1} \\
    & \hspace{2pc} - \pfbar{123}\bigl( \pfbar{123r_1r_3}\pfbar{r_2} \\
    & \hspace{7pc}   - \pfbar{123r_1r_2}\pfbar{r_3} +
      \pfbar{r_1r_2r_3}\pfbar{123} \bigr)
    \\
    & \deq - \pfbar{123}\bigl( \pfbar{r_1r_2r_3}\pfbar{123} -
      \pfbar{123r_2r_3}\pfbar{r_1} \\
    & \hspace{7pc} + \pfbar{123r_1r_3}\pfbar{r_2} 
      - \pfbar{123r_1r_2}\pfbar{r_3} \bigr) 
    \\
    & \deq - \pfbar{123}\bigl( \pfbar{r_1r_2r_3}\pfbar{123} -
      \pfbar{23r_1r_2r_3}\pfbar{1} \\
    & \hspace{7pc} + \pfbar{13r_1r_2r_3}\pfbar{2} 
      - \pfbar{12r_1r_2r_3}\pfbar{3} \bigr) \:.
  \end{align*}
  Thus $\LHS = \pm\RHS$ holds, also in this subcase.

  \textbf{\em Case II.} Assuming now that $s_3 = 4$ holds, one has
  $s \in \set{1,2,3}$ and hence
  \begin{equation}
    \label{eq:d3-2}
    \det(\partial_1[1;\overline{s_1s_2s_3}]) \deq \Sign{s}\pfbar{s} \:.
  \end{equation}
  As in Case I the argument is broken into subcases following the
  (in)equalities~\eqref{d3-c}.
  
  \textbf{\em Subcase II.a.} If $r_3 = 3$ holds, then \eqref{d3-2} and
  \lemref{d3minors-odd} yield
  \begin{equation*}
    \LHS \deq \pm(\pfbar{123})^2\pfbar{s}\:,
  \end{equation*}
  and evidently one has $\RHS = \pm(\pfbar{123})^2\pfbar{s}$.

  \textbf{\em Subcase II.b.} If $r_2 \le 3 < r_3$ hold, then
  \eqref{d3-2} and \lemref{d3minors-odd} again yield
  \begin{equation*}
    \LHS \deq \pm\pfbar{r_1r_2r_3}\pfbar{123}\pfbar{s} \:.
  \end{equation*}
  This has to be compared to
  \begin{align*}
    \RHS & \deq \pm\det
           \begin{pmatrix}
             \kd{1r_1}\pfbar{123} & 0 & \pfbar{23r_3} \\
             \kd{2r_1}\pfbar{123} & \kd{2r_2}\pfbar{123}
             & \pfbar{13r_3} \\
             0 & \kd{3r_2}\pfbar{123} & \pfbar{12r_3} \\
             \Sign{r_1-1}\pfbar{r_1} & \Sign{r_2-1}\pfbar{r_2} &
             \pfbar{r_3}
           \end{pmatrix}[s_1s_24;123] \:.
  \end{align*}
  Notice that the zeros in the matrix stand for $\kd{3r_1}\pfbar{123}$
  and $\kd{1r_2}\pfbar{123}$; the determinant is thus symmetric in the
  three possible choices of $\set{s_1,s_2} \subset \set{1,2,3}$. By
  this symmetry it is sufficient to treat the choice
  $\set{s_1,s_2} = \set{1,2}$. In this case one has $s=3$ and,
  therefore,
  \begin{equation}
    \label{eq:d3-3}
    \LHS \deq \pm\pfbar{r_1r_2r_3}\pfbar{123}\pfbar{3} \:.
  \end{equation}
  Expansion of the determinant along the first column yields
  \begin{equation}
    \label{eq:d3-4}
    \begin{aligned}
      \pm\RHS & \deq \det
      \begin{pmatrix}
        \kd{1r_1}\pfbar{123} & 0 & \pfbar{23r_3} \\
        \kd{2r_1}\pfbar{123} & \kd{2r_2}\pfbar{123} &\pfbar{13r_3} \\
        \Sign{r_1-1}\pfbar{r_1} & \Sign{r_2-1}\pfbar{r_2} &
        \pfbar{r_3}
      \end{pmatrix}
      \\
      & \deq \kd{1r_1}\pfbar{123}\bigl(\kd{2r_2}
      \pfbar{123}\pfbar{r_3} + \Sign{r_2}\pfbar{13r_3}
      \pfbar{r_2}\bigr) \\
      & \hspace{2pc} + \Sign{r_2}\kd{2r_1}\pfbar{123}
      \pfbar{23r_3}\pfbar{r_2} \\
      & \hspace{4pc} + \Sign{r_1}\kd{2r_2}\pfbar{123}
      \pfbar{23r_3}\pfbar{r_1} \:.
    \end{aligned}
  \end{equation}
  For $\set{r_1,r_2} = \set{1,2}$ one has
  $\LHS = \pm\pfbar{12r_3}\pfbar{123}\pfbar{3}$. In the next
  computation, which shows that this agrees with $\pm\RHS$, the last
  equality follows from \lemref{pf1-1}.
  \begin{align*}
    \pm\RHS & \deq \pfbar{123} \bigl(
              \pfbar{123}\pfbar{r_3}  \\
            & \hspace{5pc} + \pfbar{13r_3}\pfbar{2}\bigr) 
              -\pfbar{123}\pfbar{23r_3}\pfbar{1} \bigr)
    \\
            & \deq \pfbar{123}\bigl(\pfbar{123}\pfbar{r_3}
              - \pfbar{23r_3}\pfbar{1} + \pfbar{13r_3}\pfbar{2}\bigr)
    \\
            & \deq   \pfbar{123}(-\pfbar{3}\pfbar{12r_3}) \:.
  \end{align*}
  For $\set{r_1,r_2} = \set{1,3}$ one has
  $\LHS = \pm\pfbar{13r_3}\pfbar{123}\pfbar{3}$, see \eqref{d3-3}, and
  \eqref{d3-4} specializes to the same expression.  Similarly, for
  $\set{r_1,r_2} = \set{2,3}$ one has
  $\LHS = \pm\pfbar{23r_3}\pfbar{123}\pfbar{3}$ and \eqref{d3-4}
  specializes to the same expression.

  \textbf{\em Subcase II.c.} If $r_1 \le 3 < r_2$ hold, then
  \eqref{d3-2} and \lemref{d3minors-odd} yield
  \begin{equation*}
    \LHS \deq \pm\bigr( \pfbar{r_1r_2r_3}\pfbar{123} -\pfbar{123r_2r_3}
    \pfbar{r_1}\bigr)\pfbar{s} \:.
  \end{equation*}
  This has to be compared to
  \begin{align*}
    \RHS & \deq \pm\det
           \begin{pmatrix}
             \kd{1r_1}\pfbar{123} & \pfbar{23r_2}
             & \pfbar{23r_3} \\
             \kd{2r_1}\pfbar{123} & \pfbar{13r_2} & \pfbar{13r_3} \\
             \kd{3r_1}\pfbar{123} &  \pfbar{12r_2} & \pfbar{12r_3} \\
             \Sign{r_1-1}\pfbar{r_1} & \pfbar{r_2} & \pfbar{r_3}
           \end{pmatrix}[s_1s_24;123] \:.
  \end{align*}
  This determinant is symmetric in the three possible choices of
  $\set{s_1,s_2} \subset \set{1,2,3}$. It suffices to treat the case
  $\set{s_1,s_2} = \set{1,2}$, where one has $s=3$ and, therefore,
  \begin{equation}
    \label{eq:d3-5}
    \LHS \deq \pm\bigl( \pfbar{r_1r_2r_3}\pfbar{123} -\pfbar{123r_2r_3}
    \pfbar{r_1}\bigr)\pfbar{3} \:.
  \end{equation}
  Expanding the determinant along the first column one gets
  \begin{align*}
    \pm\RHS & \deq \det
              \begin{pmatrix}
                \kd{1r_1}\pfbar{123} & \pfbar{23r_2}
                & \pfbar{23r_3} \\
                \kd{2r_1}\pfbar{123} & \pfbar{13r_2} & \pfbar{13r_3} \\
                \Sign{r_1-1}\pfbar{r_1} & \pfbar{r_2} & \pfbar{r_3}
              \end{pmatrix}
    \\
            & \deq \kd{1r_1}\pfbar{123}\bigl(\pfbar{13r_2}\pfbar{r_3}  - 
              \pfbar{13r_3}\pfbar{r_2}\bigr)\\
            & \hspace{1.9pc} - \kd{2r_1}\pfbar{123}\bigl(\pfbar{23r_2}
              \pfbar{r_3}  -  \pfbar{23r_3}\pfbar{r_2}\bigr) \\
            & \hspace{3.8pc} + \Sign{r_1-1}\pfbar{r_1}
              \bigl(\pfbar{23r_2}\pfbar{13r_3}  - 
              \pfbar{23r_3}\pfbar{13r_2}\bigr) \:.
  \end{align*}
  
  For $r_1 = 1$ this expression specializes to
  \begin{align*}
    \pm\RHS & \deq \pfbar{123} \bigl( \pfbar{13r_2}\pfbar{r_3}  - 
              \pfbar{13r_3}\pfbar{r_2} \bigr)\\
            & \hspace{5pc} + \pfbar{1} \bigl( \pfbar{23r_2}\pfbar{13r_3}  - 
              \pfbar{23r_3}\pfbar{13r_2} \bigr)
    \\
            & \deq \pfbar{13r_2}\bigl( \pfbar{123}\pfbar{r_3}
              - \pfbar{23r_3}\pfbar{1} \bigr) \\
            & \hspace{5pc}  + \pfbar{13r_3}\bigl( \pfbar{23r_2}\pfbar{1}
              - \pfbar{123}\pfbar{r_2} \bigr)
    \\
            & \deq \pfbar{13r_2}\bigl( \pfbar{13r_3}\pfbar{2}
              - \pfbar{12r_3}\pfbar{3} ) \\
            & \hspace{5pc}  + \pfbar{13r_3}( \pfbar{12r_2}\pfbar{3}
              - \pfbar{13r_2}\pfbar{2} \bigr)
    \\
            & \deq \bigl( \pfbar{13r_3}\pfbar{12r_2}
              - \pfbar{13r_2}\pfbar{12r_3} \bigr)\pfbar{3}
    \\
            & \deq \bigl(\pfbar{1r_2r_3}\pfbar{123} -\pfbar{123r_2r_3}
              \pfbar{1}\bigr)\pfbar{3} \:,
  \end{align*}
  where the third equality follows from \lemref{pf1-1} and the last
  equality holds by \lemref{51-odd}. This matches \eqref{d3-5}.

  For $r_1 = 2$ a parallel computation using the same lemmas yields
  \begin{equation*}
    \RHS \deq \pm\bigl(\pfbar{123r_2r_3}\pfbar{2}
    - \pfbar{2r_2r_3}\pfbar{123} \bigr) \pfbar{3} \:,
  \end{equation*}
  which again matches \eqref{d3-5}.


  For $r_1 = 3$ the $\RHS$ expression specializes to
  \begin{align*}
    \pm \RHS &\deq \pfbar{3} \bigl( \pfbar{23r_2}\pfbar{13r_3}  - 
               \pfbar{23r_3}\pfbar{13r_2} \bigr)
    \\
             &\deq \bigl(\pfbar{123r_2r_3}\pfbar{3}
               - \pfbar{3r_2r_3}\pfbar{123} \bigr) \pfbar{3} \:,
  \end{align*}
  where the second equality holds by \lemref{51-odd}. This matches
  \eqref{d3-5}.

  \textbf{\em Subcase II.d.} If $3 < r_1$ holds, then \eqref{d3-2} and
  \lemref{d3minors-odd} yield
  \begin{align*}
    \LHS & \deq \pm\bigl( \pfbar{r_1r_2r_3}\pfbar{123}
           - \pfbar{23r_1r_2r_3} \pfbar{1} \\
         & \hspace{3pc} {}  + \pfbar{13r_1r_2r_3} \pfbar{2} 
           - \pfbar{12r_1r_2r_3} \pfbar{3} \bigr)\pfbar{s} \:.
  \end{align*}
  This has to be compared to
  \begin{align*}
    \RHS 
    & \deq \pm\det
      \begin{pmatrix}
        \pfbar{23r_1} & \pfbar{23r_2}
        & \pfbar{23r_3} \\
        \pfbar{13r_1} & \pfbar{13r_2} & \pfbar{13r_3} \\
        \pfbar{12r_1} & \pfbar{12r_2} & \pfbar{12r_3} \\
        \pfbar{r_1} & \pfbar{r_2} & \pfbar{r_3}
      \end{pmatrix}[s_1s_24;123] \:.
  \end{align*}
  This determinant is symmetric in the three possible choices of
  $\set{s_1,s_2} \subset \set{1,2,3}$. It is sufficient to treat the
  case $\set{s_1,s_2} = \set{1,2}$, where one has $s=3$ and,
  therefore,
  \begin{equation}
    \label{eq:d3-6}
    \begin{aligned}
      \LHS & \deq \pm\bigl( \pfbar{r_1r_2r_3}\pfbar{123}
      - \pfbar{23r_1r_2r_3} \pfbar{1} \\
      & \hspace{3pc} {} + \pfbar{13r_1r_2r_3} \pfbar{2} -
      \pfbar{12r_1r_2r_3} \pfbar{3} \bigr)\pfbar{3} \:.
    \end{aligned}
  \end{equation}
  Expansion along the third row yields
  \begin{align*}
    \pm\RHS 
    & \deq \det
      \begin{pmatrix}
        \pfbar{23r_1} & \pfbar{23r_2}
        & \pfbar{23r_3} \\
        \pfbar{13r_1} & \pfbar{13r_2} & \pfbar{13r_3} \\
        \pfbar{r_1} & \pfbar{r_2} & \pfbar{r_3}
      \end{pmatrix}
    \\
    & \deq \pfbar{r_1}\bigl( \pfbar{23r_2}\pfbar{13r_3} -
      \pfbar{23r_3}\pfbar{13r_2} \bigr) \\
    & \hspace{2pc} - \pfbar{r_2}\bigl( \pfbar{23r_1}\pfbar{13r_3} -
      \pfbar{23r_3}\pfbar{13r_1} \bigr) \\
    & \hspace{4pc} + \pfbar{r_3}\bigl( \pfbar{23r_1}\pfbar{13r_2} -
      \pfbar{23r_2}\pfbar{13r_1} \bigr)
    \\
    & \deq \pfbar{r_1}\bigl( \pfbar{123r_2r_3}\pfbar{3} -
      \pfbar{3r_2r_3}\pfbar{123} \bigr) \\
    & \hspace{2pc} - \pfbar{r_2}\bigl( \pfbar{123r_1r_3}\pfbar{3} -
      \pfbar{3r_1r_3}\pfbar{123} \bigr) \\
    & \hspace{4pc} + \pfbar{r_3}\bigl( \pfbar{123r_1r_2}\pfbar{3} -
      \pfbar{3r_1r_2}\pfbar{123} \bigr)
    \\
    & \deq \bigl( \pfbar{r_1} \pfbar{123r_2r_3} \\
    & \hspace{5pc}
      -\pfbar{r_2}\pfbar{123r_1r_3} + \pfbar{r_3}
      \pfbar{123r_1r_2}\bigr)\pfbar{3}\\
    & \hspace{2pc}- \bigl( \pfbar{r_1}\pfbar{3r_2r_3} \\
    & \hspace{7pc} - \pfbar{r_2}
      \pfbar{3r_1r_3} + \pfbar{r_3}\pfbar{3r_1r_2} \bigr)\pfbar{123}
    \\
    & \deq \bigl( \pfbar{1} \pfbar{23r_1r_2r_3} 
      - \pfbar{2}\pfbar{13r_1r_2r_3} \\
    & \hspace{5pc} + \pfbar{3}\pfbar{12r_1r_2r_3}\bigr)\pfbar{3}
      - \pfbar{3}\pfbar{r_1r_2r_3}\pfbar{123}
  \end{align*}
  where the last two equalities follow from \lemref[Lemmas~]{pf1-1}
  and \lemref[]{51-odd}.
\end{prf*}

\subsection*{Quotients of odd type}
The proofs of \lemref[]{setup-even}--\prpref[]{minors-of-d2-even}
below are, if anything, slightly simpler than the proofs of
\lemref[]{setup-odd}--\prpref[]{minors-of-d2-odd}.

\begin{lem}
  \label{lem:setup-even}
  Let $n \ge 6$ be an even number and adopt the setup from
  \thmref[]{even}. The sequence
  $0 \lra \calR^{n-3} \xra{\partial_3} \calR^{n} \xra{\partial_2}
  \calR^4 \xra{\partial_1} \calR $ is a complex.
\end{lem}

\begin{prf*}
  The product $\partial_1\partial_2$ is a $1 \times n$ matrix; the
  first three entries are evidently $0$. For $i \in \set{4,\ldots,n}$
  the $i^{\mathrm{th}}$ entry is
  \begin{equation*}
    \pm\bigl(\pfT \pfbar{123i} - \pfbar{12}\pfbar{3i} +
    \pfbar{13}\pfbar{2i} - \pfbar{23}\pfbar{1i}\bigr) \;,
  \end{equation*}
  which is zero by \lemref{pf3-1} applied with $u_1\ldots u_k = 123i$
  and $u_\ell =i$.

  The product $\partial_2\partial_3$ is a $4 \times (n-3)$ matrix. Let
  $i \in \set{4,\ldots, n}$; the entry in position $(1,i-3)$ is
  \begin{equation*}
    \sum_{j=4}^{i -1} \Sign{j-1}\tau_{ji}\pfbar{123j} -
    \sum_{j =i+1}^{n}\Sign{j-1}\tau_{ij}\pfbar{123j} \:.
  \end{equation*}
  Applied with $u_1\ldots u_k = 4\ldots n$ and $u_\ell = i$,
  \lemref{pf1-2} shows that this quantity is zero.  The entry in
  position $(2,i-3)$ is
  \begin{equation*}
    \tau_{1i}\pfbar{13} - \tau_{2i}\pfbar{23} + 
    \sum_{j=4}^{i -1} \Sign{j}\tau_{ji}\pfbar{3j} -
    \sum_{j =i+1}^{n}\Sign{j}\tau_{ij}\pfbar{3j} \:.
  \end{equation*}
  Applied with $u_1\ldots u_k = 124\ldots n$ and $u_\ell = i$,
  \lemref{pf1-2} shows that this quantity is zero. Similarly,
  \lemref{pf1-2} applied with $u_1\ldots u_k = 134\ldots n$ and
  $u_\ell = i$ shows that the entry in position $(3,i-3)$ is zero, and
  an application with $u_1\ldots u_k = 2\ldots n$ and $u_\ell = i$
  shows that the entry in position $(4,i-3)$ is zero.
\end{prf*}

\begin{lem}
  \label{lem:regseq-even}
  Let $n\ge 6$ be an even number and adopt the setup from
  \thmref[]{even}. The Pfaffians $\pfbar{12}$, $\pfbar{13}$, and
  $\pfbar{23}$ form a regular sequence in $\calR$.
\end{lem}

\begin{prf*}
  The $(n-4) \times (n-4)$ Pfaffians of the matrix
  $\calT[4\ldots n;4\ldots n]$ generate by \corcite[2.5]{TJzPPr79} an
  ideal of grade $3$ in the subring
  $\calR' = \ZZ[\grt_{ij}\mid 4 \le i < j \le n]$ of $\calR$; they are
  the Pfaffians $\pfbar{123i}$ for $4 \le i \le n$.  Applied with
  $u_1\ldots u_k = 3\ldots n$ and $\ell = 1$, \lemref{exp} yields
  \begin{equation*}
    \pfbar{12} \deq \sum_{i=4}^n
    \Sign{i}\tau_{3i}\pfbar{123i} \:.
  \end{equation*}
  Similarly, applied with $u_1\ldots u_k = 24\ldots n$ and $\ell=2$
  the same lemma yields
  \begin{equation*}
    \pfbar{13} \deq \sum_{i=4}^n
    \Sign{i}\tau_{2i}\pfbar{123i} \:.
  \end{equation*}
  Finally, with $u_1\ldots u_k = 14\ldots n$ and $\ell =1$ one gets
  \begin{equation*}
    \pfbar{23} \deq \sum_{i=4}^n
    \Sign{i}\tau_{1i}\pfbar{123i} \:.
  \end{equation*}
  Now it follows from \lemcite[6]{JAEDGN67} that $\pfbar{12}$,
  $\pfbar{13}$, and $\pfbar{23}$ form a regular sequence in $\calR$.
\end{prf*}

\begin{lem}
  \label{lem:d3-minors-even}
  Let $n \ge 6$ be an even number and adopt the setup from
  \thmref[]{even}. The ideal generated by the $(n-3) \times (n-3)$
  minors of the matrix $\partial_3$ contains the elements
  \begin{equation*}
    (\pfT)^2, \quad (\pfbar{12})^2\:, \quad (\pfbar{13})^2 \:,
    \qand (\pfbar{23})^2 \:.
  \end{equation*}
\end{lem}

\begin{prf*}
  One has $(\pfT)^2 = \det(\calT)$ and expansion of this determinant
  along the first three columns, see \cite[0.8.9]{matrixanalysis2},
  yields:
  \begin{align*}
    \det(\calT) 
    & \deq \sum_{1 \le i < j < k \le n}
      \pm \det(\calT[ijk;123])
      \det(\calT[\overline{ijk};\overline{123}]) \\
    & \deq \sum_{1 \le i < j < k \le n}
      \pm \det(\calT[ijk;123])
      \det(\partial_3[\overline{ijk};1\ldots n-3]) \:.
  \end{align*}
  Similarly, expanding along the first column one gets
  \begin{align*}
    (\pfbar{12})^2 
    & \deq \det(\calT[3\ldots n; 3\ldots n]) \deq
      \sum_{i = 3}^n \pm \calT[i;3]
      \det(\partial_3[\overline{12i};1\ldots n-3]) \:,
    \\
    (\pfbar{13})^2 
    & \deq \det(\calT[24\ldots n;24\ldots n]) \deq \mspace{-10mu}
      \sum_{\substack{2 \le i \le n \\ i \ne 3}} \mspace{-8mu} \pm \calT[i;2]
    \det(\partial_3[\overline{13i};1\ldots n-3]) \:,\\[-\baselineskip]
    \intertext{and}
    (\pfbar{23})^2 
    & \deq \det(\calT[14 \ldots n; 14\ldots n])
      \deq \mspace{-10mu} \sum_{\substack{1 \le i \le n \\ i \ne 2,3}} \mspace{-8mu}
    \pm \calT[i;1]
    \det(\partial_3[\overline{23i};1\ldots n-3]) \,.\qedhere
  \end{align*}
\end{prf*}

\begin{prp}
  \label{prp:minors-of-d2-even}
  Let $n\ge 6$ be an even number and adopt the setup from
  \thmref[]{even}. For integers $1 \le r_1 < r_2 < r_3 \le n$ and
  $1 \le s_1 < s_2 < s_3 \le 4$ one has
  \begin{equation*}
    \det(\partial_3[\overline{r_1r_2r_3};1\ldots n-3])
    \det(\partial_1[1;\overline{s_1s_2s_3}])
    \deq \pm\det(\partial_2[s_1s_2s_3;r_1r_2r_3]) \:.
  \end{equation*}
\end{prp}

\begin{prf*}
  Notice that
  $\det(\partial_3[\overline{r_1r_2r_3}; 1\ldots n-3]) =
  \det(\calT[\overline{r_1r_2r_3};\overline{123}])$ holds and set
  \begin{align*}
    \LHS &\deq \det(\calT[\overline{r_1r_2r_3};\overline{123}])
           \det(\partial_1[1;\overline{s_1s_2s_3}]) \quad\text{and} \\
    \RHS &\deq \det(\partial_2[s_1s_2s_3;r_1r_2r_3]) \:.
  \end{align*}
  The goal is now to prove that $\LHS = \pm\RHS$ holds. Set
  \begin{equation*}
    \grr \deq 1\ldots n \setminus r_1r_2r_3
    \qqand \set{s} \deq \overline{\set{s_1,s_2,s_3}} \:.
  \end{equation*}
  The possible values of $s_1$ are $1$ and $2$, and we treat these
  cases separately.

  \textbf{\em Case I.} Assuming that $s_1=1$ holds, one has
  $s \in \set{2,3,4}$. By symmetry it suffices to treat the case
  $s=4$. In this case one has
  \begin{equation}
    \label{eq:d3-2e}
    \det(\partial_1[1;\overline{s_1s_2s_3}]) \deq \pfbar{23} \:.
  \end{equation}
  Because the first three columns of the matrix $\partial_2$ are
  special, our argument depends on the size of the intersection
  $\set{1,2,3}\cap\set{r_1,r_2,r_3}$. We therefore consider four
  subcases determined by the (in)equalities
  \begin{equation}
    \label{eq:d3-ce}
    r_3 = 3\,,\qquad r_2 \le 3 < r_3\,, \qquad r_1 \le 3 < r_2\,,
    \qqand r_1 < 3 \:.
  \end{equation}
 
  \textbf{\em Subcase I.a.} If $r_3 = 3$ holds, then
  \lemref{d3minors-even} yields $\LHS = 0$, and $\partial_2$ has a
  zero row, so $\RHS = 0$ holds as well.
  
  \textbf{\em Subcase I.b.} If $r_2 \le 3 < r_3$ hold, then
  \eqref{d3-2e} and \lemref{d3minors-even} yield
  \begin{equation*}
    \LHS \deq \pfbar{123r_3}\pfbar{r_1r_2}\pfbar{23} \:.
  \end{equation*}
  Expansion of the determinant along the first row yields
  \begin{align*}
    \pm\RHS 
    & \deq \det
      \begin{pmatrix}
        0 & 0 & \pfbar{123r_3}\\
        \kd{1r_1}\pfbar{13} - \kd{2r_1}\pfbar{23} & -
        \kd{2r_2}\pfbar{23}
        & -\pfbar{3r_3} \\
        -\kd{1r_1}\pfbar{12} & \kd{3r_2}\pfbar{23} & \pfbar{2r_3}
      \end{pmatrix}
    \\
    & \deq \pfbar{123r_3}\pfbar{23}\bigl(\kd{1r_1}\bigl(\kd{3r_2}\pfbar{13} 
      - \kd{2r_2}\pfbar{12}\bigr) -
      \kd{2r_1}\kd{3r_2}\pfbar{23}\bigr)
    \\
    & \deq \pm\pfbar{123r_3}\pfbar{23}\pfbar{r_1r_2} \:.
  \end{align*}

  \textbf{\em Subcase I.c.} If $r_1 \le 3 < r_2$ hold, then
  \eqref{d3-2e} and \lemref{d3minors-even} yield
  \begin{equation*}
    \LHS \deq \pfbar{23} \cdot
    \begin{cases}
      \pfbar{12r_2r_3}\pfbar{13} - \pfbar{13r_2r_3}\pfbar{12} &
      \text{if\, $r_1 = 1$}
      \\
      \pfbar{12r_2r_3}\pfbar{23} - \pfbar{23r_2r_3}\pfbar{12} &
      \text{if\, $r_1 = 2$}
      \\
      \pfbar{13r_2r_3}\pfbar{23} - \pfbar{23r_2r_3}\pfbar{13} &
      \text{if\, $r_1 = 3\:.$}
    \end{cases}
  \end{equation*}
  This has to be compared to
  \begin{align*}
    &\pm\RHS \\
    & \deq \det
      \begin{pmatrix}
        0 & \pfbar{123r_2} & \pfbar{123r_3}  \\
        \kd{1r_1}\pfbar{13} - \kd{2r_1}\pfbar{23} & -\pfbar{3r_2}  & -\pfbar{3r_3} \\
        -\kd{1r_1}\pfbar{12} + \kd{3r_1}\pfbar{23} & \pfbar{2r_2} &
        \pfbar{2r_3}
      \end{pmatrix}
    \\
    & \deq -\kd{1r_1}\bigl(\pfbar{13}\bigl(\pfbar{123r_2}\pfbar{2r_3}
      - \pfbar{123r_3}\pfbar{2r_2} \bigr)  \\
    & \hspace{3pc} +
      \pfbar{12}\bigl(-\pfbar{123r_2}\pfbar{3r_3} + \pfbar{123r_3}\pfbar{3r_2}\bigr)\bigr)
    \\
    & \hspace{2pc} + \kd{2r_1}\pfbar{23}\bigl(\pfbar{123r_2}\pfbar{2r_3}
      - \pfbar{123r_3}\pfbar{2r_2}\bigr)
    \\
    & \hspace{4pc} + \kd{3r_1}\pfbar{23}\bigl(-\pfbar{123r_2}\pfbar{3r_3}
      + \pfbar{123r_3}\pfbar{3r_2}\bigr) \:.
  \end{align*}
  \lemref{421-even} applied with $u_1\ldots u_k = 123r_2r_3$ and
  $\ell=2$ and $\ell =3$ yields
  \begin{align*}
    - \pfbar{23r_2r_3}\pfbar{12}
    & + \pfbar{12r_2r_3}\pfbar{23}  \\
    & \deq \pfbar{123r_3}\pfbar{2r_2} - \pfbar{123r_2}\pfbar{2r_3}
  \end{align*}
  and
  \begin{align*}
    - \pfbar{23r_2r_3}\pfbar{13}
    & + \pfbar{13r_2r_3}\pfbar{23}  \\
    & \deq \pfbar{123r_3}\pfbar{3r_2} - \pfbar{123r_2}\pfbar{3r_3} \:.
  \end{align*}
  The first of these identities immediately yields $\LHS = \pm\RHS$ in
  case $r_1 = 2$, and for $r_1=3$ the second identity yields the same
  conclusion. In case $r_1=1$ one applies both identities to see that
  $\LHS = \pm\RHS$ holds.
  
  \textbf{\em Subcase I.d.} If $3 < r_1$ holds, then \eqref{d3-2e} and
  \lemref{d3minors-even} yield
  \begin{align*}
    \LHS & \deq
           \bigl(\pfbar{1r_1r_2r_3}\pfbar{23} -
           \pfbar{2r_1r_2r_3}\pfbar{13} \\
         & \hspace{3pc} {} + \pfbar{3r_1r_2r_3}\pfbar{12}
           -\pfbar{123r_1r_2r_3}\pfT \bigr)\pfbar{23} \:.
  \end{align*}
  Expansion of the determinant along the third row yields the second
  equality in the computation below. The third equality follows from
  three applications of \lemref{421-even} with
  $u_1\ldots u_k = 123r_2r_3/123r_1r_3/123r_2r_2$ and $\ell=3$. The
  fifth follows from \lemref{421-even} applied with
  $u_1\ldots u_k = 23r_1r_2r_3$ and $\ell=1$ and \lemref{pf3-1}
  applied with $u_1\ldots u_k = 123r_1r_2r_3$ and $\ell = 2$.
  \begin{align*}
    \pm\RHS
    & \deq \det
      \begin{pmatrix}
        \pfbar{123r_1} & \pfbar{123r_2}
        & \pfbar{123r_3} \\
        \pfbar{3r_1} & \pfbar{3r_2} & \pfbar{3r_3} \\
        \pfbar{2r_1} & \pfbar{2r_2} & \pfbar{2r_3}
      \end{pmatrix}
    \\
    & \deq \pfbar{2r_1}\bigl( \pfbar{123r_2}\pfbar{3r_3} -
      \pfbar{123r_3}\pfbar{3r_2} \bigr) \\
    & \hspace{2pc} - \pfbar{2r_2}\bigl( \pfbar{123r_1}\pfbar{3r_3} -
      \pfbar{123r_3}\pfbar{3r_1} \bigr) \\
    & \hspace{4pc} + \pfbar{2r_3}\bigl( \pfbar{123r_1}\pfbar{3r_2} - 
      \pfbar{123r_2}\pfbar{3r_1} \bigr)
    \\
    & \deq \pfbar{2r_1}\bigl( \pfbar{23r_2r_3}\pfbar{13} -
      \pfbar{13r_2r_3}\pfbar{23} \bigr) \\
    & \hspace{2pc} - \pfbar{2r_2}\bigl( \pfbar{23r_1r_3}\pfbar{13} -
      \pfbar{13r_1r_3}\pfbar{23} \bigr) \\
    & \hspace{4pc} + \pfbar{2r_3}\bigl( \pfbar{23r_1r_2}\pfbar{13} -
      \pfbar{13r_1r_2}\pfbar{23} \bigr)
    \\
    & \deq \pfbar{13}\bigl(\pfbar{23r_2r_3}\pfbar{2r_1} \\
    & \hspace{5pc}-
      \pfbar{23r_1r_3}\pfbar{2r_2} + \pfbar{23r_1r_2}\pfbar{2r_3} \bigr) \\
    & \hspace{2pc} - \pfbar{23}\bigl(\pfbar{13r_2r_3}\pfbar{2r_1} \\
    & \hspace{7pc} -
      \pfbar{13r_1r_3}\pfbar{2r_2} + \pfbar{13r_1r_2}\pfbar{2r_3} \bigr)
    \\
    & \deq \pfbar{13}\pfbar{2r_1r_2r_3}\pfbar{23}  - \pfbar{23}\bigl(\pfbar{3r_1r_2r_3}\pfbar{12} \\
    & \hspace{2pc}  +
      \pfbar{1r_1r_2r_3}\pfbar{23} - \pfbar{123r_1r_2r_3}\pfT \bigr) \:.
  \end{align*}
  Up to a sign, this is $\LHS$.

  \textbf{\em Case II.} Assuming that $s_1=2$ holds one has $s = 1$
  and, therefore,
  \begin{equation}
    \label{eq:d3-1e}
    \det(\partial_1[1;\overline{234}]) \deq \pfT .
  \end{equation}
  As in Case I the argument is broken into subcases following the
  (in)equalities~\eqref{d3-ce}.

  \textbf{\em Subcase II.a.} If $r_3 = 3$, then \eqref{d3-1e} and
  \lemref{d3minors-even} yield $\LHS = 0$, and one has
  \begin{align*}
    \RHS & \deq \det
           \begin{pmatrix}
             \pfbar{13} & -\pfbar{23} & 0\\
             -\pfbar{12} & 0 & \pfbar{23} \\
             0 & \pfbar{12} & -\pfbar{13}
           \end{pmatrix}
    \\
         & \deq \pfbar{13}\pfbar{23}\pfbar{12}
           - \pfbar{23}\pfbar{12}\pfbar{13} \deq 0 \:.
  \end{align*}

  \textbf{\em Subcase II.b.} If $r_2 \le 3 < r_3$ hold, then
  \eqref{d3-1} and \lemref{d3minors-even} yield
  \begin{equation*}
    \LHS \deq \pfbar{123r_3}\pfbar{r_1r_2}\pfT .
  \end{equation*}
  In the computation below, the last equality follows from
  \lemref{pf3-1} applied with $u_1\ldots u_k = 123r_3$ and $\ell = 4$;
  it shows that $\LHS$ and $\RHS$ agree up to a sign.
  \begin{align*}
    &\pm\RHS \\
    & \deq \det
      \begin{pmatrix}
        \kd{1r_1}\pfbar{13} - \kd{2r_1}\pfbar{23} & -
        \kd{2r_2}\pfbar{23}
        & \pfbar{3r_3} \\
        -\kd{1r_1}\pfbar{12} & \kd{3r_2}\pfbar{23} & -\pfbar{2r_3} \\
        \kd{2r_1}\pfbar{12} & \kd{2r_2}\pfbar{12} -
        \kd{3r_2}\pfbar{13} & \pfbar{1r_3}
      \end{pmatrix}
    \\
    & \deq\kd{1r_1}\kd{2r_2}\pfbar{12}\bigl(\pfbar{13}\pfbar{2r_3} -
      \pfbar{23}\pfbar{1r_3} - \pfbar{12}\pfbar{3r_3} \bigr)
    \\
    & \hspace{1pc} {} + \kd{1r_1}\kd{3r_2}\pfbar{13}\bigl(\pfbar{23}\pfbar{1r_3}
      - \pfbar{13}\pfbar{2r_3} + \pfbar{12}\pfbar{3r_3} \bigr)
    \\
    &  \hspace{2pc} {} - 
      \kd{2r_1}\kd{3r_2}\pfbar{23}\bigl(\pfbar{23}\pfbar{1r_3}
      - \pfbar{13}\pfbar{2r_3} + \pfbar{12}\pfbar{3r_3} \bigr)
    \\
    & \deq \pm \pfbar{r_1r_2}\bigl(\pfbar{23}\pfbar{1r_3} - \pfbar{13}\pfbar{2r_3} +
      \pfbar{12}\pfbar{3r_3} \bigr)
    \\ 
    & \deq \pfbar{r_1r_2}\pfT\pfbar{123r_3} \:.
  \end{align*}

  \textbf{\em Subcase II.c.} If $r_1 \le 3 < r_2$ hold, then
  \eqref{d3-1e} and \lemref{d3minors-even} yield
  \begin{equation*}
    \LHS \deq \pfT \cdot
    \begin{cases}
      \pfbar{12r_2r_3}\pfbar{13} - \pfbar{13r_2r_3}\pfbar{12} &
      \text{if\, $r_1 = 1$}
      \\
      \pfbar{12r_2r_3}\pfbar{23} - \pfbar{23r_2r_3}\pfbar{12} &
      \text{if\, $r_1 = 2$}
      \\
      \pfbar{13r_2r_3}\pfbar{23} - \pfbar{23r_2r_3}\pfbar{13} &
      \text{if\, $r_1 = 3\:.$}
    \end{cases}
  \end{equation*}
  This has to be compared to
  \begin{align*}
    \pm\RHS 
    & \deq \det
      \begin{pmatrix}
        \kd{1r_1}\pfbar{13} - \kd{2r_1}\pfbar{23} & \pfbar{3r_2} & \pfbar{3r_3}  \\
        -\kd{1r_1}\pfbar{12} + \kd{3r_1}\pfbar{23} & -\pfbar{2r_2}  & -\pfbar{2r_3} \\
        \kd{2r_1}\pfbar{12} - \kd{3r_1}\pfbar{13} & \pfbar{1r_2} &
        \pfbar{1r_3}
      \end{pmatrix}
    \\
    & \deq\kd{1r_1}\bigl(\pfbar{13}\bigl(-\pfbar{2r_2}\pfbar{1r_3} + \pfbar{1r_2}\pfbar{2r_3} \bigr)\\
    & \hspace{6pc} {} +
      \pfbar{12}\bigl(\pfbar{3r_2}\pfbar{1r_3} - \pfbar{1r_2}\pfbar{3r_3}\bigr)\bigr)
    \\
    & \:-\: \kd{2r_1}\bigl(\pfbar{23}\bigl(-\pfbar{2r_2}\pfbar{1r_3} + \pfbar{1r_2}\pfbar{2r_3}\bigr) \\
    & \hspace{6pc} {} -
      \pfbar{12}\bigl(-\pfbar{3r_2}\pfbar{2r_3} + \pfbar{2r_2}\pfbar{3r_3}\bigr)\bigr)
    \\
    & \:-\: \kd{3r_1}\bigl(\pfbar{23}\bigl(\pfbar{3r_2}\pfbar{1r_3} - \pfbar{1r_2}\pfbar{3r_3}\bigr) \\
    & \hspace{6pc} {} -
      \pfbar{12}\bigl(-\pfbar{3r_2}\pfbar{2r_3} + \pfbar{2r_2}\pfbar{3r_3}\bigr)\bigr) \:.
  \end{align*}
  For $r_1=1$ it follows from two applications of \lemref{pf3-1},
  namely with $u_1\ldots u_k = 12r_2r_3/13r_2r_3$ and $\ell =1$, that
  $\LHS$ and $\RHS$ agree up to a sign. For $r_1=2$ one gets the same
  conclusion by applying \lemref{pf3-1} with
  $u_1\ldots u_k = 12r_2r_3/23r_2r_3$ and $\ell = 1$. For $r_1=3$ one
  gets the desired conclusion from \lemref{pf3-1} applied with
  $u_1\ldots u_k = 13r_2r_3/23r_2r_3$ and $\ell = 1$.
  
  \textbf{\em Subcase II.d.} If $3 < r_1$ holds, then \eqref{d3-1e}
  and \lemref{d3minors-even} yield
  \begin{align*}
    \LHS & \deq
           \bigl(\pfbar{1r_1r_2r_3}\pfbar{23} -
           \pfbar{2r_1r_2r_3}\pfbar{13} \\
         & \hspace{3pc} {} + \pfbar{3r_1r_2r_3}\pfbar{12}
           -\pfbar{123r_1r_2r_3}\pfT\bigr)\pfT \,.
  \end{align*}
  Expansion of the determinant along the first column yields the
  second equality in the computation below. The third equality follows
  from three applications of \lemref{pf3-1}. The fifth follows from
  two applications of \lemref{pf3-1} with
  $u_1\ldots u_k = 123r_1r_2r_3/123r_1$ and $\ell = 4$. The last
  equality follows from \lemref{42-even}.
  \begin{align*}
    &\pm\RHS \\
    & \deq \det
      \begin{pmatrix}
        \pfbar{3r_1} & \pfbar{3r_2}
        & \pfbar{3r_3} \\
        \pfbar{2r_1} & \pfbar{2r_2} & \pfbar{2r_3} \\
        \pfbar{1r_1} & \pfbar{1r_2} & \pfbar{1r_3}
      \end{pmatrix}
    \\
    & \deq \pfbar{3r_1}\bigl( \pfbar{2r_2}\pfbar{1r_3} -
      \pfbar{2r_3}\pfbar{1r_2} \bigr) \\
    & \hspace{2pc} - \pfbar{2r_1}\bigl( \pfbar{3r_2}\pfbar{1r_3} -
      \pfbar{3r_3}\pfbar{1r_2} \bigr) \\
    & \hspace{4pc} + \pfbar{1r_1}\bigl( \pfbar{3r_2}\pfbar{2r_3} - 
      \pfbar{3r_3}\pfbar{2r_2} \bigr)
    \\
    & \deq \pfbar{3r_1}\bigl(\pfbar{12r_2r_3}\pfT -
      \pfbar{12}\pfbar{r_2r_3} \bigr) \\
    & \hspace{2pc} - \pfbar{2r_1}\bigl(\pfbar{13r_2r_3}\pfT -
      \pfbar{13}\pfbar{r_2r_3} \bigr) \\
    & \hspace{4pc} + \pfbar{1r_1}\bigl(\pfbar{23r_2r_3}\pfT - 
      \pfbar{23}\pfbar{r_2r_3} \bigr)
    \\
    & \deq \bigl(\pfbar{1r_1}\pfbar{23r_2r_3} - \pfbar{2r_1}\pfbar{13r_2r_3}
      + \pfbar{3r_1}\pfbar{12r_2r_3} \bigr)\pfT \\
    & \hspace{2pc} {} + \bigl( -\pfbar{1r_1}\pfbar{23} + \pfbar{2r_1}\pfbar{13}
      - \pfbar{3r_1}\pfbar{12} \bigr)\pfbar{r_2r_3}
    \\
    & \deq \bigl(\pfbar{123r_1r_2r_3}\pfT - \pfbar{r_1r_2}\pfbar{123r_3} + \pfbar{r_1r_3}\pfbar{123r_2} \\
    & \hspace{2pc} {} - \pfbar{r_2r_3}\pfbar{123r_1}\bigr)\pfT
    \\
    & \deq \bigl(\pfbar{123r_1r_2r_3}\pfT - \pfbar{12}\pfbar{3r_1r_2r_3} + \pfbar{13}\pfbar{2r_1r_2r_3} \\
    & \hspace{2pc} {} - \pfbar{23}\pfbar{1r_1r_2r_3}\bigr)\pfT .
  \end{align*}
  Up to a sign, this is $\LHS$.
\end{prf*}


\section*{Acknowledgments}

\noindent
We thank Ela Celikbas, Luigi Ferraro, and Jai Laxmi for helpful
comments on earlier versions of this paper.

\def\soft#1{\leavevmode\setbox0=\hbox{h}\dimen7=\ht0\advance \dimen7
  by-1ex\relax\if t#1\relax\rlap{\raise.6\dimen7
  \hbox{\kern.3ex\char'47}}#1\relax\else\if T#1\relax
  \rlap{\raise.5\dimen7\hbox{\kern1.3ex\char'47}}#1\relax \else\if
  d#1\relax\rlap{\raise.5\dimen7\hbox{\kern.9ex \char'47}}#1\relax\else\if
  D#1\relax\rlap{\raise.5\dimen7 \hbox{\kern1.4ex\char'47}}#1\relax\else\if
  l#1\relax \rlap{\raise.5\dimen7\hbox{\kern.4ex\char'47}}#1\relax \else\if
  L#1\relax\rlap{\raise.5\dimen7\hbox{\kern.7ex
  \char'47}}#1\relax\else\message{accent \string\soft \space #1 not
  defined!}#1\relax\fi\fi\fi\fi\fi\fi}
  \providecommand{\MR}[1]{\mbox{\href{http://www.ams.org/mathscinet-getitem?mr=#1}{#1}}}
  \renewcommand{\MR}[1]{\mbox{\href{http://www.ams.org/mathscinet-getitem?mr=#1}{#1}}}
  \providecommand{\arxiv}[2][AC]{\mbox{\href{http://arxiv.org/abs/#2}{\sf
  arXiv:#2 [math.#1]}}} \def\cprime{$'$}
\providecommand{\bysame}{\leavevmode\hbox to3em{\hrulefill}\thinspace}
\providecommand{\MR}{\relax\ifhmode\unskip\space\fi MR }
\providecommand{\MRhref}[2]{%
  \href{http://www.ams.org/mathscinet-getitem?mr=#1}{#2}
}
\providecommand{\href}[2]{#2}

\end{document}